\documentclass{amsart} 
\usepackage{epsf} 
\usepackage{amssymb, 
latexsym, 
amsmath, 
amscd, 
array, 
graphicx}
\swapnumbers 
\numberwithin{equation}{section} 

\theoremstyle{plain} 
 
\newtheorem{thm}{Theorem}[section]

\newtheorem{prop}[thm]{Proposition}

\theoremstyle{definition} 
\newtheorem{defin}[thm]{Definition} 
 
\newtheorem{remark}[thm]{Remark}

\def\ANE{\protect\operatorname{ANE}}
\def\ANR{\protect\operatorname{ANR}}
 
\def\ch{\protect\operatorname{ch}}

\def\Fr{\protect\operatorname{Fr}}

\def\Im{\protect\operatorname{Im}}

\def\mod{\protect\operatorname{mod}} 
\def\Or{\protect\operatorname{Or}}

\def\sign{\protect\operatorname{sign}}

\def\pr{\protect\operatorname{pr}}

\def\pt{\protect\operatorname{pt}}

\def\bor{\protect\operatorname{\Omega}}
\def\borst{\protect\operatorname{\Omega}}

\def\C{{\mathbb C}} 
\def\Z{{\mathbb Z}} 
\def\Q{{\mathbb Q}} 
 
\def\N{{\mathbb N}}

\def\1{\hbox{\rm\rlap {1}\hskip.03in{\rom I}}} 
\def\Bbbone{{\rm1\mathchoice{\kern-0.25em}{\kern-0.25em} 
{\kern-0.2em}{\kern-0.2em}I}} 

\def\p{\partial}

\def\pp{\medskip{\parindent 0pt \it Proof.\ }} 
\def\wt{\widetilde} 
\def\wh{\widehat} 

\def\m{\medskip} 
\def\ov{\overline} 
\long\def\forget#1\forgotten{} % 
\def\m{\medskip}

\begin{document} 
%\date{July 10, 2005} 
%\leftline{ } 
%\centerline{ } 
\title[Graded Poisson algebras on bordism groups of garlands] 
{Graded Poisson algebras on bordism groups of garlands} 
\author[V.~Chernov]{Vladimir V. 
Chernov} 
\address{V. Chernov, Department of Mathematics, 
6188 Kemeny Hall, Dartmouth College, Hanover NH 03755-3551, 
USA} 
\email{Vladimir.Chernov@dartmouth.edu} 
%\address{Yu. Rudyak, Department of Mathematics, University 
%of Florida, 358 
%Little Hall, Gainesville, FL 32611-8105, USA} 
%\email{rudyak@math.ufl.edu} 

\begin{abstract}
Let $M^{\frak m}$ be an oriented manifold and let $\frak N$ be a set consisting of oriented closed manifolds of the same odd dimension $\frak n$. 
We consider a topological space $G_{\frak N, M}$ of commutative diagrams. Each commutative diagram consists of a few manifolds from $\frak N$ that are mapped to $M$ and of one point spaces $\pt$ that are each mapped to a pair of manifolds from $\frak N$. 

We consider the oriented bordism group $\bor_*(G_{\frak N, M})=\oplus_{i=0}^{\infty} \Omega_i(G_{\frak N, M}).$
We introduce the operations $\star$ and $[\cdot, \cdot]$ on  $\bor_*(G_{\frak N,M})\otimes \Q,$  that make $\bor_*(G_{\frak N,M})\otimes \Q$
into a $\Z$-graded $(\frak m - 2\frak n)$-Poisson algebra. 

For  $\frak N=\{S^1\}$ and a surface $M=F^2,$ the subalgebra $\bor_0(G_{\{S^1\}, F^2})\otimes \Q$ of our algebra  is related to the Andersen-Mattes-Reshetikhin Poisson algebra of chord-diagrams. 
%Our $\bor_*(G_{\{S^1\}, M})\otimes %\Q$ algebra is also related to the Chas-Sullivan string homology algebra, we plan to discuss the %relation with the Chas-Sullivan algebra in the future work.

\end{abstract}

\maketitle

\section{Basic definitions and main results} 

In this paper the word ``smooth'' means $C^{\infty}.$ All the manifolds in this paper are assumed to be smooth and oriented and $\pt$ denotes the one-point space. 
For a manifold $M$ and a compact smooth manifold $N$ we write $C^{\infty}(N, M)$ for the standard topological space of smooth maps $N\to M.$ (When $N$ is compact the strong and the weak topologies on the space of smooth maps $N\to M$ coincide.) We write $C(N, M)$ for the topological space of 
continuous maps $N\to M.$ For a set $X$ and an commutative ring $R$ we write $RX$ for the free $R$-module over the set $X.$ Given $q\in \Q$ and $s\in RX$ we write $qs$ to denote the element $s\otimes q\in RX \otimes_{\Z}\Q$. We also write $RX_{ \Q}$ instead of $RX \otimes_{\Z}\Q$.

Fix an oriented connected smooth manifold $M$ of dimension $\frak m$ and fix a set $\frak N$ of oriented connected smooth closed oriented manifolds of the same odd dimension $\frak n$.

Below we define the space $G_{\frak N, M}$ that we call {\em the space of $\frak N$-garlands in $M$.\/} Essentially it is the topological space of finite commutative diagrams that look as follows. 

Each diagram consists of $M$, a finite number of copies of $\pt$-spaces, and a finite number of manifolds $N\in \frak N$. We do allow the manifolds from $\frak N$ to participate more than once in a diagram. All the $N$-manifolds in a diagram are ordered with the enumeration starting from $1,$ and we denote by $N_i$ the manifold enumerated by $i$ in a given diagram. Each manifold $N_i$ in a diagram is continuously mapped to $M$ by exactly one map. 
Each $\pt$-space in a diagram is mapped to exactly two of $N_i$-manifolds. We do allow diagrams such that for some of the $N_i$-manifolds in them the diagram contains no map from $\pt$ to these $N_i.$

To a $\pt$-space in such a diagram we correspond its {\em index\/} $I$ which is the ordered sequence of indices $i$ of manifolds $N_i$ to which the $\pt$ is mapped. For example if $\pt$ is mapped to $N_1, N_3$ then its index is $\{1, 3\}.$ 

Each such commutative diagram $D$ gives rise to the oriented graph $\Gamma(D)$ with vertices denoted by $M$, $\pt_I,$ $N_i$ and with a manifold from $\frak N$ associated to each vertex $N_i.$
(It may be useful to think of the manifold from $\frak N$ associated to a vertex $N_i$ as a 
label of the vertex.) We require that the unoriented bipartite graph obtained by deletion of the $M$ vertex and all edges leading to it has no loops. In particular we get that all the indices of the $\pt$ vertices are distinct.

%{\it For simplicity of exposition, in this work we consider only the ``tree-like'' commutative diagrams,\/} ie the commutative diagrams as above for which the unoriented graph $\ov \Gamma(D)$ obtained from $\Gamma(D)$ by forgetting the orientation on all the edges and deleting the $M$-vertex (together with all  the edges leading to it) is a disjoint union of tree graphs. 

Note that the resulting graph $\Gamma (D)$  is rather special, for example there are no directed edges starting at $M,$ or directed edges from $N$ to $\pt.$

{\em In this paper, unless the opposite is explicitly stated, we will consider only the graphs 
that could be obtained from some commutative diagram of the type describe above. We will refer to these graphs as allowed graphs or just graphs when confusion does not arise.\/}

It is easy to see that the graph $\Gamma(D)$ corresponding 
to a diagram $D$ is completely determined by the number of $N$-vertices in it, the indices of all the $\pt$-vertices in it, and the actual manifolds from $\frak N$ associated to each one of the $N_i$-vertices in the graph.

\begin{defin}[space of garlands $G_{\frak N, M}$]
The topological space $G_{\frak N, M}$ of $\frak N$-garlands in $M$ that we are about to define is the disjoint union of topological spaces 
corresponding to different allowed graphs. 

Let $\Gamma$ be an allowed graph. Since all the $N_i$-vertices in $\Gamma$ are ordered and for each one of them there is a unique edge $N_i\to M$ in $\Gamma,$ we get that all the edges $N\to M$ in $\Gamma$ are ordered. 
Since all the indices of the $\pt$-vertices are different, all the $\pt$-vertices in $\Gamma$ are ordered by the 
lexicographical order of their indices. (With respect to this order a $\pt$-vertex with the  index $\{1,5\}$ is less than a $\pt$-vertex with the  index $\{1,6\}.$)

Since the index of a $\pt$ vertex consists of two different integers, the $\pt\to N_i$ edges in $\Gamma$ starting from the same $\pt$-vertex are also ordered by the index $i$ of the $N_i$-vertex to which they go.  
Thus all the edges $\pt\to N$ in $\Gamma$ are also ordered. For a graph $\Gamma,$ we put $\nu(\Gamma)$ to be the number of the $N$-vertices in the graph $\Gamma$ and we put $\pi(\Gamma)$ to be the number of the $\pt\to N$ edges in $\Gamma.$ For $j\in \{1, \cdots, \pi(\Gamma)\}$ we denote by $i(j)\in \{1, \cdots ,\nu(\Gamma)\}$ the index of the $N$-vertex to which the $j$-th edge $\pt\to N$ is going.

Consider the space of all (not necessarily commutative) diagrams
obtained from $\Gamma$ by equipping all the directed edges $\pt\to N_i$ in $\Gamma$ by some maps from $\pt$ to the manifold $N_i\in \frak N$ associated to the vertex $N_i,$ and by equipping all the directed edges $N_i\to M$ in $\Gamma$ with some  continuous maps into $M$  of the manifold $N_i\in \frak N$ associated to the vertex $N_i.$ Each map $\pt\to N$ is completely described by its image point in $N\in \frak N.$ Thus the space of all such diagrams corresponding to $\Gamma$ can be identified with 
\begin{equation}\label{hugeproduct}\bigl(\prod_{j=1}^{\pi(\Gamma)}N_{i(j)}\bigr)\times \bigl(\prod_{i=1}^{\nu(\Gamma)} C(N_i, M)\bigr).
\end{equation}(Recall that by the definition of the allowed graphs $\Gamma,$ $\nu(\Gamma)$  coincides with the number of $N\to M$ edges in $\Gamma.$) This identification gives topology on the space of all such not-necessarily commutative diagrams.

{\it Put $G_{\frak N, M}^{\Gamma}$ to be the subspace of $\bigl(\prod_{j=1}^{\pi(\Gamma)}N_{i(j)}\bigr)\times \bigl(\prod_{i=1}^{\nu(\Gamma)} C(N_i, M)\bigr)$ that is formed by the diagrams that are actually commutative. We equip $G_{\frak N,M}^{\Gamma}$ with the induced topology.

The space $G_{\frak N, M}$ is the disjoint union of the subspaces $G_{\frak N, M}^{\Gamma}$ over all the allowed graphs $\Gamma,$ i.e. $G_{\frak N, M}=\sqcup_{\Gamma}G_{\frak N, M}^{\Gamma}.$\/}
\end{defin}

\begin{remark}
If the graph $\Gamma$ has no $\pt$-vertices, then $G_{\frak N, M}^{\Gamma}$ is a product of spaces $C^{\infty}(N, M).$

If $\{\frak N\}=\{\pt\}$, then $G_{\frak N, M}^{\Gamma}$ are related to completions of the configuration space. For example in the case where $\Gamma$ is the graph with one $\pt$-vertex and two $N$-vertices (that are all equal to  $\pt$ as well) the resulting topological space is the stratum of codimension one in the configuration space of two points in $M.$ The spaces resulting from more complicated $\Gamma$ can be interpreted similarly.
\end{remark}

\begin{figure}[htbp] 
\begin{center} 
\epsfxsize 8 cm 
\includegraphics[width=8cm]{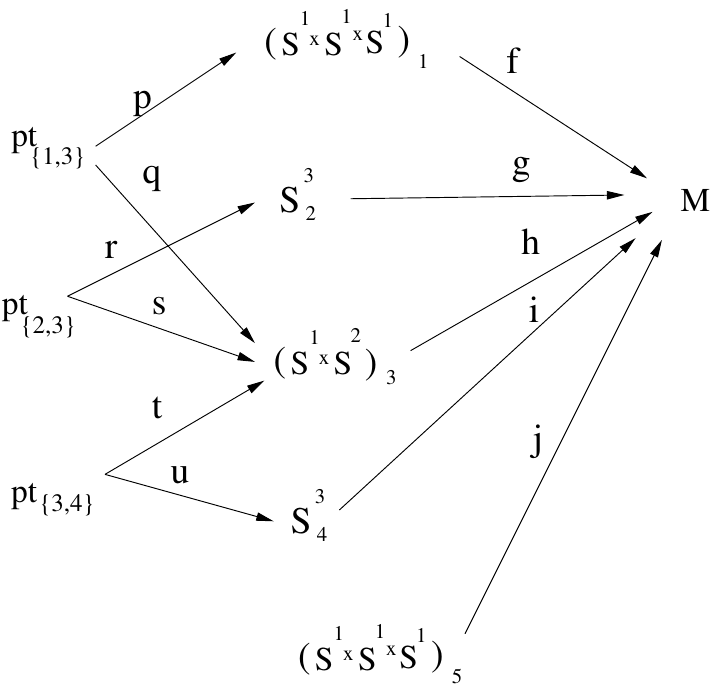} 
%\hepsffile{garlandtreediag.eps} 
\end{center} 
\caption{Example of a garland for $\frak N=\{S^1\times S^1\times S^1,S^3, S^1\times  S^2\}$}\label{garlandtreediag.fig} 
\end{figure}

\begin{defin}[Oriented bordism group $\borst_*(Y)$]
For a topological space $Y$, we denote by $\borst_n(Y)$ the $n$-dimensional 
oriented bordism group of $Y$. Recall that $\borst_n(Y)$ is 
the set of 
equivalence classes of continuous maps $g: V^n \to X$ 
where $V$ is a 
closed oriented smooth $n$-manifold. Here two maps $g_1: V_1^n \to Y$ 
and $g_2: V_2^n 
\to Y$ are equivalent if there exists a pair $(Z, h)$ where $Z$ is a compact oriented smooth $(n+1)$-manifold whose oriented boundary $\p 
Z$ is 
diffeomorphic to $V_1\sqcup (-V_2)$ and $h:Z\to Y$ is a continuous map such that 
$h|_{\p Z}=f_1\sqcup f_2.$ (Here we use a slightly ambiguous notation and denote by $f_2$ both the maps from $V_2$ and from $-V_2$ to $Y.$)
The 
operation of disjoint union gives $\borst_n(Y)$ the structure of a
commutative group, see~\cite{Rudyak, Stong, Switzer} for details. 
For a manifold $V$ and a map $g:V^n\to Y$ we will denote by $[V,g]$ the corresponding 
bordism class in $\borst_n(Y).$

A $0$-dimensional manifold is a disjoint union of finitely many one-point-spaces, 
and an orientation of a zero-dimensional manifolds is an assignment of a sign $\pm 1$ to each of the one-point-subspaces of it. We will write $\pt^+$ and $\pt^-$ respectively for the positively and for the negatively oriented space $\pt.$

For a space $Y$, the group $\borst_0(Y)=H_0(Y)=\Z \pi_0(Y)$ is the 
free commutative 
group with the basis $\pi_0(Y)$. So, every element of 
$\borst_0 
(Y)$ can be 
represented as a finite linear combination $\sum \gamma_kP_k$ 
with $\gamma_k\in 
\Z$ and $P_k\in Y$. Conversely every such linear combination gives 
us an element 
of $\borst_0(Y)$.

Let $[V]\in H_n(V)$ be the fundamental class of a closed 
oriented 
$n$-dimensional manifold $V$. Every map $f: V \to X$ gives 
us an 
element $f_*[V]\in H_n(X)$, and the correspondence $(V,f) 
\mapsto 
f_*[V]$ yields the Steenrod--Thom homomorphism 
$$ 
\tau: \borst_n(X) \to H_n(X). 
$$ 
This homomorphism is an isomorphism 
for $n\le 3$ 
and an epimorphism for $n\le 6$. The bordism group $\borst_*(\pt)$ has the natural ring structure and $\borst_*(Y)\otimes _{\borst_*(\pt)}\Q=H_*(Y,\Q).$ See~\cite{Thom} and 
\cite[Theorem 7.37]{Rudyak}, \cite[Theorem p. 209]{Stong} 
for modern proofs. 
\end{defin}

For $\eta\in \bor_k(G_{\frak N, M})\otimes \Q=\bor_k(G_{\frak N, M})_{\Q}$ we put $|\eta|=k$ and call it the {\it degree of the element $\eta$.\/} We will often write $(-1)^{|\eta|}$ instead of $(-1)^k.$ Since $\bor_*(G_{\frak N, M})\otimes \Q=\oplus_{k=0}^{\infty} (\bor_k(G_{\frak N, M})\otimes \Q)$ we see that $\bor_*(G_{\frak N, M})_{ \Q}=\oplus_{k=0}^{\infty} (\bor_k(G_{\frak N, M})_{\Q})$ is a $\Z$-graded $\Q$-vector space.

Below we list some of the main results of this work. 
We construct a bilinear operation $\star$ on $\bor_*(G_{\frak N, M})_\Q.$

\begin{thm}\label{main0}
Operation $\star$ makes $\bor_*(G_{\frak N, M})_\Q$ into an associative, graded commutative algebra. Namely for all $\frak N, \eta_i\in \bor_{j_i}(G_{\frak N, M})_\Q, i=1,2,3,$ we have:
\begin{description}
\item[1] $\eta_1\star \eta_2=(-1)^{|\eta_1||\eta_2|}\eta_2\star \eta_1;$ 
\item[2] $(\eta_1\star \eta_2)\star \eta_3=\eta_1\star (\eta_2\star \eta_3).$
\item[3] The $\star$-product is of degree zero, i.e. for all $j_1, j_2\in \N$  we have $$\star: (\bor_{j_1}(G_{\frak N, M})_\Q)\otimes (\bor_{j_2}(G_{\frak N, M})_ \Q)\to \bor_{j_1+j_2}(G_{\frak N, M})_\Q.$$
\end{description}
\end{thm}
See Theorems~\ref{multiplicationdimension} and Theorem~\ref{associativecommutative}.

Recall now the definition of the graded Poisson algebras, see the work of Cattaneo, Fiorenza and Longoni for the good exposition~\cite{CFL} of the subject.

\begin{defin}
A {\em graded Poisson algebra of degree $n$, or $n$-Poisson algebra\/} is a triple $(A, \cdot, \{,\})$ consisting of $\Z$-graded vector space $A$ endowed with a degree zero graded commutative product and degree $-n$ Lie bracket. The bracket is required to be a biderivation of the product. Namely $\{a, b\cdot c\}=\{a, b\}\cdot c+(-1)^{|b|(|a|+n)}b\cdot \{a, c\}.$ The usual Poisson algebras appear as a particular case of the $0$-Poisson algebras, and when $n=1$ one gets the famous Gestenhaber algebras or Schouten algebras.
\end{defin}

We construct a bilinear operation $[\cdot, \cdot]$ on $\bor_*(G_{\frak N, M})_\Q.$

The operations $\{\cdot, \cdot\}$ and $[\cdot, \cdot]$ are related as follows. For elements $\eta_i\in \bor_{j_i}(G_{\frak N, M})_\Q,$ $ i=1,2,$ we put $\{\eta_1, \eta_2\}=(-1)^{|\eta_2|(\frak m+1)}[\eta_1, \eta_2].$ We extend $\{\cdot, \cdot\}$ to the whole $\bor_*(G_{\frak N, M})_\Q$ by linearity.

\begin{thm}\label{main1}
The operations $\{\cdot, \cdot\}$ and $\star$ make 
$\bor_*(G_{\frak N, M})_\Q$ into a $\Z$-graded $(\frak m - 2\frak n)$-Poisson algebra.
Namely,  $\eta_i\in \bor_{j_i}(G_{\frak N, M})_\Q, i=1,2,3,$ we have the following identities and properties:
\begin{description}
\item[1] graded skew-symmetricity $\{\eta_1, \eta_2\}=-(-1)^{(|\eta_1|-\frak m)(|\eta_2|-\frak m)}\{\eta_2, \eta_1\};$ 
\item[2] graded Jacobi identity
$$
\{\eta_1, \{\eta_2, \eta_3\}\}=\{\{\eta_1, \eta_2\}, \eta_3\}+(-1)^{(|\eta_1|-\frak m)(|\eta_2|-\frak m)}\{\eta_2, \{\eta_1, \eta_3\}\}.
$$
\item[3] graded Leibnitz identity 
$$
\{\eta_1, \eta_2\star\eta_3\}=\{\eta_1, \eta_2\}\star \eta_3+(-1)^{|\eta_2|(|\eta_1|-\frak m)}\eta_2\star\{\eta_1, \eta_3\}.
$$
\item[4] The $\{\cdot, \cdot\}$ operation is of degree $2\frak n-\frak m$, i.e. for all $j_1, j_2\in \N$ we have 
$$
\{\cdot, \cdot\}:(\bor_{j_1}(G_{\frak N, M})_\Q)\otimes (\bor_{j_2}(G_{\frak N, M})_\Q)\to \bor_{j_1+j_2+2\frak n-\frak m}(G_{\frak N, M})_\Q.
$$
\end{description}
\end{thm}

Theorem~\ref{main1} is a straightforward corollary of Theorems ~\ref{dimension},~\ref{skewsymmetry},~\ref{Jacobi}~\ref{Poisson} after one reformulates them in terms of the $\{\cdot, \cdot\}$ operation.

\begin{remark}[possible degree shifts and similarities with the Chas-Sullivan loop homology Gerstenhaber algebra]
It may be useful to change the definition of the degree so that for for $\eta\in \bor_j(G_{\frak N, M})\otimes \Q$ we have $\deg(\eta)=|\eta|+2\frak n-\frak m.$
In this case $(-1)^{|\eta|-\frak m}$ in the identities of Theorem~\ref{main1} becomes $(-1)^{\deg (\eta)},$ and the Lie bracket 
$\{\cdot, \cdot\}$ becomes an operation of degree zero. However, then $\star$ becomes an operation of degree $\frak m-2\frak n$ rather than of degree zero and instead of $\eta_1\star\eta_2=(-1)^{|\eta_1||\eta_2|}\eta_2\star \eta_1$ in Theorem~\ref{main0} we would get 
$\eta_1\star \eta_2=(-1)^{(\deg (\eta_1)+\frak m)(\deg (\eta_2)+\frak m)}\eta_2\star \eta_1.$  For these reasons we prefer not to do this degree shift.

Note that if $\frak m$ is odd, then the identities of Theorem~\ref{main1} are similar to the identities of the loop homology Gerstenhaber algebra of Chas-Sullivan~\cite[Theorem 4.7]{ChasSullivan1}. 
However our and Chas-Sullivan Lie brackets have different degrees. Our Lie bracket is of degree $2\frak n-\frak m,$ while Chas-Sullivan Lie bracket is of degree $1-\frak m$ (or of degree $1$ after the appropriate shift of the degrees). Our product is of degree $0,$ while the Chas-Sullivan product is of degree $-\frak m$  (or of degree $0$ after the appropriate shift of the degrees). For these reasons we believe that there is no immediate relation between our and the Chas-Sullivan algebras.

There is a clear relation between the Chas-Sullivan bracket and our bracket, but it exists only in the case where you pair bordism classes for the very simple graphs $\Gamma$ consisting of only one $\mathfrak N$-manifold which is $S^1$ and the arrow connecting it to the target manifold $M.$ In this case our bracket results in a bordism class corresponding to a graph with one $\pt$-vertex connected by arrows to two copies of $S^1$ which are then connected by arrows to $M.$ One interprets such a commutative diagram as a loop product based at the image of the $\pt$-vertex in $M$ and gets a bordism class in the space $C^{\infty}(S^1, M).$ This is essentially the main idea of the Chas-Sullivan bracket. However this analogy does not extend to more complicated graphs $\Gamma.$

It is also useful to note that the original Chas-Sullivan argument had significant gaps that used transversality of chains of maps of circles and this gap was unfortunately never fixed. Of course the algebra of Chas and Sullivan is indeed well defined but the justification requires complicated homotopy theory techniques developed by Cohen-Jones~\cite{CohenJones} and Cohen-Voronov~\cite{CohenVoronov}. One may also refer to the book of Kreck~\cite{Kreck} that discusses the new mathematical methods developed in part to deal with similar transversality issues.

Bordism theory that we use here to define our algebra is natural to use in the transversality arguments but as one can see from the Appendix~\ref{transversality} overcoming the transversality issues is still very difficult.
\end{remark}

\begin{remark}[Relation to the previously known algebras]
As it is discussed in Section~\ref{SectionLieLoops}, for $\frak N=\{S^1\}$ and an oriented  surface $F^2,$ the subalgebra $\bor_0(G_{\frak \{S^1\}, F^2})_\Q$ of our algebra $\bor_*(G_{\frak \{S^1\}, F^2})_\Q$ formed by 0-dimensional bordisms
is closely related to the Andersen-Mattes-Reshetikhin~\cite{AMR1}, \cite{AMR2} Poisson algebra of chord diagrams. 
%and to the Goldman-Turaev \cite{Goldman}, \cite{Turaevskein} algebra of loops on a surface. 
\end{remark}

%\begin{remark}[Getting rid of the enumeration of the $N$-vertices]

%In section~\ref{noenumeration} we introduce a possible version of the space $\wt G_{\frak N, M}$ %and define the operations $\star$ and $[\cdot, \cdot]$ on $\bor_*(\wt G_{\frak N, M})\otimes \Q.$
%However, we believe that there are better topological spaces of commutative diagrams with the %$N$-manifolds in them not enumerated.  
%\end{remark}

\begin{remark}[Relation to the algebra in our joint work with Yu. Rudyak~\cite{ChernovRudyakGarlands}]
The Lie bracket and many other operations on $\bor_*(G_{\frak N, M})$ for one element sets $\frak N$ 
were first introduced in our joint preprint with Yuli Rudyak~\cite{ChernovRudyakGarlands}. In that preprint  we were able to get only the $\mod 2$ proof of the Jacobi identity. This was due to the major difficulties we had in computing orientations of the manifolds parameterizing the bordism classes of the terms in the Jacobi identity. This work shows that for $\frak N$ consisting of  manifolds of the same odd dimension $\frak n$, the operation $[\cdot, \cdot]$ is indeed a graded Lie bracket.

The main ingredient of $[\cdot, \cdot]$, first appeared (under the name $\mu$) in our work with Rudyak~\cite{ChernovRudyakGT} where we used it to generalize the linking number invariant to the case of arbitrary nonzero homologous linked submanifolds $\phi_1(N_1^{n_1}), \phi_2(N_2^{n_2})$ of $M^{n_1+n_2+1}.$ We applied these generalzied linking numbers to the study of causality in spacetimes~\cite{ChernovRudyakCausality}.

The idea of the $\mu$ operation is: that given two maps $f_i:V^i\to (C^{\infty}(N_i, M))$ one can consider the pull back $W$ of the adjoint maps $V_i\times N\to M$. This pull back is naturally mapped to the space of commutative diagrams which look as $\pt$ space mapped to $N_1, N_2$, which are both mapped to $M.$

Essentially the Lie bracket constructed in the current paper is the (symmetrized) iterated version of the $\mu$-operation.

%Note also that for odd-dimensional $M,$ the signs in statements $1$ of the two %Theorems above are different from those expected by the analogy with the %intersection pairing on $H_*(M).$ These signs result from our orientation %conventions, see Remark~\ref{strangesign}. This may indicate the existence of one %more graded Lie bracket on $\bor_*(G_{\frak N, M})\otimes \Q$ that behaves %similarly to the 
%intersection pairing under the permutation of the two terms.
\end{remark}

\begin{remark}[More general allowed graphs]
It is rather easy to construct the topological space of $\frak N$-garlands in $M$ corresponding to more general allowed graphs, for which a $\pt$-vertex can  be connected to the same $N$-vertex by many edges, and for which the unoriented graph $\ov \Gamma$ corresponding to  a graph $\Gamma$ is not necessarily a disjoint union of trees. This can be done by enumerating all the $\pt\to N$ edges in a more general allowed graph. The $\star$ operation can be easily generalized to tensor product of $\Q$ and the oriented bordism group of such more general commutative diagrams. 

However we can not currently generalize the transversality result of Appendix~\ref{transversality} to the more general commutative diagrams. If this is possible, then $[\cdot, \cdot]$ also can be generalized to the tensor product of $\Q$ and the bordism group of the space of such more general commutative diagrams. 
\end{remark}

\begin{remark}[Counting double points] In the simplest case where $[\cdot, \cdot]$ pairs two $0$-dimensional bordism classes of two loops on a surface (considered as trivial commutative diagrams with one copy of $S^1$ and no $\pt$-vertexes) the sum of the absolute values of the coefficients of the result gives a lower estimate on the number of double points of the two loops. Since this sum does not change under the homotopy of the loops, we get that it estimates from below the minimal number of double points of two loops homotopic to the given ones. In fact, as we proved with Cahn~\cite{CahnChernov}, this sum always equals to the minimal number of such double points. (Hence in particular the Andersen-Mattes-Resheikhin bracket always gives the exact answer to such a problem.)

One should keep in mind that the estimate coming from a related Goldman Lie bracket~\cite{Goldman} is, as it was shown by Chas~\cite{Chas}, generally not equal to the minimal number of double points. 

For higher dimensional situations our operation $[\cdot, \cdot]$ also provides interesting estmates on the number of double points (and more generally on the algebraic topology complexity of the set parameterizing the double points). 
\end{remark}

\begin{remark}[Reasons for enumerating the garland pieces and possible ways
to avoid this]
We used enumeration of the $N$-vertices in the graphs to define topology on $G_{\frak N, M}^{\Gamma}.$ However if $\Gamma'$ is obtained from $\Gamma$ by reenumerating the vertices, then the corresponding topological spaces of all {\it noncommutative\/} diagrams corresponding to $\Gamma$ and $\Gamma'$ are homeomorphic. These spaces are products of mapping spaces described by equation~\eqref{hugeproduct} and the homeomorphism exchanges the factors in the products. This homeomorphism induces the homeomorphism between the subspaces $G_{\frak N, M}^{\Gamma}$ and $G_{\frak N, M}^{\Gamma'}$ of commutative diagrams corresponding to the graphs $\Gamma$ and $\Gamma'$. Thus the topology on the space of commutative diagrams corresponding to a graph $\Gamma$ actually does not depend on the enumeration of the $N$-vertices in it.

The reason why we can not currently bypass enumerating the vertices is that in the construction of algebra operations and in the verification of the identities they satisfy we need to be able to easily identify the vertices we refer to. This is easy if $N$-vertices are enumerated. However an unfortunate side effect of this enumeration is that the algebra is defined on $\bor_*(G_{\frak N, M})\_\Q$ rather than on $\bor_*(G_{\frak N, M}).$ 

If one manages to get rid of the enumeration of the $N$-vertices in the graphs, then the algebra will be defined on $\bor_*(G_{\frak N, M})$ rather than on $\bor_*(G_{\frak N, M})\otimes \Q.$ In the cases when $\frak N$ consists of manifolds of the same even or odd dimension the algebra operations would satisfy identities similar to those of Theorem~\ref{main1}. Moreover a similar construction should give a Poisson algebra on the unoriented bordism group.
\end{remark}

\begin{remark}[Some possible applications for the garland algebras and relations to Vassiliev invariants]
The Lie bracket between $0$- and $1$-dimensional bordism classes for the subspaces corresponding to graphs with no $\pt$-vertices allows one to construct generalizations of linking number to the case of arbitrary submanifolds $M_1^{m_1}, M_2^{m_2}\subset N^{n}$ in the appropriate dimensions $m_1+m_2+1=n,$ see~\cite{ChernovRudyakGT}. These generalized ``linking numbers'' were successfully applied by Rudyak and the author to the study of causality in spacetimes and admit a nice physics interpretation as the algebraic number of times an observer sees light from an event when traveling along a trajectory with the speed less or equal than the speed of light, see~\cite{ChernovRudyakCausality}.

In general as it was shown by Nemirovski and the author the Legendrian (as opposed to the Topological) linking allows to completely reconstruct causal relation for all reasonable spacetimes~\cite{ChernovNemirovski1, ChernovNemirovski2, ChernovCausality} thus proving a few conjectures on the subject inspired by the former students of Penrose (Low~\cite{Low0, Low1, Low2}, Natario and Tod~\cite{NatarioTod}). 

However there is still a lot of new results that can be obtained along the approach developed in our works with Rudyak. In particular the gardland algebra approach should allow to construct secondary invariants defined when the generalized (affine) linking numbers are zero and to estimate the number of times one sees the same star on the night sky due to gravitational lensing. 

Also the garland approach allows one to construct invariants of linked submanifolds under link homotopy without the dimension restriction that $n_1+n_2+1=m$ and it would be interesting to relate the resulting invariants to those introduced by Koschorke~\cite{Koschorke1, Koschorke2}.

The generalized linking numbers are Vassiliev~\cite{Vassiliev} invariants of degree one.
The garland algebras for more complicated graphs when there are multiple $\pt$-vertices should allow one to introduce and study Vassieliev invariants of higher degree for linked submanifolds in all manifolds. They most likely also have nice applications to the study of causality in spacetimes.
\end{remark}

\m {\bf The paper is organized as follows.} In Section~\ref{sectionAoperation} we introduce the $A$-operation which is the main ingredient of the Lie bracket $\{\cdot, \cdot\}.$ In Section~\ref{bracket} we construct the Lie bracket and in Section~\ref{ProofJacobi} we show that it satisfies the Jacobi identity. In Section~\ref{SectionCoperation} we introduce the $C$-operation which is the main ingredient of the $\star$-product defined in Section~\ref{multiplication}. The main ingredient of the Leibnitz identity is proved in Section~\ref{proofPoissoningredient}, while the identity itself is proved in Section~\ref{differentialalgebras}. In Section~\ref{SectionLieLoops} we relate our algebra to the Andersen-Mattes-Reshetikhin algebra.

In Appendix~\ref{transversality} we prove the transversality result needed to show that the operations are well defined.  In Appendix~\ref{EvenDimensional} we mention the algebras that will appear if the manifolds in $\frak N$ are of the same even dimension. In Appendices~\ref{proofskewsymmetryingredient},~\ref{proofJacobiingredient}
%,~\ref{proofPoissoningredient} 
we prove the main technical orientation related identities  that would imply respectively the graded skew symmetry, graded Jacobi 
%and graded Poisson 
identities for our bracket operation. Finally, in Appendix~\ref{coincedenceLiebracket} we outline the construction and properties of one more Lie bracket, which is well defined only if $\frak N$ is a one lement set. This Lie bracket is related to coincidence theory, since the pull back used to define is intersected with the $N\times N$-diagonal  and hence takes into account only the double points of the two maps whose two preimages are equal.

\section{Operation $A$ on $\bor_*(G_{\frak N, M})$ }\label{sectionAoperation}
In this section we introduce the operation $A$ on $\bor_*(G_{\frak N, M}).$ In Section~\ref{bracket} we will use the symmetrization of $A$ to construct the graded Lie algebra structure on $\bor_*(G_{\frak N, M})_\Q.$

\begin{defin}[operation $B_{k_1,k_2}$ and actions of permutation groups on graphs]\label{Boperation}
Let $\Gamma_1$ and $\Gamma_2$ be two allowed graphs. Let $\nu(\Gamma_i)$ be the number of $N$-vertices in the graph $\Gamma_i, i=1,2.$ For positive integer $k_i\leq \nu(\Gamma_i), i=1,2,$ we define the allowed graph $B_{k_1, k_2}(\Gamma_1, \Gamma_2)$ to be the graph resulting after the following sequence of operations: 
\begin{description}
\item[1] Take the disjoint union of the graphs $\Gamma_1$ and $\Gamma_2,$ preserving the manifolds in $\frak N$ that were associated to the $N$-vertices in the graphs $\Gamma_1$ and $\Gamma_2.$
\item[2] Keep the indices of the $N$-vertices that came from $\Gamma_1$ and increase by $\nu(\Gamma_1)$ the indices of the $N$-vertices that came from $\Gamma_2.$ 
\item[3] Change the indices of the $\pt$-vertices in $\Gamma_2$ so that the new index of each $\pt$-vertex gives the indices (with respect to the shifted enumeration of $N$-vertices) of the $N$-vertices in $\Gamma_2$ connected to this $\pt$-vertex by an oriented edge.
\item[4] Identify the two $M$-vertices in $\Gamma_1\sqcup \Gamma_2$ to get just one $M$ vertex and redirect 
all the $N\to M$ edges in the two graphs to it.
\item[5]Finally add one new $\pt$-vertex that has index $\{k_1, k_2+\nu(\Gamma_1)\}$ and hence is connected by the oriented edges to the $N_{k_1}$- and $N_{k_2+\nu (\Gamma_1)}$-vertices in the resulting graph. Note that the $N_{k_2+\nu(\Gamma_1)}$-vertex in the resulting graph corresponds to the $N_{k_2}$-vertex in the graph $\Gamma_2.$
\end{description}

For an allowed graph $\Gamma$ and a permutation $\alpha\in S_{n}$ of the ordered sequence of numbers $\{1, 2, \cdots, n\},$ {\it we define the allowed graph\/} $\alpha \cdot \Gamma$ as follows: 
\begin{description}
\item[1] If $n\neq \nu(\Gamma),$ then $\alpha\cdot \Gamma=\Gamma.$
\item[2] If $n= \nu(\Gamma),$ then to get $\alpha \cdot \Gamma$ keep the manifolds from $\frak N$ associated to the $N$-vertices unchanged and keep all the edges in the graph $\Gamma$ unchanged.
Change the indices of the $N$-vertices as it is described by the permutation $\alpha.$ 
Change the  indices of the $\pt$-vertices in the graph, so that an index of a $\pt$-vertex is the ordered sequences of the new indices of the $N$-vertices to which the $\pt$-vertex is connected by an edge. 
\end{description}

For $n_1, n_2\in \N$ {\em define a permutation\/} $(n_1, n_2)\in S_{n_1+n_2}$ by putting its value on $\{1, 2, \cdots, n_1+n_2\}$ to be $\{n_1+1, n_1+2, \cdots, n_1+n_2, 1, 2, \cdots, n_1\}.$ 
\end{defin}

The following Proposition follows immediately from the definition of $B_{k_1, k_2}.$

\begin{prop}\label{propgraph}
Let $\Gamma_1, \Gamma_2,\Gamma_3$ be three allowed graphs, then 
\begin{description}
\item[1] $B_{\nu (\Gamma_1)+\ov k_2, k_3}\bigl (B_{k_1, k_2}(\Gamma_1, \Gamma_2), \Gamma_3\bigr)=  \alpha_1 \cdot \bigl( B_{k_2, k_1} \bigl (B_{\ov k_2, k_3} (\Gamma_2, \Gamma_3), \Gamma_1\bigr)\bigr),$ for all positive integer $k_1\leq \nu (\Gamma_1), k_2, \ov k_2\leq \nu (\Gamma_2), k_3\leq \nu (\Gamma_3),$ and for the permutation \\ $\alpha_1=\bigl (\nu(\Gamma_1), \nu(\Gamma_2)+\nu(\Gamma_3) \bigr)\in  S_{(\nu(\Gamma_1)+\nu(\Gamma_2)+\nu(\Gamma_3))}$ (defined in~\ref{Boperation}). 
\item[2] $B_{k_1, k_2}(\Gamma_1, \Gamma_2)=
\alpha_2\cdot B_{k_2, k_1}(\Gamma_2, \Gamma_1),$ for all positive integer $k_1\leq  \nu (\Gamma_1), k_2\leq  \nu (\Gamma_2)$ and for the permutation $\alpha_2=\bigl( \nu(\Gamma_1), \nu(\Gamma_2)\bigr)\in S_{(\nu(\Gamma_1)+\nu(\Gamma_2))}$ (defined in~\ref{Boperation}).
\end{description}
\qed
\end{prop}

\begin{defin}[Nice maps]
Let us define {\it nice maps} of a smooth manifold $V$ into $G_{\frak N, M}.$ We need them  to simplify the transversality arguments in this work.

Let $V$ be a connected smooth manifold possibly with boundary, and let $f:V\to G_{\frak N, M}$ be a continuous map. Since $V$ is connected, $f$ maps it to some subspace $G_{\frak N, M}^{\Gamma}.$
Thus for every $v\in V,$ $f(v)$ is a commutative diagram from $G_{\frak N, M}^{\Gamma}.$
We say that $f$ is {\it nice\/} if the following conditions hold:
\begin{description}
\item[a]
Choose a $\pt\to N$ edge in $\Gamma.$ Let $N\in \frak N$ be the manifold corresponding to the vertex $N.$ The map $f:V\to G_{\frak N, M}^{\Gamma}$ defines the adjoint map $V\times \pt=V\to N$ that sends $v\in V$ to the image point of the $\pt \to N$ map in the commutative diagram $f(v).$ We require that this map $V\to N$ is smooth for all edges $\pt\to N.$
\item[b] Choose an $N$-vertex in $\Gamma$ and let $N\in \frak N$ be the manifold corresponding to this vertex. The map $f:V\to G_{\frak N, M}^{\Gamma}$ gives rise to the adjoint map $V\times N\to M$ defined by the condition that its restriction to each $v\times N, v\in V,$ is the map $N\to M$ from the commutative diagram $f(v).$ We require that this map $V\times N\to M$ is smooth, for every vertex $N$ in $\Gamma.$
\end{description}

We say that a continuous map $f:V\to G_{\frak N, M}$ of a not necessarily connected smooth manifold is {\it nice,\/} if its restriction to every connected component of $V$ is nice.

Every $f:V\to G_{\frak N, M}$ is homotopic to a nice map. Moreover if $f|_{\partial V}$ is nice then $f$ is homotopic relative $\partial V$ to a nice map, see Theorem~\ref{transversalitytheorem} and Remark~\ref{nicebordismmaps}.

\end{defin}

\begin{defin}[operation $A_{k_1, k_2}^{j_1, j_2}$]\label{Aoperation}
For bordism classes $\omega_i\in \bor_{j_i}(G_{\frak N, M}^{\Gamma_i}), i=1,2,$ 
choose nice maps $f_i:V_i^{j_i}\to G_{\frak N, M}^{\Gamma_i}$   realizing 
$\omega_i.$  
Thus each $f_i(v_i), v_i\in V_i,$ is a commutative diagram from $G_{\frak N, M}^{\Gamma_i}, i=1,2.$ 

Fix $k_1\in \{1, \cdots, \nu(\Gamma_1)\}$ and $k_2\in \{1, \cdots, \nu(\Gamma_2)\}.$ Everywhere below we denote by $N_{k_i}$ the manifold from $\frak N$ corresponding to the 
$N_{k_i}$-vertex in the graph $\Gamma_i, i=1,2.$

For $v_i\in V_i, i=1,2,$ let $f_{i, k_i}(v_i): N_{k_i}\to M$ be the map from the commutative diagram $f_i(v_i).$ Let $\ov f_{i, k_i}:V_i\times N_{k_i}\to M, i=1,2,$ be the smooth adjoint maps given by $\ov f_{i, k_i}|_{v_i \times N_{k_i}}=f_{i, k_i}(v_i), v_i\in V_i.$ 

By statement 3 of Theorem~\ref{transversalitytheorem}, we can deform $f_2$ by a homotopy so that $\ov f_{1, k_1}$ and $\ov f_{2, k_2}$ become transverse.  
Then  the map $\ov f_{1, k_1}\times \ov f_{2, k_2}:(V_1\times N_{k_1})\times (V_2\times N_{k_2})\to M\times M$ is transverse to the diagonal submanifold $\Delta=\Delta_{M\times M}\subset M\times M$
defined as $\Delta=\{(m,m)\in M\times M| m\in M\}.$ 

{\it For simplicity of exposition we equip all the manifolds in the construction with Riemannian metrics. Thus the normal bundle $\nu(W)$ of a submanifold $W$ is identified with the orthogonal compliment bundle $(TW)^{\perp}$ of $TW.$ The products of oriented Riemannian manifolds will be equipped with the product Riemannian metric and with the product orientation. The product $\pt\times Y$ of an oriented one-point-space $\pt$ and a manifold $Y$ is oriented 
so that the diffeomorphism $\pt\times Y\to Y, \pt\times y\to y,$ is orientation preserving (respectively reversing) if the orientation of $\pt$ is $+1$ (respectively $-1$). Products of an oriented $0$-dimensional manifold and an oriented manifold $Y$ are oriented similarly.\/}
  
In particular, the normal bundle $\nu \Delta$  to $\Delta$ in $M\times M$ is identified with the 
$\frak m$-dimensional bundle $(T\Delta)^{\perp}$ orthogonal in $T(M\times M)$ to the tangent bundle $T\Delta.$ We identify $\Delta$ with $M$ and {\it we orient the normal bundle $\nu \Delta=(TM)^{\perp}$ in such a way that a positive orientation frame of $T\Delta=TM$ followed by a positive orientation frame of $\nu \Delta=(TM)^{\perp}$ gives a positive orientation of $T(M\times M)|_{\Delta}.$\/}  We denote by $\pr$ the orthogonal projection operator $\pr:T(M\times M)|_{\Delta}=T\Delta\oplus \nu \Delta=T\Delta\oplus (T\Delta)^{\perp}\to \nu \Delta=(T\Delta)^{\perp}.$

Consider the $(j_1+j_2+2\frak n-\frak m)$-dimensional pullback submanifold $W\subset (V_1\times N_{k_1}\times V_2\times N_{k_2})$ that consists of all $(v_1, n_{k_1}, v_2, n_{k_2})$ such that $\ov f_{1, k_1}(v_1, n_{k_1})\times \ov f_{2, k_2}(v_2, n_{k_2})\in \Delta.$ We have the natural bundle map $\pr\circ (\ov f_{1, k_1}\times \ov f_{2, k_2})_*|_{(TW)^{\perp}}:TW^{\perp}=\nu W\to (TM)^{\perp}=\nu \Delta$ and {\it we orient $\nu W=(TW)^{\perp}$ so that this bundle map is orientation preserving.  We orient the manifold $W$ so that for $w\in W$ a positive orientation frame of $T_wW$ followed by a positive orientation frame of $\nu_w W=(T_wW)^{\perp}$ gives a positive orientation frame of  $T_w(V_1\times N_{1, k_1}\times V_2\times N_{2, k_2}).$

If $W$ is $0$-dimensional, then we orient its components as follows: the orientation of $w\in W$ is $+1$ if a positive orientation frame of $\nu_w W=(T_wW)^{\perp}$ is a positive orientation frame of $T_w(V_1\times N_{1, k_1}\times V_2\times N_{2, k_2});$ and the orientation of $w$ is $-1$ otherwise.\/}

Each point $w=(v_1, n_{1, k_1}, v_2, n_{2, k_2})\in W$ gives rise to a commutative diagram $g(w)\in G_{\frak N, M}^{B_{k_1, k_2}(\Gamma_1, \Gamma_2)}$ as follows. Recall that the graph $B_{k_1, k_2}(\Gamma_1, \Gamma_2),$ see Definition~\ref{Boperation}, consists of three parts: the part coming from $\Gamma_1,$ the part coming from $\Gamma_2,$ and the extra $\pt$-vertex that is connected by edges to the vertices $N_{k_1}$ and $N_{\nu (\Gamma_1)+k_2}.$ 
In the commutative diagram $g(w)$ all the maps in the part of the diagram corresponding to $\Gamma_1$ are exactly those from the commutative diagram $f_1(v_1).$ All the maps in the part of $g(w)$ corresponding to $\Gamma_2$ are exactly those from the commutative diagram $f_2(v_2).$ The new $\pt$-space is mapped to $n_{1, k_1}\in N_{k_1}$ and to $n_{2, k_2}\in N_{\nu(\Gamma_1)+k_2}.$ The only part of $g(w)$ where commutativity is not inherited from the commutativity of $f_1(v_1)$ and $f_2(v_2)$ is the part coming from the new $\pt$-space. In this part the diagram $g(w)$ is commutative by the definition of $W.$

We get the continuous map $g:W\to G_{\frak N, M}^{B_{k_1, k_2}(\Gamma_1, \Gamma_2)}$ defined by $w\to g(w).$ 

For each point $w=(v_1, n_{k_1}, v_2, n_{k_2})\in W\subset V_1\times N_{k_1}\times V_2\times N_{k_2}$ the maps in the commutative diagram $g(w)$ are defined through the projections of 
$w$ to the corresponding coordinates. Since the projections of $W$ to the coordinates $V_1, V_2, N_{k_1}, N_{k_2}$ are smooth, we see that $g$ is a nice map.

The pair $(W, g)$ defines an element $[W,g]\in \bor_{j_1+j_2+2\frak n-\frak m}(G_{\frak N, M}^{B_{k_1, k_2}(\Gamma_1, \Gamma_2)}).$ Using Theorem~\ref{transversalitytheorem} and standard bordism theory type arguments we get that the bordism class $[W,g]$ depends only on $\omega_1, \omega_2.$ We put $A_{k_1, k_2}^{j_1, j_2}(\omega_1, \omega_2)=[W,g].$

{\it When the dimensions of the bordism classes to which $A$ is applied are obvious, we will often omit the upper indices in $A$ and write $A_{k_1, k_2}$ rather than $A_{k_1, k_2}^{j_1, j_2}.$\/}
\end{defin}

\begin{defin}[action of permutation groups on $G_{\frak N,M}$]\label{groupaction}
A permutation $\alpha\in S_{n}$ induces a homeomorphism $\wt \alpha_{\Gamma}: G_{\frak N,M}^{\Gamma}\to G_{\frak N, M}^{\alpha \cdot \Gamma}$ that is defined as follows:
\begin{description}
\item[1] If $n\neq \nu(\Gamma),$ then by Definition~\ref{Boperation} $\alpha \cdot\Gamma=\Gamma$ and 
$\wt \alpha_{\Gamma}$ is put to be the identity auto-homeomorphism $G_{\frak N,M}^{\Gamma}\to G_{\frak N,M}^{\Gamma}=G_{\frak N, M}^{\alpha \cdot \Gamma}.$
\item[2] If $n=\nu(\Gamma),$ then the homeomorphism $\wt \alpha_{\Gamma}$ is given 
by re-enumerating all the $N$-manifolds in the commutative diagrams constituting $G_{\frak N,M}^{\Gamma}$ according to the permutation $\alpha$ and by changing the labels of the $\pt$-spaces accordingly. The manifolds from $\frak N$ corresponding to the vertices and the actual maps in the commutative diagrams are unchanged.
\end{description} 

Since $G_{\frak N,M}=\sqcup_{\Gamma} G_{\frak N,M}^{\Gamma}$ and $\alpha\cdot \Gamma_1\neq \alpha\cdot \Gamma_2$ for $\Gamma_1\neq \Gamma_2,$ we get that the collection of homeomorphisms 
$\wt \alpha_{\Gamma}: G_{\frak N,M}^{\Gamma}\to G_{\frak N, M}^{\alpha \cdot \Gamma}$ defines an auto-homeomorphism $\wt \alpha:G_{\frak N,M}\to G_{\frak N,M}$ with $\wt \alpha|_{G_{\frak N,M}^{\Gamma}}=\wt \alpha_{\Gamma},$ for every allowed graph $\Gamma.$ Thus $\wt \alpha$ induces the automorphism $\alpha_*:\bor_*(G_{\frak N,M})\to \bor_*(G_{\frak N,M}),$ defined by $\alpha_*([V, f])=[V, \wt \alpha \circ f].$
%Note that if $f:V\to G_{\frak N, M}$ is nice, then $\wt \alpha\circ f:V\to %G_{\frak N, M}$ is also nice. 

{\it Since $\wt \alpha_{\Gamma}=\wt \alpha|_{G_{\frak N,M}^{\Gamma}},$ when there is no confusion we will write $\wt \alpha$ instead of $\wt \alpha_{\Gamma}$.  Also we will use $\alpha_*$ instead of $(\alpha_{\Gamma})_*$ to denote the isomorphism 
$\bor_*(G_{\frak N,M}^{\Gamma})\to \bor_*(G_{\frak N,M}^{\alpha \cdot \Gamma}).$\/}
\end{defin}

\begin{figure}[htbp] 
\begin{center} 
\includegraphics[width=11cm]{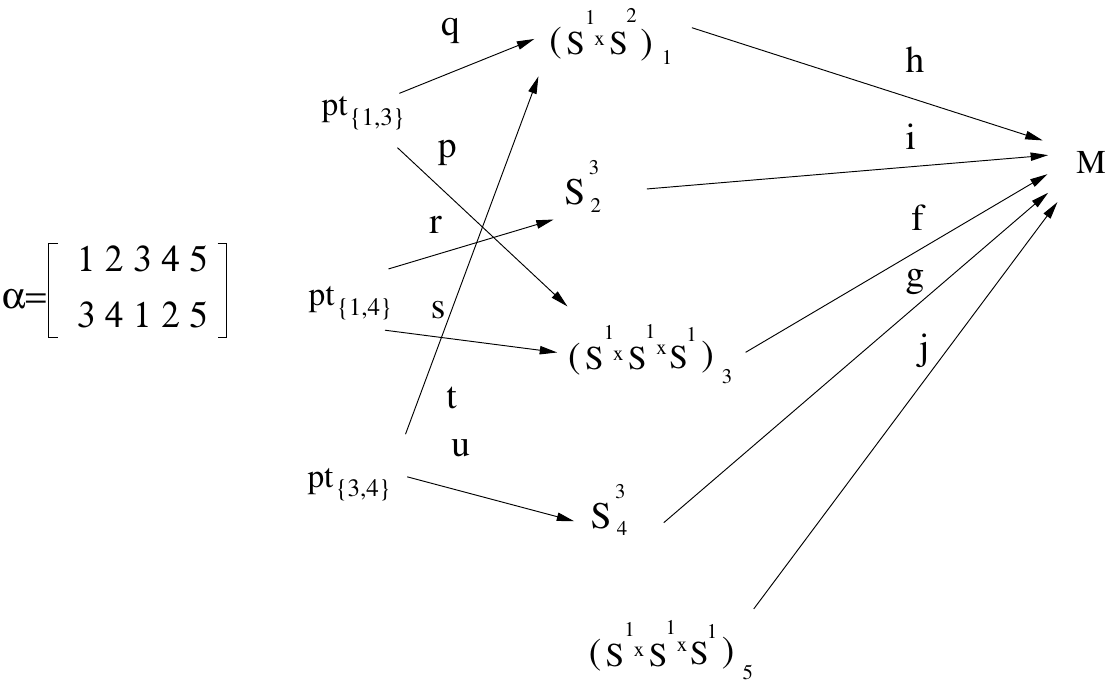} 
%\hepsffile{garlandtreepermutation.eps} 
\end{center} 
\caption{The value of $\wt \alpha$ on the garland in Figure~\ref{garlandtreediag.fig}}\label{garlandtreepermutation.fig} 
\end{figure}

The following Theorem will be essential in the proof of Theorem~\ref{skewsymmetry} stating that the Lie bracket is skew symmetric in the graded sense.

\begin{thm}\label{skewsymmetryingredient}
Let $\omega_i\in \bor_{j_i}(G_{\frak N,M}^{\Gamma_i}), i=1,2$ be bordism classes and let
$k_1\leq \nu (\Gamma_1), k_2\leq \nu (\Gamma_2)$ be positive integers. 
 Then 
$$A_{k_1, k_2}(\omega_1, \omega_2)= (-1)^{(j_1+\frak n)(j_2+\frak n)+\frak m}\alpha_* \bigl( A_{k_2, k_1}(\omega_2, \omega_1)\bigr)$$ for the permutation $\alpha=\bigl(\nu (\Gamma_1), \nu(\Gamma_2)\bigr)\in S_{(\nu(\Gamma_1)+\nu(\Gamma_2))}$ (defined in~\ref{Boperation}).
\end{thm}

The Proof of this Theorem is given in Appendix~\ref{proofskewsymmetryingredient}

\begin{remark}\label{strangesign} 
Let $h_1:X_1^{l_1}\to M^{\frak m}$ and $h_2:X_2^{l_2}\to M^{\frak m}$ be smooth transverse maps of closed oriented manifolds.
According to our orientation convention the bordism class in $\bor_{l_1+l_2-\frak m}(M)$ is $(h_1\times h_2)^{-1}(\Delta)$ which is $(-1)^{l_1l_2+\frak m}(h_2\times h_1)^{-1}(\Delta).$
For this reason the sign in Theorem~\ref{skewsymmetry} is indeed $(-1)^{(j_1+\frak n)(j_2+\frak n)+\frak m},$ rather than $(-1)^{(\frak m-(j_1+\frak n))(\frak m-(j_2+\frak n))}$ that one might expect by the analogy with the intersection pairing on $H_*(M).$
\end{remark}

The following Theorem will be essential in the proof of Theorem~\ref{Jacobi} stating that the Lie bracket satsfies the graded Jacobi identity.

\begin{thm}\label{Jacobiingredient}
Let $\omega_i\in \bor_{j_i}(G_{\frak N,M}^{\Gamma_i}), i=1,2,3,$ be bordism classes and 
let $k_1\leq \nu (\Gamma_1),$ $k_2,\ov k_2\leq \nu (\Gamma_2), k_3\leq \nu (\Gamma_3)$ be positive integers. 
Then 
\begin{equation}
\begin{split}
A_{\nu (\Gamma_1)+\ov k_2, k_3}\bigl (A_{k_1, k_2}(\omega_1, \omega_2), \omega_3\bigr)=(-1)^{\sigma}\alpha_*  \Bigl( A_{k_2, k_1} \bigl (A_{\ov k_2, k_3} (\omega_2, \omega_3), \omega_1\bigr)\Bigr)
\end{split}
\end{equation}
for $\sigma=(\frak m+j_1)j_3+(j_1+\frak n_{k_1})(j_2+\frak n_{k_2})+\frak m j_1$ and
the permutation $\alpha=\bigl(\nu(\Gamma_1), \nu (\Gamma_2)+\nu(\Gamma_3) \bigr)\in S_{(\nu(\Gamma_1)+\nu(\Gamma_2)+\nu(\Gamma_3))}$ 
(defined in~\ref{Boperation}).
\end{thm}

The Proof of this Theorem is given in Appendix~\ref{proofJacobiingredient}

\section{Graded Lie algebra on $\bor_*(G_{\frak N, M})$}\label{bracket}

Recall that we use the following {\bf convention:\/} given $q\in \Q$ and $\omega\in \bor_*(G_{\frak N, M})$ we will write $q\omega$ to denote the element $\omega\otimes q \in \bor_*(G_{\frak N, M})_\Q=\bor_*(G_{\frak N, M})\otimes \Q.$ Note that if $k\in \Z\subset \Q,$ then, depending on the context, $k\omega$ could stand either for $\omega\otimes k\in \bor_*(G_{\frak N, M})\otimes \Q$ or for $k\omega\in \bor_*(G_{\frak N, M}).$
We would also often write $1\omega$ to stress that we are talking about $\omega\otimes 1\in \bor_*(G_{\frak N, M})_\Q$ rather than about $\omega\in \bor_*(G_{\frak N, M}).$ Note that $1\omega$ is zero if $\omega$ is an element of finite order.

\begin{defin}[Lie bracket] Below we define $[\cdot, \cdot]:\bor_*(G_{\frak N,M})_\Q \otimes_{\Q} \bor_*(G_{\frak N,M})_\Q)\to \bor_*(G_{\frak N,M})_\Q.$

Since $\bor_*(G_{\frak N,M})_\Q=\oplus_{j=0}^{\infty} (\bor_j(G_{\frak N,M})_\Q),$ it is a graded $\Q$-vector space with the grading of $\bor_{j}(G_{\frak N,M})\otimes \Q$ equal to $j.$ 

Since the addition in $\bor_*(G_{\frak N,M})$ is by taking the disjoint union, it suffices to define $[1\omega_1, 1\omega_2]$ for the case where $\omega_1, \omega_2$ can each be realized by a map of a connected manifold. Then one extends $[\cdot, \cdot]$ to the whole 
$\bor_*(G_{\frak N,M})_\Q\otimes_{\Q} \bor_*(G_{\frak N,M})_\Q$
by linearity.

Let $f_i:V_i^{j_i}\to G_{\frak N,M}, i=1,2,$ be maps of connected manifolds realizing 
$\omega_{i}\in \bor_{j_i}(G_{\frak N,M}).$ Since $V_i$ are connected, $\Im f_i\subset G_{\frak N,M}^{\Gamma_i},$ for some $\Gamma_i, i=1,2.$ 
Thus $\omega_i$ can be regarded as elements of $\bor_{j_i}(G_{\frak N, M}^{\Gamma_i}), i=1,2.$
Put 
\begin{equation}
[1\omega_1,1\omega_2]=\frac{1}{(\nu(\Gamma_1)+\nu(\Gamma_2))!}\bigl(\sum_{k_1=1}^{\nu(\Gamma_1)}\sum_{k_2=1}^{\nu(\Gamma_2)} \sum_{\beta\in S_{(\nu(\Gamma_1)+\nu(\Gamma_2))}} \beta_*\bigl(A_{k_1, k_2}(\omega_1, \omega_2)\bigr)\bigr).
\end{equation}
\end{defin}

The following Theorem follows immediately from the definitions of $[\cdot, \cdot]$ and of $A_{k_1, k_2}.$

\begin{thm}\label{dimension}
$[\cdot, \cdot]: \bor_{j_1}(G_{\frak N,M})_\Q\otimes_{\Q} \bor_{j_2}(G_{\frak N,M})_\Q\to \bor_{j_1+j_2+2\frak n -\frak m}(G_{\frak N,M})_\Q$ is a bilinear operation.
\qed
\end{thm}

\begin{thm}\label{skewsymmetry} Let $\eta_{i}\in \bor_{j_i}(G_{\frak N,M})_\Q,i=1,2$ be elements, then 
$$[\eta_1, \eta_2]=(-1)^{(j_1+1)(j_2+1)+\frak m}[\eta_2, \eta_1].$$
\end{thm}

\pp It suffices to prove the theorem for $\eta_i=1\omega_i,$ where both of $\omega_i\in \bor_{j_i}(G_{\frak N}, M),$ $ i=1,2,$ are realizable by a map of a connected oriented manifold. 
Let $f_i:V_i^{j_i}\to G_{\frak N,M}$ be maps of connected manifolds realizing $\omega_i, i=1,2.$
Since $V_i$ are connected $\Im f_i\subset G_{\frak N,M}^{\Gamma_i}$ for some $\Gamma_i, i=1,2.$
For brevity we write
$\nu_1, \nu_2$ instead of $\nu(\Gamma_1), \nu(\Gamma_2).$

By definition $[1\omega_1, 1\omega_2]= \frac{1}{(\nu_1+\nu_2)!}\sum_{k_1, k_2=1}^{\nu_1, \nu_2}\sum_{\beta\in S_{(\nu_1+\nu_2)}} \beta_*(A_{k_1, k_2}(\omega_1, \omega_2)).$
By Theorem~\ref{skewsymmetryingredient} we get 
$$[1\omega_1, 1\omega_2]=\frac{1}{(\nu_1+\nu_2)!}\sum_{k_1, k_2=1}^{\nu_1, \nu_2} \sum_{\beta\in S_{(\nu_1+\nu_2})} \beta_*\Bigl((-1)^{(j_1+1)(j_2+1)+\frak m}\alpha_*\bigl(A_{k_2, k_1}(\omega_2, \omega_1)\bigr)\Bigr),$$ 
for the permutation $\alpha=\bigl(\nu_1, \nu_2\bigr)\in S_{(\nu_1+\nu_2)}.$ Since $(\beta\alpha)_*=\beta_*\alpha_*,$ we make a substitution $\beta'=\beta\alpha$ and get that 
\begin{equation}
\begin{split}
[1\omega_1, 1\omega_2]=\frac{(-1)^{(j_1+1)(j_2+1)+\frak m}}{(\nu_2+\nu_1)!}\bigl(\sum_{k_2, k_1=1}^{\nu_2, \nu_1} \sum_{\beta'\in S_{(\nu_2+\nu_1)}} \beta'_*\bigl(A_{k_2, k_1}(\omega_2, \omega_1)\bigr)\bigr)\\=(-1)^{(j_1+1)(j_2+1)+\frak m}[1\omega_2, 1\omega_1].\qed
\end{split}
\end{equation}

\m

Theorems~\ref{dimension},~\ref{skewsymmetry} and the following Theorem~\ref{Jacobi} say that  the operation $[\cdot, \cdot]$ gives $\bor_*(G_{\frak N,M})_\Q$ the structure of a graded Lie algebra. 

\begin{thm}\label{Jacobi}
We have the following graded Jacobi identity 
\begin{equation}
\begin{split}
(-1)^{(\frak mj_3+j_1j_3+j_1)}[[\eta_1, \eta_2],\eta_3]+(-1)^{(\frak mj_1+j_1j_2+j_2)}[[\eta_2, \eta_3], \eta_1]+\\(-1)^{(\frak mj_2+j_2j_3+j_3)}[[\eta_3, \eta_1],\eta_2]=0
\end{split}
\end{equation}
for  all $\eta_i\in \bor_{j_i}(G_{\frak N,M})_\Q, i=1,2,3.$ 
\end{thm}

\section{Proof of Theorem~\ref{Jacobi}}\label{ProofJacobi}

To prove Theorem~\ref{Jacobi} we will need the following Proposition.

\begin{prop}\label{technical}
Let $\omega_i\in \bor_{j_i}(G_{\frak N,M}^{\Gamma_i}), i=1,2,3,$ be bordism classes, and let $\nu_i=\nu(\Gamma_i), i=1,2,3.$ Then 
\begin{equation}
\begin{split}
[[1\omega_1, 1\omega_2], 1\omega_3]=\\ \frac{1}{(\nu_1+\nu_2+\nu_3)!}\Bigl(\sum_{\beta\in S}\beta_*\Bigl(\sum_{\ov  k_1=1}^{\nu_1+\nu_2}
\sum_{k_3=1}^{\nu_3}\sum_{k_1=1}^{\nu_1}\sum_{k_2=1}^{\nu_2} A_{\ov  k_1, k_3}\bigl (A_{k_1, k_2} (\omega_1, \omega_2), \omega_3\bigr)\Bigr)\Bigr),
\end{split}
\end{equation}
where $S$ is the permutation group $S_{(\nu_1+\nu_2+\nu_3)}.$
\end{prop}

\pp  
We denote $S_{\nu_1+\nu_2}$ by $\wt S$ and
we will often identify $\wt S$ with the subgroup of $S=S_{(\nu_1+\nu_2+\nu_3)}$ that consists of group elements acting trivially on the last $\nu_3$ elements in $\{1, 2, \cdots, (\nu_1+\nu_2+\nu_3)\}.$ For a permutation $\gamma\in \wt S,$ we denote by $\ov \gamma$ the corresponding permutation in $S.$

By definition of $[\cdot, \cdot]$ we have 
\begin{equation}\label{eqproptechnical}
\begin{split}(\nu_1+\nu_2)!(\nu_1+\nu_2+\nu_3)!\times [[1\omega_1, 1\omega_2], 1\omega_3]=\\ 
\sum_{\ov  k_1, k_3=1}^{\nu_1+\nu_2, \nu_3} \sum_{\beta\in S}\beta_* \Bigl( A_{\ov  k_1, k_3}\Bigl (\sum_{k_1, k_2=1}^{\nu_1, \nu_2}\sum_{\gamma\in \wt S}
 \gamma_* \bigl(A_{k_1, k_2}(\omega_1, \omega_2)\bigr), \omega_3\Bigr)\Bigr)=\\ \sum_{\beta\in S}\sum_{\ov  k_1, k_3=1}^{\nu_1+\nu_2, \nu_3} \beta_* \sum_{\gamma\in \wt S}
\ov \gamma_*
\Bigl( A_{\gamma^{-1}(\ov  k_1), k_3}\Bigl( \Bigl(\sum_{k_1, k_2=1}^{\nu_1, \nu_2}\bigl(A_{k_1, k_2}(\omega_1, \omega_2)\bigr)\Bigr), \omega_3\Bigr)\Bigr) =\\
\sum_{\beta\in S, \gamma\in \wt S}(\beta\ov \gamma)_*\sum_{\ov  k'_1, k_3=1}^{\nu_1+\nu_2, \nu_3} 
\Bigl( A_{\ov  k'_1, k_3}\Bigl (\sum_{k_1, k_2=1}^{\nu_1, \nu_2}\bigl(A_{k_1, k_2}(\omega_1, \omega_2)\bigr), \omega_3\Bigr)\Bigr),
\end{split} 
\end{equation}
where for a fixed $\gamma\in \wt S$ we made a substitution $\ov  k'_1=\gamma^{-1}(\ov k_1).$ Each permutation $\beta'$ in $S$ is realized in $(\nu_1+\nu_2)!$
ways as $\beta\ov \gamma,$ for $\beta\in S$ and $\gamma\in \wt S.$ Thus the right hand side of~\eqref{eqproptechnical} is equal to 

\begin{equation}\label{eq2proptechnical}
\begin{split}
(\nu_1+\nu_2)!
\Bigl(\sum_{\beta'\in S}\beta'_*\Bigl( \sum_{\ov  k'_1, k_3=1}^{\nu_1+\nu_2, \nu_3} 
\sum_{k_1, k_2=1}^{\nu_1, \nu_2}A_{\ov  k'_1, k_3}\bigl(A_{k_1, k_2}(\omega_1, \omega_2), \omega_3\bigr)\Bigr)\Bigr). 
\end{split}
\end{equation}
To get the proof divide the left hand side of~\eqref{eqproptechnical} and the equal to it expression~\eqref{eq2proptechnical} by $(\nu_1+\nu_2)!(\nu_1+\nu_2+\nu_3)!.$\qed

\m
{\it Let us prove Theorem~\ref{Jacobi}.\/} We observe that it suffices to prove the Theorem for $\eta_i=1\omega_i$ where $\omega_i, i=1,2,3,$ are realizable by maps of connected manifolds. 

Let $f_i:V_i^{j_i}\to G_{\frak N,M}$ be maps of connected manifolds realizing $\omega_i, i=1,2,3.$ Since $V_i$ are connected, $\Im f_i\subset G_{\frak N,M}^{\Gamma_i}$ for some $\Gamma_i, i=1,2,3.$ For brevity we put $\nu_i=\nu(\Gamma_i), i=1,2,3.$ Recall that $S$ is the permutation group 
$S_{(\nu_1+\nu_2+\nu_3)}.$

Put 
\begin{equation}\label{AB}
\begin{split}
\mathfrak A =\sum_{\beta\in S}\beta_*\Bigl(\sum_{\ov  k_1, k_3=1}^{\nu_1, \nu_3}
\sum_{k_1, k_2=1}^{\nu_1, \nu_2}A_{\ov  k_1, k_3}(A_{k_1, k_2} (\omega_1, \omega_2), \omega_3)\Bigr)\\
\text{ and put }\mathfrak B =\sum_{\beta\in S}\beta_*\Bigl(\sum_{\ov  k_1, k_3=1}^{\nu_2, \nu_3}
\sum_{k_1, k_2=1}^{\nu_1, \nu_2} A_{\nu_1+\ov  k_1, k_3}(A_{k_1, k_2} (\omega_1, \omega_2), \omega_3)\Bigr).
\end{split}
\end{equation}

By Proposition~\ref{technical} the first term of the graded Jacobi identity is
\begin{equation}
\begin{split}\label{firsttermeq}
(-1)^{(\frak mj_3+j_1j_3+j_1)}[[1\omega_1, 1\omega_2],1\omega_3]=\\\frac{(-1)^{(\frak mj_3+j_1j_3+j_1)}}{(\nu_1+\nu_2+\nu_3)!}\sum_{\beta\in S}\beta_*\Bigl(\sum_{\ov  k_1, k_3=1}^{\nu_1+\nu_2, \nu_3}
\sum_{k_1, k_2=1}^{\nu_1, \nu_2}A_{\ov  k_1, k_3}\bigl (A_{k_1, k_2} (\omega_1, \omega_2), \omega_3\bigr)\Bigr)=\\
\frac{(-1)^{(\frak mj_3+j_1j_3+j_1)}}{(\nu_1+\nu_2+\nu_3)!}
\Bigl( \mathfrak A + \mathfrak B\Bigr),
\end{split}
\end{equation}
where we split the huge sum into two parts $\frak A$ and $\frak B$ based on whether $\ov  k_1$ was indexing a manifold $N$ that came from $\omega_1$ or a manifold $N$ that came from $\omega_2.$

Similarly put 
\begin{equation}\label{CD}
\begin{split}
\mathfrak C =\sum_{\beta\in S}\beta_*\Bigl(\sum_{\ov  k_2, k_1=1}^{\nu_2, \nu_1}
\sum_{k_2, k_3=1}^{\nu_2, \nu_3}A_{\ov  k_2, k_1}(A_{k_2, k_3} (\omega_2, \omega_3), \omega_1)\Bigr)\\
\text{ and put }\mathfrak D =\sum_{\beta\in S}\beta_*\Bigl(\sum_{\ov  k_2, k_1=1}^{\nu_3, \nu_1}\sum_{k_2, k_3=1}^{\nu_2, \nu_3} A_{\nu_2+\ov  k_2, k_1}(A_{k_2, k_3} (\omega_2, \omega_3), \omega_1)\Bigr).
\end{split}
\end{equation}

By Proposition~\ref{technical} the second term of the graded Jacobi identity is
\begin{equation}\label{secondtermeq}
\begin{split}
(-1)^{(\frak mj_1+j_1j_2+j_2)}[[1\omega_2, 1\omega_3],1\omega_1]=\\
\frac{(-1)^{(\frak mj_1+j_1j_2+j_2)}}
{(\nu_1+\nu_2+\nu_3)!}
\sum_{\beta\in S}\beta_*\Bigl(\sum_{\ov  k_2, k_1=1}^{\nu_2+\nu_3, \nu_1}
\sum_{k_2, k_3=1}^{\nu_2, \nu_3} A_{\ov  k_2, k_1}(A_{k_2, k_3} (\omega_2, \omega_3), \omega_1)\Bigr)=\\
\frac{(-1)^{(\frak mj_1+j_1j_2+j_2)}}{(\nu_1+\nu_2+\nu_3)!}
\Bigl( \mathfrak C+ \mathfrak D\Bigr),
\end{split}
\end{equation}
where we split the huge sum into two parts $\frak C$ and $\frak D$ based on whether $\ov  k_2$ was indexing a manifold $N$ that came from $\omega_2$ or a manifold $N$ that came from $\omega_3.$

Similarly put 
\begin{equation}\label{EF}
\begin{split}
\mathfrak E =\sum_{\beta\in S}\beta_*\Bigl(\sum_{\ov  k_3, k_2=1}^{\nu_3, \nu_2}
\sum_{k_3, k_1=1}^{\nu_3, \nu_1}A_{\ov  k_3, k_2}(A_{k_3, k_1} (\omega_3, \omega_1), \omega_2)\Bigr)\\
\text{ and put }\mathfrak F =\sum_{\beta\in S}\beta_*\Bigl(\sum_{\ov  k_3, k_2=1}^{\nu_1, \nu_2}
\sum_{k_3, k_1=1}^{\nu_3, \nu_1} A_{\nu_3+\ov  k_3, k_2}(A_{k_3, k_1} (\omega_3, \omega_1), \omega_2)\Bigr).
\end{split}
\end{equation}

By Proposition~\ref{technical} the third term of the graded Jacobi identity is
\begin{equation}\label{thirdtermeq}
\begin{split}
(-1)^{(\frak mj_2+j_2j_3+j_3)}[[1\omega_3, 1\omega_1],1\omega_2]=\\
\frac{(-1)^{(\frak mj_2+j_2j_3+j_3)}}{(\nu_1+\nu_2+\nu_3)!}
\Bigl(\sum_{\beta\in S}\beta_*\Bigl(\sum_{\ov  k_3, k_2=1}^{\nu_3+\nu_1, \nu_2}
\sum_{k_3, k_1=1}^{\nu_3, \nu_1}A_{\ov  k_3, k_2}(A_{k_3, k_1} (\omega_3, \omega_1), \omega_2)\Bigr)\Bigr)=\\
\frac{(-1)^{(\frak mj_2+j_2j_3+j_3)}}{(\nu_1+\nu_2+\nu_3)!}
\Bigl( \mathfrak E+ \mathfrak F\Bigr),
\end{split}
\end{equation}
where we split the huge sum into two parts $\frak E$ and $\frak F$ based on whether $\ov  k_3$ was indexing a manifold $N$ that came from $\omega_3$ or a manifold $N$ that came from $\omega_1.$

Apply identities~\eqref{firsttermeq},~\eqref{secondtermeq},~\eqref{thirdtermeq} to the left hand side of the graded Jacobi identity to get 
\begin{equation}\label{Jacobisplit}
\begin{split}
(-1)^{(\frak mj_3+j_1j_3+j_1)}[[1\omega_1, 1\omega_2],1\omega_3]+(-1)^{(\frak mj_1+j_1j_2+j_2)}[[1\omega_2, 1\omega_3], 1\omega_1]+\\
(-1)^{(\frak mj_2+j_2j_3+j_3)}[[1\omega_3, 1\omega_1],1\omega_2]=\\
\frac{1}{(\nu_1+\nu_2+\nu_3)!}
\Bigl( (-1)^{(\frak mj_3+j_1j_3+j_1)}(\mathfrak A + \mathfrak B)+
(-1)^{(\frak mj_1+j_1j_2+j_2)}(\mathfrak C+ \mathfrak D)+\\(-1)^{(\frak mj_2+j_2j_3+j_3)}(\mathfrak E+\mathfrak F)\Bigr).
\end{split}
\end{equation}

Using Proposition~\ref{identities} we get that this expression is zero. This finishes the proof of Theorem~\ref{Jacobi} modulo the proof of Proposition~\ref{identities}.\qed

\begin{prop}\label{identities}
The following identities hold
\begin{equation}
\begin{split}
(-1)^{(\frak mj_3+j_1j_3+j_1)}\mathfrak A =- (-1)^{(\frak mj_2+j_2j_3+j_3)}\mathfrak F\\
(-1)^{(\frak mj_3+j_1j_3+j_1)}\mathfrak B= - (-1)^{(\frak mj_1+j_1j_2+j_2)}\mathfrak C\\
(-1)^{(\frak mj_1+j_1j_2+j_2)}\mathfrak D= -(-1)^{(\frak mj_2+j_2j_3+j_3)}\mathfrak E.
\end{split}
\end{equation}
\end{prop}

\pp Let us prove the first identity. 

Apply Theorem~\ref{Jacobiingredient} to the ordered triple 
$\omega_3, \omega_1, \omega_2$ of bordism classes to get
$$A_{\nu (\Gamma_3)+\ov k_1, k_2}\bigl (A_{k_3, k_1}(\omega_3, \omega_1), \omega_2\bigr)= (-1)^{\sigma_1}\alpha_*  \Bigl( A_{k_1, k_3} \bigl (A_{\ov k_1, k_2} (\omega_1, \omega_2), \omega_3\bigr)\Bigr),$$
for $\sigma_1=(\frak m+j_3)j_2+(j_3+1)(j_1+1)+\frak mj_3,$ the permutation $\alpha= \bigl(\nu_3, \nu_1+\nu_2\bigr)\in S=S_{\nu_1+\nu_2+\nu_3},$ 
and for $k_1, \ov k_1\leq \nu_1, k_2 \leq \nu_2, k_3\leq \nu_3.$

Apply this identity with $\ov k_3$ instead of $\ov k_1$ to the equation~\eqref{EF} defining $\frak F$ to get

\begin{equation}
\begin{split}
- (-1)^{(\frak mj_2+j_2j_3+j_3)}\mathfrak F =
\\(-1)^{(1+\frak mj_2+j_2j_3+j_3)}\sum_{\beta\in S}\beta_*\Bigl(\sum_{\ov  k_3, k_2=1}^{\nu_1, \nu_2}
\sum_{k_3, k_1=1}^{\nu_3, \nu_1} A_{\nu_3+\ov  k_3, k_2}(A_{k_3, k_1} (\omega_3, \omega_1), \omega_2)\Bigr)=
\\(-1)^{\sigma_1+(1+\frak mj_2+j_2j_3+j_3)}\sum_{\beta\in S}\beta_*\Bigl(\sum_{\ov  k_3, k_2=1}^{\nu_1, \nu_2}
\sum_{k_3, k_1=1}^{\nu_3, \nu_1}\alpha_*  \Bigl( A_{k_1, k_3} \bigl (A_{\ov k_3, k_2} (\omega_1, \omega_2), \omega_3\bigr)\Bigr) \Bigr)=\\
(-1)^{\frak mj_3+j_1j_3+j_1}\sum_{\beta\in S}(\beta_*\alpha_*)\Bigl(\sum_{\ov  k_3, k_2=1}^{\nu_1, \nu_2}
\sum_{k_3, k_1=1}^{\nu_3, \nu_1}\Bigl( A_{k_1, k_3} \bigl (A_{\ov k_3, k_2} (\omega_1, \omega_2), \omega_3\bigr)\Bigr) \Bigr)
\end{split}
\end{equation}
Make a substitution $\beta'=\beta\alpha,$ rearrange the order of addition, and use equation~\eqref{AB} defining 
$\frak A$ 
to get that this expression equals
\begin{equation}
\begin{split}
%(-1)^{\frak mj_3+j_1j_3+j_1}\Bigl(\sum_{\beta'\in S}\beta'_*\Bigl(\sum_{\ov  k_3, k_2=1}^{\nu_1, \nu_2}
%\sum_{k_3, k_1=1}^{\nu_3, \nu_1}\Bigl( A_{k_1, k_3} \bigl (A_{\ov k_3, k_2} (\omega_1, \omega_2), %\omega_3\bigr)\Bigr) \Bigr)\Bigr)=\\
(-1)^{(\frak mj_3+j_1j_3+j_1)}\Bigl(\sum_{\beta'\in S}\beta'_*\Bigl(\sum_{k_1, k_3=1}^{\nu_1, \nu_3}\sum_{\ov  k_3, k_2=1}^{\nu_1, \nu_2}
\Bigl( A_{k_1, k_3} \bigl (A_{\ov k_3, k_2} (\omega_1, \omega_2), \omega_3\bigr)\Bigr) \Bigr)\Bigr)=\\
(-1)^{(\frak mj_3+j_1j_3+j_1)}\mathfrak A.
\end{split}
\end{equation}
This proves the first of the three identities. The proofs of the other two identities are obtained similarly. This finishes the proof of Proposition~\ref{identities} and of the Theorem~\ref{Jacobi}.
\qed

\section{Operation $C$ on $\bor_*(G_{\frak N, M})$}\label{SectionCoperation}
In this section we introduce the operation $C$ on $\bor_*(G_{\frak N, M}).$ In Section~\ref{multiplication} we will use the symmetrization of $C$ to construct the product $\star$ on $\bor_*(G_{\frak N, M})_\Q.$
Operations $\star$ and $[\cdot, \cdot]$ give a graded Poisson algebra structure on $\bor_*(G_{\frak N, M})_\Q,$  see Theorem~\ref{associativecommutative} and \ref{Poisson}.

\begin{defin}[operation $D$ on graphs and one more important permutation]\label{Doperation}
Let $\Gamma_1$ and $\Gamma_2$ be two allowed graphs. Let $\nu(\Gamma_i)$ be the number of $N$-vertices in the graph $\Gamma_i, i=1,2.$ We define the allowed graph $D(\Gamma_1, \Gamma_2)$ as the graph resulting after the following sequence of operations: 
\begin{description}
\item[1]
Take the disjoint union of the graphs $\Gamma_1$ and $\Gamma_2,$ preserving the manifolds in $\frak N$ that were associated to the $N$-vertices in the graphs $\Gamma_1$ and $\Gamma_2.$
\item[2]
Keep the indices of the $N$-vertices that came from $\Gamma_1$ and increase by $\nu(\Gamma_1)$ the indices of the $N$-vertices that came from $\Gamma_2.$ 
\item[3]
Change the  indices of the $\pt$-vertices in $\Gamma_2,$ so that the new index of each $\pt$-vertex gives the indices (with respect to the shifted enumeration of $N$-vertices) of the $N$-vertices in $\Gamma_2$ connected to this $\pt$-vertex by an oriented edge.
\item[4] Identify the two $M$-vertices in $\Gamma_1\sqcup \Gamma_2$ to get just one $M$ vertex and redirect 
all the $N\to M$ edges in the two graphs to it.
\end{description}

For $n_1, n_2, n_3\in \N$ {\em define a permutation\/} $(n_1, n_2, n_3)\in S_{n_1+n_2+n_3}$ by putting its value on $\{1, 2, \cdots, n_1+n_2+n_3\}$ to be $\{n_1+1, n_1+2, \cdots, n_1+n_2, 1, 2, \cdots, n_1, n_1+n_2+1, n_1+n_2+2, \cdots, n_1+n_2+n_3\}.$

\end{defin}

The following Proposition follows immediately from the definition of the operation $D.$

\begin{prop}\label{propgraphD}
Let $\Gamma_1, \Gamma_2,\Gamma_3$ be three allowed graphs, then 
\begin{description}
\item[1] $D\bigl (D(\Gamma_1, \Gamma_2), \Gamma_3\bigr)=D\bigl (\Gamma_1, D(\Gamma_2, \Gamma_3)\bigr).$
\item[2] $D(\Gamma_1, \Gamma_2)=\alpha_1\cdot D(\Gamma_2, \Gamma_1),$ for the permutation $\alpha_1=\bigl(\nu(\Gamma_1), \nu(\Gamma_2)\bigr)\in S_{(\nu(\Gamma_1)+\nu(\Gamma_2))}$ (defined in~\ref{Boperation}).
\item[3] $B_{k_1, k_2}\bigl(\Gamma_1, D(\Gamma_2, \Gamma_3)\bigr)=D\bigl(B_{k_1, k_2}(\Gamma_1, \Gamma_2), \Gamma_3\bigr),$ for all positive integer $k_1\leq \nu (\Gamma_1), k_2\leq \nu(\Gamma_2).$
\item[4] $B_{k_1, \nu(\Gamma_2)+k_3}\bigl(\Gamma_1, D(\Gamma_2, \Gamma_3)\bigr)=\alpha_2\cdot D\bigl(\Gamma_2, B_{k_1, k_3}(\Gamma_1, \Gamma_3)\bigr),$  for the permutation $\alpha_2=\bigl(\nu(\Gamma_1), \nu (\Gamma_2), \nu(\Gamma_3)\bigr)\in S_{(\nu(\Gamma_1)+\nu(\Gamma_2)+\nu(\Gamma_2))}$ (defined in~\ref{Coperation}) and all positive integer $k_1\leq \nu (\Gamma_1), k_3\leq \nu(\Gamma_3).$
\end{description}
\qed
\end{prop}

\begin{defin}[operation $C^{j_1, j_2}:\bor_{j_1}(G_{\frak N, M}^{\Gamma_1})\otimes \bor_{j_2}(G_{\frak N, M}^{\Gamma_2})\to \bor_{j_1+j_2}\bigl(G_{\frak N, M}^{D(\Gamma_1, \Gamma_2)}\bigr)$]\label{Coperation}

Let $f_i:V_i^{j_i}\to \bor_{j_i}(G_{\frak N, M}^{\Gamma_i})$ be maps of oriented manifolds realizing oriented bordism classes $\omega_i\in \bor_{j_i}(G_{\frak N, M}^{\Gamma_i}), i=1,2.$ 
Thus each $f_i(v_i), v_i\in V_i,$ is a commutative diagram from $G_{\frak N, M}^{\Gamma_i}, i=1,2.$

Consider the $(j_1+j_2)$-dimensional manifold $X=V_1\times V_2$ oriented as the product of oriented manifolds. 
Each point $x=(v_1, v_2)\in X$ gives rise to a commutative diagram $g(x)\in G_{\frak N, M}^{D(\Gamma_1, \Gamma_2)}$ as follows. Recall that by Definition~\ref{Doperation}, the graph $D(\Gamma_1, \Gamma_2)$  consists of two parts: the part coming from $\Gamma_1$ and the part coming from $\Gamma_2.$  
In the commutative diagram $g(x)$ all the maps in the part of the diagram corresponding to $\Gamma_1$ are exactly those from the commutative diagram $f_1(v_1).$ All the maps in the part of $g(x)$ corresponding to $\Gamma_2$ are exactly those from the commutative diagram $f_2(v_2).$ 

We get the continuous map $g:X\to G_{\frak N, M}^{D(\Gamma_1, \Gamma_2)}, x\to g(x).$
The pair $(X, g)$ defines an element $[X,g]\in \bor_{j_1+j_2}\bigl (G_{\frak N, M}^{D(\Gamma_1, \Gamma_2)}\bigr).$
Standard bordism theory arguments show that $[X,g]$ depends only on the bordism classes $\omega_1, \omega_2$ and we put $C^{j_1, j_2}(\omega_1, \omega_2)=[X,g].$ 

%For every $x=(v_1, v_2)\in X=V_1\times V_2$ all the %maps %in the commutative diagram $g(x)$ are defined by the %projections of $x$ to the coordinates $V_1$ and $V_2.$ %Since the projections $V_1\times V_2\to V_1$ and %$V_1\times V_2\to V_2$ are smooth, 
%we get that $g$ is a nice map.

{\it When the dimensions of the bordism classes to which $C$ is applied are obvious, we will often omit the indices in $C.$ So we will write $C$ rather than $C^{j_1, j_2}.$\/}
\end{defin}

The following Theorem will be needed in the proof of Theorem~\ref{associativecommutative} saying that  $\star$ is an associative graded commutative product on $\bor_*(G_{\frak N, M})_\Q,$ for all $\frak N.$

\begin{thm}\label{multiplicationingredient}
Let $\omega_i\in \bor_{j_i}(G_{\frak N,M}^{\Gamma_i}), i=1,2,3,$ be bordism classes. Then
\begin{description}
\item[1] $C(\omega_1, \omega_2)=(-1)^{j_1+j_2}\alpha_* (C(\omega_2, \omega_1)),$ where $\alpha=(\nu(\Gamma_1), \nu(\Gamma_2))$ is the permutation (defined in~\ref{Boperation}).
\item[2] $C\bigl(\omega_1, C(\omega_2, \omega_3)\bigr)=C\bigl(C(\omega_1, \omega_2), \omega_3\bigr).$
\end{description}
\end{thm}

The proof of this Theorem is straighforward. \qed

The following Theorem will be needed in the proof of Theorem~\ref{Poisson} saying that 
the $\star$-product satisfies a graded Leibniz rule with respect to $[\cdot, \cdot].$
%\m

\begin{thm}\label{Poissoningredient}
Let $\omega_i\in \bor_{j_i}(G_{\frak N,M}^{\Gamma_i}), i=1,2,3,$ be bordism classes and 
let $k_i\leq \nu (\Gamma_i), i=1,2,3,$ be positive integers. 
Then 
\begin{description}
\item[1]
$A_{k_1, \nu(\Gamma_2)+k_3}\bigl (\omega_1, C(\omega_2, \omega_3)\bigr)=(-1)^{j_2(j_1+\frak n)}\alpha_*\bigl(C\bigl(\omega_2, A_{k_1, k_3}(\omega_1, \omega_3)\bigr)\bigr),$
for the permutation $\alpha=\bigl(\nu(\Gamma_1), \nu (\Gamma_2), \nu(\Gamma_3)\bigr)\in S_{(\nu(\Gamma_1)+\nu(\Gamma_2)+\nu(\Gamma_3))}$ (defined in~\ref{Doperation}).
\item[2]
$A_{k_1, k_2}\bigl(\omega_1, C(\omega_2, \omega_3)\bigr)=(-1)^{j_3(\frak m+\frak n)}C\bigl(A_{k_1, k_2}(\omega_1, \omega_2), \omega_3\bigr);$
\end{description}
\end{thm}

%The proof of Theorem~\ref{Poissoningredient} is given in the Appendix~\ref{proofPoissoningredient}.

The proof of Theorem~\ref{Poissoningredient} is similar to the proof of Theorem~\ref{Jacobiingredient} and is omitted. The interested reader can find a detailed proof in the fourth version of the arXiv preprint version of this work.

\section{Product $\star$ and the graded Leibnitz identity}\label{multiplication}
We define an operation $\star:\bor_*(G_{\frak N,M})_\Q \otimes_{\Q} \bor_*(G_{\frak N,M})_\Q\to \bor_*(G_{\frak N,M})_\Q.$ 
The operation $\star$ will give $\bor_*(G_{\frak N,M})_\Q$ the structure of an associative graded commutative 
algebra. 

Since the addition in $\bor_*(G_{\frak N,M})$ is by taking the disjoint union, it suffices to define $(1\omega_1)\star (1\omega_2)$ for the case where the bordism classes $\omega_1, \omega_2$ can each be realized by a nice map of a connected manifold. Then one extends the $\star$-product
by linearity. For brevity we will often write $1\omega_1\star1\omega_2$ rather than $(1\omega_1)\star(1\omega_2).$

\begin{defin}\label{multiplicationdefin}
Let $f_i:V_1^{j_i}\to G_{\frak N,M}, i=1,2,$ be nice maps of connected manifolds realizing oriented bordism classes 
$\omega_{i}\in \bor_{j_i}(G_{\frak N,M}).$ Since $V_1$ is connected, all the elements $f_1(v_1), v_1\in V_1,$ are in a subspace $G_{\frak N,M}^{\Gamma_1}$ for some $\Gamma_1.$ So that $f_1$ is in fact a map to 
$G_{\frak N,M}^{\Gamma_1}\subset \sqcup_{\Gamma} G_{\frak N,M}^{\Gamma}=G_{\frak N,M}.$ Similarly $f_2$ can be regarded as a map to some $G_{\frak N,M}^{\Gamma_2}\subset G_{\frak N,M}.$ Put $\nu_i=\nu(\Gamma_i), i=1,2,$  and put 
\begin{equation}
1\omega_1\star 1\omega_2=\frac{1}{(\nu_1+\nu_2)!}\bigl (\sum_{\beta\in S_{(\nu_1+\nu_2)}} \beta_*(C(\omega_1, \omega_2))\bigr).\end{equation}
\end{defin}

The following Theorem follows immediately from the definition of $\star.$

\begin{thm}\label{multiplicationdimension}
$\star: \bor_{j_1}(G_{\frak N,M})_\Q\otimes_{\Q} \bor_{j_2}(G_{\frak N,M})_\Q\to \bor_{j_1+j_2}(G_{\frak N,M})_\Q$ is a bilinear operation.\qed
\end{thm}

\begin{thm}\label{associativecommutative}
The product $\star$ gives $\bor_*(G_{\frak N, M})_\Q$ the structure of an associative, graded commutative algebra.  Namely  for all $\eta_{i}\in \bor_{j_i}(G_{\frak N,M})_\Q, i=1,2,3,$ the following properties hold:
\begin{description}
\item[1] $\eta_1\star \eta_2=(-1)^{j_1j_2}\eta_2\star \eta_1;$ 
\item[2] $(\eta_1\star \eta_2)\star \eta_3=\eta_1\star (\eta_2\star \eta_3).$
\end{description}

\end{thm}

The proof of the Theorem~\ref{associativecommutative} is similar to the proof of the Theorem~\ref{skewsymmetry} and is omitted. An interested reader is referred to the fourth version of the arXiv preprint of this work for the detailed proof.

%\section{Leibnitz Rule for $[\cdot, \cdot]$ and $\star.$}\label{differentialalgebras}

\begin{thm}\label{Poisson}
The $\star$ product satisfies a graded Leibniz identity with respect to $[\cdot, \cdot].$ 
Namely for all $\eta_i\in \bor_{j_i}(G_{\frak N, M})_\Q, i=1,2,3,$ the following statements hold:
\begin{equation}
\begin{split}
[\eta_1, \eta_2\star \eta_3]=
(-1)^{j_3(\frak m+1)}[\eta_1, \eta_2]\star \eta_3+
(-1)^{j_2(j_1+1)}\eta_2\star [\eta_1, \eta_3].
\end{split}
\end{equation}
\end{thm}

%The following Corollary follows immediately from %Theorems~\ref{dimension},~\ref{skewsymmetry},~\ref{Jacobi},~\ref{plicationdim%e%nsion},~\ref{associativecommutative},~\ref{Poisson}.
%
%\begin{cor}
%\begin{description}
%\item[1] If $\frak N$ consists of manifolds of the same odd dimension, then %$\bigl(\bor_*(G_{\frak N, M})\otimes \Q, \star, [\cdot, \cdot]\bigr)$ is a %$\Z$-graded Poisson algebra (Gerstenhaber-like 
%algebra).
%\item[2] If $\frak N$ consists of odd-dimensional manifolds of different %dimensions, then $\bigl(\bor_*(G_{\frak N, M})\otimes \Q, \star, [\cdot, %\cdot]\bigr)$ is a $\Z_2$-graded Poisson algebra.
%\end{description}
%\end{cor}

The proof of the Theorem~\ref{Poisson} is similar to the proof of the Theorem~\ref{Jacobi} and is omitted. An interested reader is referred to the fourth version of the arXiv preprint of this work for the detailed proof.

\section{The $\bor_0(G_{\{S^1\}, F^2})\otimes \Q$ and the Andersen-Mattes-Reshetikhin Poisson algebras }\label{SectionLieLoops}

Since $\frak N=\{S^1\}$ is just one manifold and $F^2$ is fixed, we write $G$ rather than  $G_{\{S^1\}, F^2}$.
Using Theorems~\ref{dimension},~\ref{skewsymmetry},~\ref{Jacobi},~\ref{multiplicationdimension},~\ref{associativecommutative},~\ref{Poisson}  we get that $\bigl(\bor_0(G)_\Q, [\cdot, \cdot], \star\bigr)$ is a subalgebra of $\bigl(\bor_*(G)_\Q, [\cdot, \cdot], *\bigr)$ which is an ordinary Poisson algebra.

In the description of the Andersen-Mattes-Reshetikhin Poisson algebra we follow their works~\cite{AMR1},~\cite{AMR2}.

{\it A chord diagram\/} is a topological space that consists of some number of disjoint oriented circles $S^1_i, i=1, \cdots, q,$ and disjoint arcs $A_j, j=1, \cdots, r,$ such that the end points of the arcs are distinct and $\cup _j\partial A_j=(\cup _i S^1_i)\cap (\cup_j A_j).$ The circles $S^1_i$ 
are called {\it core components\/} and the arcs $A_j$ are called {\it chords\/} of the diagram.

A {\it geometric chord diagram on an oriented surface $F$\/} is a smooth map of a chord diagram to the surface, mapping each chord to a point.  {\it A generic geometric chord diagram on $F$\/} is a geometric chord diagram with all the circles immersed and all the multiple points between them being transverse double points. A {\it chord diagram on $F$\/} is an equivalence class of  geometric chord diagrams modulo homotopy.

%Note that diffeomorphism classes of chord diagrams correspond to chord diagrams on $S^2.$

The commutative {\it product of two chord diagrams on $F$\/} is defined to be their union. 

\begin{figure}[htbp] 
\begin{center} 
% \epsfxsize\hsize\advance\epsfxsize -0.5cm 
%\epsfxsize 8cm 
%\hepsffile{4Trelation.eps} 
\includegraphics[width=10cm]{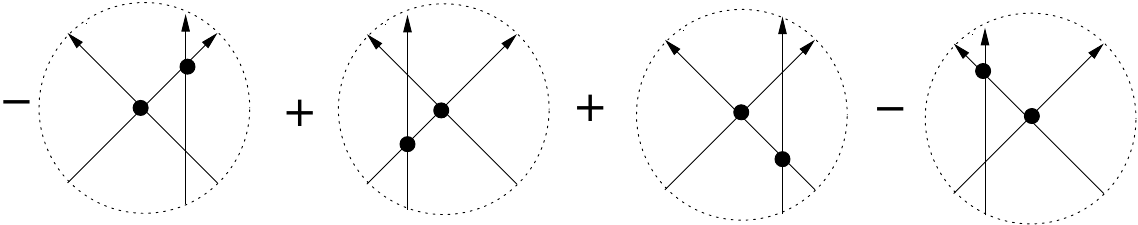} 
\end{center} 
\caption{$4T$-relation. The diagrams supposed to be identical outside of the dotted circles, and the parts that are drawn are immersed in $F.$}\label{4Trelation.fig} 
\end{figure}

Consider the complex vector space $X$ whose basis is the set of chord diagrams on $F.$
Let $Y$ be the subspace of $X$ generated by the linear combinations that are called $4T$-relations. One of the $4T$-relations is depicted in Figure~\ref{4Trelation.fig}.  
The others are obtained by reversing the orientations of strands in Figure~\ref{4Trelation.fig}
following the rule that for each chord that intersects a component whose orientation is reversed 
we get the factor of $(-1)$ in front of the diagram

The subspace $Y$ is an ideal with respect to the product, {\it and the Andersen-Mattes-Reshetikhin chord algebra is $\ch(F)=X/Y.$\/}

Given two chord diagrams $[D_1], [D_2]$ on $F,$ pick two geometric chord diagrams $D_1,D_2$ representing them so that $D_1\cup D_2$ is a generic chord diagram. Put $P=D_1\cap D_2$ to be the set of intersection points of $D_1$ and $D_2.$ Each $p\in P$ is the intersection point of a curve from 
$D_1$ with a curve  from $D_2.$
To $p\in P$ we correspond a  $2$-frame in $T_pF$ whose first vector is the velocity vector of the intersecting curve from $D_1$ and whose second vector is the velocity vector of the intersecting curve from $D_2.$ We put $\sign(p)=+1$ if this $2$-frame gives a positive orientation of $T_pF,$ and we put $\sign(p)=-1$ otherwise.
For $p\in P$ we put $D_1\cup_p D_2$ to be the chord diagram on $F$ obtained by joining $D_1^{-1}(p)$ and $D_2^{-1}(p)$ by a chord.

{\it The Andersen-Mattes-Reshetikhin Poisson bracket $\{\cdot, \cdot\}$ of two chord diagrams $[D_1], [D_2]$ on $F$\/} is given by 
\begin{equation}\label{AMRPoisson}
\{[D_1], [D_2]\}=\sum _{p\in D_1\cap D_1}\sign(p) [D_1\cup_p D_2], 
\end{equation}
where a geometric diagram in square brackets denotes the chord diagram equivalence class realized by it. 

Below we relate the Poisson subalgebra  $\bor_0(G)\otimes \C$ of our graded Poisson algebra $\bor_*(G)\otimes \C$
to the Andersen-Mattes-Reshetikhin algebra. 

Declare two points of a chord diagram to be equivalent if they belong to the same core component. A chord diagram is {\it tree-like\/} if the quotient space of it by this equivalence relation is a disjoint union of topological tree graphs. Similarly, we get the notions of {\it geometric tree-like chord diagrams on $F$\/} and of {\it tree-like chord diagrams on $F.$\/}

Consider the vector subspace $X_T\subset X$ whose basis is the set of all tree-like chord diagrams on $F.$ It is easy to see that if one of the four geometric chord diagrams in Figure~\ref{4Trelation.fig} is tree-like, then so are the other three. Put $Y_T$ to be the subspace of $X_T$ generated by the $4T$-relations consisting of tree like diagrams. Clearly $Y_T$ is an ideal of $X_T$ with respect to product.

One verifies that if the diagrams $D_1$ and $D_2$ in equation~\eqref{AMRPoisson} are tree-like, 
then all the diagrams $D_1\cup _p D_2$ are also tree-like.
Thus we get the subalgebra $\ch_T(F)$ of $\ch (F)$ that is the image of the vector subspace $X_T\subset X$ under the quotient homomorphism $X\to X/Y=\ch(F).$ Clearly $\ch_T(F)=X_T/Y_T.$

A geometric tree-like chord diagram $D$ on $F$ with all its core components enumerated in some way
defines a point $x_D\in G.$  The $S^1$-spaces in the commutative diagram $x_D$ are the enumerated circles of $D$ and the maps $S^1\to M$ in $x_D$ are the maps of the corresponding core components in the geometric diagram $D.$ The $\pt$-spaces in $x_D$ are in one-to-one correspondence with the chords of $D,$ and the image points of the $\pt\to S^1$ maps in $x_D$ are the end points of the corresponding chords. The labels of the $\pt$-spaces in $x_D$ are determined by the indices of the circles to which they are mapped.

Take a basis element $[D]\in X_T$ and a geometric chord diagram $D$ realizing $[D].$
Enumerate all the core components of $D$ in some way and let $x_D\in G$ be the resulting point. 
Let $\rho_{x_D}:\pt^+\to G$ be the map of the  positively oriented one-point space $\pt^+$ such that $\rho_{x_D}(\pt^+)=x_D.$
Put $\chi(D)=\frac{1}{|D|!}\bigl(\sum_{\alpha\in S_{|D|}}\alpha_*([\pt^+, \rho_{x_D}])\bigr)\in \bor_0(G)\otimes \C.$ Here $|D|$ is the number of core components in the diagram $D$ that equals to the number of the $S^1$-vertices in the allowed graph corresponding 
to the commutative diagram $x_D.$

Since $\chi(D)$ is constructed through the symmetrization by the action of the permutation group 
$S_{|D|},$ the element $\chi(D)$ does not depend on the enumeration of the core components of $D$ we chose to construct it. It is easy to see that if one uses the same enumeration of 
the core components in two homotopic diagrams $D_1$ and $D_2,$ then $x_{D_1}$ and $x_{D_2}$ are in the same path connected component of $G.$ Thus $\chi(D)$ is determined by the chord diagram $[D],$
and we will write $\chi([D])$ instead of $\chi(D).$

Since the basis of $X_T$ is given by the tree-like chord diagrams $[D]$ on $F,$ we have the linear homomorphism $\chi: X_T\to \bor_0(G)\otimes \C,$ defined by $[D]\to \chi([D])$ on the basis of $X_T.$

Consider a linear combination of four chord diagrams depicted in Figure~\ref{4Trelation.fig}.
If one uses the same enumeration of the core components of the four chord diagrams, then it is easy to see that the images under $\chi$ of the first and second chord diagrams are the same. Similarly the images under $\chi$ of the third and fourth chord diagrams in Figure~\ref{4Trelation.fig} are the same. Taking the signs in front of the diagrams in Figure~\ref{4Trelation.fig} into account, we see that $\chi$ maps to zero all the $4T$ linear combinations generating the ideal $Y_T\subset X_T.$

Hence we have the linear homomorphism of vector spaces $\chi:\ch_T=X_T/Y_T\to \borst_0(G)\otimes \C.$ The following proposition follows from the definitions of the operations in the Andersen-Mattes-Reshetikhin Poisson algebra and in our algebra.

\begin{prop}
$\chi$ is a Poisson algebra homomorphism of the $\ch_T(F)$ subalgebra of the Andersen-Mattes-Reshetikhin algebra $\ch(F)$ to the Poisson subalgebra $\bor_0(G)\otimes \C$ of our graded Poisson algebra $\bor_*(G)\otimes \C.$ \qed
\end{prop}

This homomorphism is not injective, since the image under $\chi$ of the first 
half of the expression depicted in Figure~\ref{4Trelation.fig} is zero, while the first half of the expression is generally nonzero in $\ch_T(F).$

\begin{remark} We do not know how to generalize the $4T$ relation to the case where the target manifold $M\neq F^2$ or to the case of higher dimensional bordism groups.
\end{remark}

\appendix
\section{Some transversality facts and homotopies between continuous and nice 
maps to $G_{\frak N, M}^{\Gamma}$}\label{transversality}
{\it In this paper the word graph stands for the allowed graph, however in this Appendix it will also stand for a subgraph of an allowed graph.\/}

\begin{thm}\label{transversalitytheorem}
Let $f:V\to G_{\frak N, M}^{\Gamma}$ be a continuous map of a smooth manifold $V.$
Then the following statements hold
\begin{enumerate}
\item $f$ is homotopic to a nice map. 
\item If $f|_{\partial V}:\partial V\to G_{\frak N, M}^{\Gamma}$ is nice, then $f$ is homotopic relative $\partial V$ to a nice map.
\item If $g:V'\to G_{\frak N, M}^{\Gamma'}$ is a nice map and 
$i\leq \nu(\Gamma), j\leq \nu(\Gamma')$ are positive integers, then $f$ is homotopic to a nice map $f'$ such that the adjoint maps $\ov f'_{i}:V\times N_i\to M$ and $\ov g_{j}:V'\times N_j\to M$ (introduced in Definition~\ref{Aoperation}) are transverse. 
\item If $g:V'\to G_{\frak N, M}^{\Gamma'}$ is a nice map, $f|_{\partial V}$ is nice, $i\leq \nu(\Gamma), j\leq \nu(\Gamma')$ are positive integers, such that $\ov f_{i}|_{\partial V\times N_{i}}:\partial V\times N_i\to M$ is transverse to $\ov g_{j}:V'\times N_j\to M,$ then $f$ is homotopic relative $\partial V$ to a nice map $f'$ such that $\ov f'_{i}:V\times N_i\to M$ and $\ov g_{j}:V\times N_j\to M$  are transverse. 
\end{enumerate}
\end{thm}

\pp
{\bf We prove statement $1$ and then explain how to modify the proof to get the other statements.\/}

Choose $\kappa\leq \nu(\Gamma).$ Put $\Gamma_{0}$ to be the $N_{\kappa}\to M$ subgraph of $\Gamma.$ (This subgraph has no $\pt$-vertices.)

For integer $l>0,$ define $\Gamma_{l}$ inductively to be the minimal subgraph of $\Gamma$ that contains: $\Gamma_{l-1}$; all the $\pt$-vertices from which there is an edge to some $N$-vertex in $\Gamma_{l-1};$ all the $N$-vertices that are the end points of the  $\pt\to N$ edges starting from such $\pt$-vertices; and all the $\pt\to N$ and $N\to M$ edges whose end points are in $\Gamma_l.$

Note that for all $l,$ the graph $\Gamma_{l}$ is a subgraph of $\Gamma_{l+1}.$
Since the graph $\Gamma$ is finite we can and do choose $\ov l\in \N$ such that $\Gamma_{l_1}=\Gamma_{l_2},$ for any $l_1, l_2\geq \ov l.$ Note that $\Gamma_{\ov l}$ is generally not equal to $\Gamma.$ 

Recall that, as it is explained in the Introduction, the graph $\Gamma$ gives rise to a unoriented graph $\ov \Gamma.$ By the definition of allowed graphs, $\ov \Gamma$ is a disjoint union of trees. It is easy to check that $\Gamma_{\ov l}$ will contain exactly those of the $\pt$-vertices, $N$-vertices, and edges between them, that are in the tree component of $\ov \Gamma$ containing the vertex $N_{\kappa}.$

For $l\geq 0,$ we define $\wt \Gamma^{l}$ to be $\Gamma_{l}$ with all the $\pt$-vertices, $N$-vertices, $\pt\to N$ edges, and $N\to M$ edges of $\Gamma,$ that are absent in $\Gamma_{\ov l},$ added to $\Gamma_{l}.$

We will construct homotopies $\mathcal H_l:[0,1]\times V\to G_{\frak N, M}^{\Gamma}, l=0, \cdots, \ov l,$ and $\mathcal H'_l:[0,1]\times V\to G_{\frak N, M}^{\Gamma}, l=1, \cdots, \ov l.$ These homotopies have the following properties 
\begin{enumerate} 
\item $\mathcal H_0(0, v)=f(v),$ for all $v\in V.$
\item $\mathcal H'_l(0, v)=\mathcal H_{l-1}(1, v)$ for all $v\in V$ and $l\in \{1, \cdots, \ov l\}.$
\item $\mathcal H_l(0, v)=\mathcal H'_{l}(1, v),$ for all $v\in V$ and $l\in \{1, \cdots, \ov l\}.$
\item for all $l\in \{1, \cdots, \ov l\}, t\in [0,1], v\in V,$ all the maps in the commutative diagram $\mathcal H'_l(t, v)$ that correspond to the $N_k\to M$ and $\pt \to N_k$ edges contained in $\wt \Gamma_{l-1}$ are exactly those from $\mathcal H_{l-1}(1, v),$ i.e. these maps do not change during the homotopy $\mathcal H'_l.$
\item for every $l\in \{1, \cdots, \ov l\},$ all the adjoint maps $V\times \pt=V\to N_k$ and $V\times N_k\to M$ arising from $\mathcal H'_l|_{1\times V}:V\to G_{\frak N, M}^{\Gamma}$ that correspond to $\pt \to N_k$ edges in $\Gamma_l$ and to $N_k\to M$ edges in $\Gamma_{l-1}$ are smooth.
\item for all $l\in \{1, \cdots, \ov l\}, t\in [0,1], v\in V,$ all the maps in the commutative diagram $\mathcal H_l(t, v)$ that correspond to $\pt \to N_k$ edges in $\wt \Gamma_{l}$ and $N_k\to M$ edges in $\wt \Gamma_{l-1}$ are exactly those from $\mathcal H'_{l}(1, v),$ i.e. these maps do not change during the homotopy $\mathcal H_l.$ Moreover for $l=0$ all the maps in the commutative diagram $\mathcal H_0(t, v)$ that correspond to $\pt \to N_k$ edges in $\wt \Gamma_{0}$ and all the $N_k\to M$ edges in $\wt \Gamma_{0},$ except of $N_{\kappa}\to M,$ are exactly those from $f(v),$ i.e. these maps do not change during the homotopy $\mathcal H_0.$
\item for every $l\in \{0, \cdots, \ov l\},$ and all $N_k\to M, \pt \to N_k$ edges in $\Gamma_l$ we have that the adjoint maps $V\times N_k\to M$ and $V\times \pt =V\to N_k$ that correspond to $\mathcal H_l|_{1\times V}:V\to G_{\frak N, M}^{\Gamma}$ are smooth.
\end{enumerate}

Put $\wt {\mathcal H}_{\kappa}:[0,1]\times V\to G_{\frak N, M}^{\Gamma}$ to be a homotopy of $f$ that is the composition of the sequence of homotopies $\mathcal H_0, \mathcal H'_1, \mathcal H_1, \cdots, \mathcal H'_{\ov l}, \mathcal H_{\ov l}.$ 
Clearly the map $\wt {\mathcal H}_{\kappa}|_{1\times V}:V\to G_{\frak N, M}^{\Gamma}$ is such that all the adjoint maps 
$V\times \pt=V\to N_k$ and $V\times N_k\to M$ that correspond to the $\pt\to N_k$  and $N_k\to M$ edges in $\Gamma_{\ov l}$ are smooth. Note also that for every $v\in V, t\in [0,1],$ the maps $N_k\to M, \pt \to M$ in the commutative diagram 
$\wt {\mathcal H}_{\kappa}(t, v)$ that correspond to the $\pt\to N_k$ edges and $N_k$ vertices in $\Gamma\setminus \Gamma_{\ov l}$ are equal to the corresponding maps in the commutative diagram $f(v).$
Put $f_{\kappa}=\wt {\mathcal H}_{\kappa}|_{1\times V}:V\to G_{\frak N, M}^{\Gamma}$ to be the map obtained in the end of the homotopy $\wt {\mathcal H}_{\kappa}$ of $f.$

Let $r$ be the number of connected components in $\ov \Gamma.$ Take one vertex $N_{\kappa_i}, i=1, \cdots, r,$ in each connected component of $\ov \Gamma.$ Apply the homotopy $\wt H_{\kappa_1}$ of $f$ to get $f_{\kappa_1}.$ For $i=2, \cdots, r$ inductively construct the homotopies 
$\wt{\mathcal H}_{\kappa_i}$ using $f_{\kappa_{i-1}}=
\wt{\mathcal H}_{\kappa_{i-1}}|_{1\times V}:V\to G_{\frak N, M}^{\Gamma}$ instead of $f$ in their construction.

Then the composition of 
the sequence of homotopies $\wt {\mathcal H}_{\kappa_1}, \wt {\mathcal H}_{\kappa_2}, \cdots, \wt {\mathcal H}_{\kappa_r}$ is the desired homotopy of $f$ to a nice map. This finishes the proof modulo the fact that we still have to construct the homotopies $\mathcal H_l$ and $\mathcal H'_l$. We do this below.

{\bf Construction of $\mathcal H_l.$\/} Fix $l\in \{0, \cdots, \ov l\}.$
Even though the enumeration of the $N$-vertices in graphs $\wt \Gamma_{a}, a=l, \cdots, \ov l$ contains gaps, it still allows one to define the space $G_{\frak N, M}^{\wt \Gamma_{a}}$ of commutative diagrams corresponding to $\wt \Gamma_{a}.$ Nice maps $V\to G_{\frak N, M}^{\wt \Gamma_{a}}$ are defined similarly to how the nice maps to $G_{\frak N, M}^{\Gamma}$ were defined.

We inductively define a sequence of homotopies $H_{l, a}:[0,1]\times V\to G_{\frak N, M}^{\wt \Gamma_{a}}, a\in \{l, \cdots, \ov l\},$ such that: 
\begin{description}
\item[1] All the adjoint maps $V\times \pt=V \to N_k$ and $V\times N_k\to M$ arising from $H_{l, l}|_{1\times V}:V\to G_{\frak N, M}^{\wt \Gamma_l}$ that correspond to $\pt\to N_k$ and $N_k\to M$ edges of $\Gamma$ contained in $\Gamma_l$ are smooth.
\item[2] If $l>0,$ then for every $t\in [0,1], v\in V$ all the maps $N_k\to M, \pt\to N_k$ in the commutative diagram $H_{l, l}(t, v)$ that correspond to the $N_k\to M$ edges in $\wt \Gamma_{l-1}$ and $\pt\to N_k$ edges in $\wt \Gamma_l$ are exactly those from the commutative diagram $\mathcal H'_{l}(1, v).$ If $l=0$ then for every $v\in V, t\in [0,1],$ all the $\pt\to N_k$ and $N_k\to M$ maps corresponding to the $\pt\to N_k$ and $N_{k}\to M$ edges in $\wt \Gamma_0,$ except of $N_{\kappa}\to M,$ are exactly those from the commutative diagram $f(v).$ 
\item[3]
For all $t\in [0,1], v\in V,$ and $a\in\{l+1,\cdots, \ov l\},$ the part of the commutative diagram $H_{l, a}(t, v)$ corresponding to the subgraph $\wt \Gamma^{a-1}\subset \wt \Gamma^{a}$ is the commutative diagram $H_{l, a-1}(t, v);$ 
\item[4] If $l>0$ then for all $v\in V$ and $a\in\{l, \cdots, \ov l\},$ the commutative diagram $H_{l, a}(0, v)$ is the part of the commutative diagram $\mathcal H'_l(1, v)$ corresponding to the subgraph $\wt \Gamma^{a}\subset \Gamma.$ If $l=0$ then for all $v\in V$ and $a\in\{0, \cdots, \ov l\},$ the commutative diagram $H_{0, a}(0, v)$ is the part of the commutative diagram $f(v)$ corresponding to the subgraph $\wt \Gamma^{a}\subset \Gamma.$ 
\end{description}

It is easy to see that we can {\bf take the homotopy $\mathcal H_l$ to be $H_{l, \ov l}.$\/}

{\bf Let us construct the homotopies $H_{l, a}.$\/}

{\it Since the homotopy $H_{0,0}$ is essential for the proof of the transversality statements of the Theorem we describe it separately.\/}
The graph $\wt \Gamma_0$ consists of the subgraph $N_{\kappa}\to M$ and of the part of $\Gamma$ that is not in $\Gamma_{\ov l}.$
Put $h:[0,1] \to C(V\times N_{\kappa}, M)$ to be a homotopy of $\ov f_{\kappa}:V\times N_{\kappa}\to M$ to a smooth map and put $\ov h:[0,1]\times V\times N_{\kappa}\to M$ to be the adjoint map. 
For all $t\in [0,1], v\in V,$ we define the commutative diagram $H_{0, 0}(t, v)\in G_{\frak N, M}^{\wt \Gamma^{0}}$ as follows:
\begin{description}
\item[a] On the part of the diagram $H_{0, 0}(t, v)$ corresponding to the subgraph $N_{\kappa}\to M$  we define the map $N_{\kappa}\to M$ in the diagram to be $\ov h|_{t\times v\times N_{\kappa}}:N_{\kappa}\to M.$
\item[b] In the part of the diagram $H_{0, 0}(t, v)$ that comes from $\Gamma\setminus \Gamma^{\ov l}$ we put all the maps 
to be the corresponding maps from the commutative diagram $f(v)\in G_{\frak N, M}^{\Gamma}.$
\end{description}

{\it Let us construct the homotopies $H_{l,l}$ for $(l,l)\neq (0,0).$ \/} We have to define all the maps in the commutative diagrams $H_{l,l}(t, v), t\in [0,1], v\in V.$
For all $t\in [0,1], v\in V$ put the maps $\pt \to N, N\to M$ in the diagram $H_{l,l}(t, v)$ that correspond to $\pt \to N$ edges in $\wt \Gamma_{l}$ and $N\to M$ edges in $\wt \Gamma_{l-1}$ to be exactly those from the commutative diagram $\mathcal H'_{l}(1,v).$

We have to define the maps $N\to M$ corresponding to $N$-vertices in $\Gamma_l\setminus \Gamma_{l-1}.$ Take such a vertex $\wh N.$  Since $\Gamma$ is an allowed graph, $\ov \Gamma$ is a union of trees. Hence there is a unique $\pt$-vertex $\check \pt$ in $\wt \Gamma_{l}$ from which there is an edge to $\wh N,$ and moreover there is a unique vertex $\check N$ in $\wt \Gamma_{l-1}$ connected to the $\check \pt$-vertex by an edge. Let $\check N$ and $\wh N$ be the manifolds corresponding to these vertices. 

For $v\in V$ put $\wh \rho(v)\in \wh N$ and $\check \rho(v)\in \check N$ to be the images of $\check \pt$ under the maps $\check \pt\to \wh N$ and $\check \pt \to \check N$ in the commutative diagram $\mathcal H'_{l}(1, v).$ By definition of $\mathcal H'_l$ the maps $V\to V\times \wh N,$ $v\to (v, \wh \rho(v)),$ and $V\to V\times \check N,$ $v\to (v, \check \rho(v))$ are smooth. Put $\wh P=\{(v, \wh \rho(v)), v\in V\}\subset V\times \wh N$ and $\check P=\{(v, \check \rho(v)), v\in V\}\subset V\times \check N.$ Then $\wh P$ and $\check P$ are graphs of smooth functions $V\to \wh N$ and $V\to \check N,$ and hence they are smooth submanifolds of $V\times \wh N$ and of $V\times \check N$ that are naturally diffeomorphic to $V.$ 

Define a continuous $\wh F:V\times \wh N\to M$ by the requirement that for all $v\in V$ the map $F|_{v\times \wh N}:\wh N\to M$ is the map corresponding to the edge $\wh N\to M$ in the diagram $\mathcal H'_l(1, v).$ By the commutativity of the diagrams, the restriction $\wh F|_{\wh P}$ equals to the composition $\wh P\to V\to \check P\to M,$ where $\wh P\to V, V\to \check P$ are diffeomorphisms and $\check P\to M$ is the restriction of the adjoint map  $V\times \check N\to M$ corresponding to the edge $\check N\to M$ that arises from $\mathcal H'_{l}|_{1\times V}:V\to G_{\frak N, M}^{\Gamma}.$  By definition of $\mathcal H'_l$ the map $V\times \check N\to M$ is smooth and hence $\wh F|_{\wh P}$ is smooth.

By Whitney Approximation on Manifolds~\cite[Theorem 10.21]{Lee} there is a homotopy  $\wh F_t:V\times \wh N\to M$ relative to $\wh P$ of $\wh F=\wh F_0$ to a smooth map $\wh F_1.$ Define the map $\wh N\to M$ in the commutative diagram $H_{l, l}(t, v)$ to be $\wh F_t|_{v\times N}.$ Perform this procedure for all the vertices $\wh N$ in $\Gamma_l\setminus \Gamma_{l-1}$ to get $H_{l, l}.$

{\it Let us construct the homotopies $H_{l,a}$ for $a>l\geq 0.$ Note that in the case $l=0$ the construction is slightly different from the one for $l>0.$\/}
Since $H_{l,a}$ are defined inductively, let $H_{l, a-1}:[0,1]\times V\to G_{\frak N, M}^{\wt \Gamma_{a-1}}$ be the already constructed homotopy. Let us construct $H_{l,a}.$

For all $t\in [0,1], v\in V,$ we define all the maps in the part of the diagram $H_{l, a}(t, v)$ corresponding to $\wt \Gamma^{a-1}$ to be the corresponding maps in the commutative diagram $H_{l, a-1}(t, v).$

Take a $\pt$-vertex $\check \pt$ in $\wt \Gamma_{a}\setminus \wt \Gamma_{a-1}.$ Since $\Gamma$ is an allowable graph, $\ov \Gamma$ is a disjoint union of trees. Then it is easy to see that there exists a unique $N$-vertex $\check N$ in $\wt \Gamma_{a-1}$ connected by an edge to $\check \pt.$ As usual, we put $\check N$ to be the corresponding manifold. 

If $l>0$ then for all $t\in[0,1], v\in V,$ we put the map $\check\rho:\pt\to \check N$ in $H_{l,a}(t,v)$ to be the corresponding map from the commutative diagram $\mathcal H'_l(1, v).$ If $l=0,$ then we put the map $\check \rho:\pt\to \check N$ in $H_{l,a}(t,v)$ to be the corresponding map from the commutative diagram $f(v).$ 

It may be that $\check \pt$-vertex is connected by an edge to some of the $N$-vertices in $\wt \Gamma_a\setminus\wt \Gamma_{a-1}.$ Take such a vertex $\wh N$ and let $\wh N$ be the corresponding manifold.

If $l>0$ then for all $t\in[0,1], v\in V,$ we put the map $\wh \rho:\pt\to \wh N$ in $H_{l,a}(t,v)$ to be the corresponding map from the commutative diagram $\mathcal H'_l(1, v).$ If $l=0,$ then we put the map $\wh\rho: \pt\to \wh N$ in $H_{l,a}(t,v)$ to be the corresponding map from the commutative diagram $f(v).$ 

Now for all $t\in [0,1], v\in V$ we have to define the map $\wh N\to M$ in the diagram $H_{l,a}(t,v).$ Put $\wh P=\{(v, \wh \rho(v))|v\in V\}\subset V\times \wh N.$ By the definition of the space $G_{\frak N, M},$ the set $\wh P$ is the graph of a continuous function $V\to \wh N, v\to \wh \rho(v),$ and hence it is naturally homeomorphic to $V.$

{\it If $l>0,$\/} then for $(t, v, \wh \rho(v))\in ([0,1]\times \wh P)\subset [0,1]\times V\times N, v\in V, t\in [0,1],$ we define $G'(t, v, \wh \rho(v))\in M$ to be the value of the map $\wh N\to M$ in the commutative diagram $\mathcal H'_l(1, v)$ computed at the point $\wh \rho(v)\in \wh N.$ For all $(0, v, \wh n)\in 0\times V\times \wh N\subset [0,1]\times V\times \wh N$ we put $G'(0, v, \wh n)\in M$ to be the value of the map $\wh N\to M$ from the commutative diagram $\mathcal H'_{l}(1, v)$ computed at the point $\wh n.$ One verifies that the map $G'$ is well-defined and continuous on $([0,1]\times \wh P)\cup (0\times V\times N)\subset [0,1]\times V\times \wh N.$ 

{\it If $l=0,$\/} then for $(t, v, \wh \rho(v))\in [0,1]\times \wh P, v\in V, t\in [0,1],$ we define 
$G'(t, v, \wh \rho(v))\in M$ to be the value of the map $\wh N\to M$ in the commutative diagram $f(v)$ computed at the point $\wh \rho(v)\in \wh N.$ For all $(0, v, \wh n)\in 0\times V\times \wh N$ we put $G'(0, v, \wh n)\in M$ to be the value of the map $\wh N\to M$ from the commutative diagram $f(v)$ computed at the point $\wh n.$ One verifies that the map $G'$ is well-defined and continuous on $([0,1]\times \wh P)\cup (0\times V\times N)\subset [0,1]\times V\times N.$

Since $[0,1]\times \wh P$ is the graph of the continuous function $[0,1]\times V\to \wh N, (v, t)\to \wh \rho(v),$ Proposition~\ref{propretract} says that $([0,1]\times \wh P)\cup (0\times V\times \wh N)$ is a retract of $[0,1]\times V\times N.$ Hence we can extend the map $G'$ to a continuous $G:[0,1]\times V\times \wh N\to M.$

For all $v\in V, t\in [0,1]$
we put the map $\wh N\to M$ in the commutative diagram $H_{l,a}(t, v)$ to be $G|_{t\times v\times \wh N}:\wh N\to M.$

Perform this construction for all vertices $\wh N$ in 
$\wt \Gamma_{l}\setminus \wt \Gamma_{l-1}$ that are connected to $\check \pt$ vertex by an edge. Then perform a similar process for all the other vertices $\check \pt$ in $\wt \Gamma_{l}\setminus \wt \Gamma_{l-1}.$ As a result we get the desired homotopy $H_{l,a}.$

This finishes the construction of homotopies $H_{l, a}$ and hence of homotopies $\mathcal H_l.$

{\bf Construction of $\mathcal H'_l.$\/} Fix $l\in \{1, \cdots, \ov l\}.$
We inductively define a sequence of homotopies $H'_{l, a}:[0,1]\times V\to G_{\frak N, M}^{\wt \Gamma_{a}}, a\in \{l, \cdots, \ov l\},$ such that: 
\begin{description}
\item[1] All the adjoint maps $V\times \pt=V \to N_k$ and $V\times N_k\to M$ arising from $H'_{l, l}|_{1\times V}:V\to G_{\frak N, M}^{\wt \Gamma_l}$ that correspond to $\pt\to N_k$ edges of $\Gamma_l$ and $N_k\to M$ edges of $\Gamma_{l-1}$ are smooth. 
\item[2] For every $t\in [0,1], v\in V,$ all the maps $N_k\to M, \pt\to N_k$ in the commutative diagram $H'_{l, l}(t, v)$ that correspond to the $N_k\to M$ and $\pt\to N_k$ edges in $\wt \Gamma_{l-1}$ are exactly those from the commutative diagram $\mathcal H_{l-1}(1, v).$ 
\item[3] For all $t\in [0,1], v\in V, a\in\{l+1,\cdots, \ov l\},$ the part of the commutative diagram $H'_{l, a}(t, v)$ corresponding to the subgraph $\wt \Gamma^{a-1}\subset \wt \Gamma^{a}$ is the commutative diagram $H'_{l, a-1}(t, v).$ 
\item[4] For all $a \in\{l, \cdots, \ov l\}, v\in V$ the commutative diagram $H'_{l, a}(0, v)$ is the part of the commutative diagram $\mathcal H_{l-1}(1, v)$ corresponding to the subgraph $\wt \Gamma^{a}\subset \Gamma.$ 
\end{description}

It is easy to see that we can {\bf take the homotopy $\mathcal H'_l$ to be $H'_{l, \ov l}.$ Let us construct the homotopies $H'_{l, a}.$\/}

{\it First we construct the homotopies $H'_{l,l}.$\/}
We have to define all the maps in the commutative diagrams $H'_{l,l}(t, v), t\in [0,1], v\in V.$
For all $t\in [0,1], v\in V,$ put the maps $\pt \to N, N\to M$ in the diagram $H'_{l,l}(t, v)$ that correspond to $\pt \to N, N\to M$ edges in $\wt \Gamma_{l-1}$ to be those from the commutative diagram $\mathcal H_{l-1}(1,v).$
We have to define the maps $\pt\to N, N\to M$ corresponding to the $\pt\to N$ and $N\to M$ edges that start at a vertex in $\Gamma_l\setminus \Gamma_{l-1}.$ 

Take a vertex $\check \pt\in \Gamma_l\setminus \Gamma_{l-1}.$ Let $\check N$ be the unique $N$-vertex in $\Gamma_{l-1}$ connected to $\check \pt$ by an edge. Put $\check N$ to be the manifold corresponding to the vertex.
Put $\check \rho:V\to \check N$ to be the map that sends $v\in V$ to the image of $\check \pt$ under the map $\check \pt\to \check N$ in the commutative diagram 
$\mathcal H'_{l-1}(1,v).$ By definition of $G_{\frak N, M}$ we get that $\check \rho$ is continuous and hence there is a homotopy $\check \rho_t, t\in [0,1],$ of $\check \rho=\check \rho_0$ to a smooth $\check \rho_1:V\to \check N.$ For $t\in [0,1], v\in V$
we put the map $\check \pt\to \check N$ in $H'_{l, l}(t,v)$ to be the map that sends $\check \pt$ to $\check \rho_t(v)\in\check N.$

It may be that the vertex $\check \pt$ is connected by an edge to some $N$-vertex $\wh N$ in $\Gamma_l\setminus \Gamma_{l-1}.$ Let $\wh N$ be the manifold corresponding to the vertex $\wh N.$
Put $\wh \rho:V\to \wh N$ to be the map that sends $v\in V$ to the image of $\check \pt$ under the map $\check \pt\to \wh  N$ in the commutative diagram 
$\mathcal H'_{l-1}(1,v).$ By definition of $G_{\frak N, M}$ we get that $\wh \rho$ is continuous and hence there is a homotopy $\wh \rho_t, t\in [0,1],$ of $\wh \rho=\wh \rho_0$ to a smooth $\wh \rho_1:V\to \wh N.$ For $t\in [0,1], v\in V$
we put the map $\check \pt\to \wh N$ in $H'_{l, l}(t,v)$ to be the map that sends $\check \pt$ to $\wh \rho_t(v)\in\wh N.$

Let us define the maps $\wh N\to M$ in $H'_{l,l}(t,v).$
Put $\check E:[0,1]\times V\times \check N\to M$ to be the continuous map such that for all $t\in [0,1], v\in V$ the map $\check E|_{t\times v\times \check N}:\check N\to M$ is the $\check N\to M$ map from $\mathcal H_{l-1}(1,v).$
Put $\check Q=\{(t, v, \check n)| \check n=\check \rho_t(v)\}\subset [0,1]\times V\times \check N$ and put $\wh Q=\{(t, v, \wh n)| \wh n=\wh \rho_t(v)\}\subset [0,1]\times V\times \wh N.$ Since $\check Q$ and $\wh Q$ are graphs of continuous functions $V\times [0,1]\to \check N$ and $V\times [0,1]\to \wh N$ respectively, they are homeomorphic to $V\times [0,1]$ and the map 
$\wh Q\to \check Q, (t, v, \wh \rho_t(v))\to (t, v, \check \rho_t(v))$ is a homeomorphism.

We 
define the map $\wh E'$ of $(0\times V\times \wh N)\cup \wh Q\subset [0,1]\times V\times \wh N$ to $M,$ as follows.
We put $\wh E'|_{0\times V\times \wh N}:V\times \wh N \to M$ to be the adjoint map corresponding to the edge $\wh N\to M$ arising from $\mathcal H_{l-1}|_{1\times V}.$ We define $\wh E'|_{\wh Q}$ to be the composition of the homeomorphism $\wh Q\to \check Q, (t, v, \wh \rho_t(v))\to (t, v, \check \rho_t(v))$ and the continuous $\check E.$ From the commutativity of the diagrams $\mathcal H_{l-1}(1,v)$ we get that $\wh E'$ is continuous. 

Since $\wh Q$ is the graph of a continuous function $V\times [0,1]\to \wh N,$ we get that $(0\times V\times \wh N)\cup \wh Q$ is a retract of $[0,1]\times V\times \wh N$ by Proposition~\ref{propretract}.
Hence we can extend $\wh E'$ to a continuous $\wh E:[0,1]\times V\times \wh N\to M.$ Finally for all $t\in [0,1], v\in V$ we put the $\wh N\to M$ map in the diagram $H'_{l,l}(t,v)$ to be 
$\wh E|_{t\times v\times \wh N}:\wh N\to M.$ 

Perform this procedure for all vertices $\wh N$ in $\Gamma_l\setminus \Gamma_{l-1}$ that are connected to $\check \pt$ by an edge. Repeat the whole process for all the $\check \pt$ vertices in $\Gamma_l\setminus \Gamma_{l-1}$ to get the homotopy $H'_{l,l}.$

{\it Homotopies $H'_{l,a}, a>l.$\/} Since $a>l,$ we assume that we have constructed $H'_{l,a-1}.$ To construct $H'_{l,a}:[0,1]\times V\to G_{\frak N, M}^{\wt \Gamma_a}$ we have to define all the $\pt\to N$ and $N\to M$ maps in the diagrams $H'_{l,a}(t,v), t\in [0,1], v\in V.$

For all $t\in [0,1], v\in V,$ put the maps $\pt \to N, N\to M$ in the diagram $H'_{l,a}(t, v)$ that correspond to $\pt \to N, N\to M$ edges in $\wt \Gamma_{a-1}$ to be those from the commutative diagram $\mathcal H'_{l,a-1}(t,v).$
We have to define the maps $\pt\to N, N\to M$ corresponding to the $\pt\to N$ and $N\to M$ edges that start at a vertex in $\Gamma_a\setminus \Gamma_{a-1}.$ 

Take a vertex $\check \pt\in \Gamma_a\setminus \Gamma_{a-1}.$ Let $\check N$ be the unique $N$-vertex in $\Gamma_{a-1}$ connected to $\check \pt$ by an edge. Put $\check N$ to be the manifold corresponding to the vertex.
Put $\check \rho: V\to \check N$ to be the continuous map that sends $v\in V$ to the image of $\check \pt$ under the map $\check \pt\to \check N$ in the commutative diagram $\mathcal H_{l-1}(1,v).$ 
For $t\in [0,1], v\in V$
we put the map $\check \pt\to \check N$ in $H'_{l, a}(t,v)$ to be the map that sends $\check \pt$ to $\check \rho (v)\in\check N.$

It may be that the vertex $\check \pt$ is connected by an edge to some $N$-vertex $\wh N$ in $\Gamma_a\setminus \Gamma_{a-1}.$ Let $\wh N$ be the manifold corresponding to the vertex $\wh N.$
Put $\wh \rho:V\to \wh N$ to be the continuous map that sends $v\in V$ to the image of $\check \pt$ under the map $\check \pt\to \wh  N$ in the commutative diagram 
$\mathcal H_{l-1}(1,v).$  For $t\in [0,1], v\in V$
we put the map $\check \pt\to \wh N$ in $H'_{l, a}(t,v)$ to be the map that sends $\check \pt$ to $\wh \rho (v)\in\wh N.$

Let us define the maps $\wh N\to M$ in $H'_{l,a}(t,v).$
Put $\check F:[0,1]\times V\times \check N\to M$ to be the continuous map such that for all $t\in [0,1], v\in V$ the map $\check F|_{t\times v\times \check N}:\check N\to M$ is the $\check N\to M$ map from $H'_{l, a-1}(t,v).$
Put $\check P=\{(v, \check n)| \check n=\check \rho (v)\}\subset V\times \check N$ and put $\wh P=\{(v, \wh n)| \wh n=\wh \rho (v)\}\subset V\times \wh N.$ 
Clearly $\check P$ and $\wh P$ are the graphs of continuous maps $V\to \check N$ and $V\to \wh N$ respectively. Hence $\check P$ and $\wh P$ are homeomorphic to $V$ and $\wh P\to \check P, (v, \wh \rho(v))\to (v, \check \rho(v))$ is a homeomorphism.

We 
define the map $\wh F'$ of $(0\times V\times \wh N)\cup ([0,1]\times \wh P)\subset [0,1]\times V\times \wh N$ to $M,$ as follows.
We put $\wh F'|_{0\times V\times \wh N}:V\times \wh N \to M$ to be the adjoint map corresponding to the edge $\wh N\to M$ arising from $\mathcal H_{l-1}|_{1\times V}.$ We define $\wh F'|_{[0,1]\times \wh P}$ to be the composition of the homeomorphism $[0,1]\times \wh P\to [0,1]\times \check P, (t, v, \wh \rho (v))\to (t, v, \check \rho (v))$ and the continuous $\check F.$ From the commutativity of the diagrams $\mathcal H_{l-1}(1,v)$ we get that $\wh F'$ is continuous. 

Since $[0,1]\times \wh P$ is identified with the graph of the continuous function $[0,1]\times V\to \wh N, (t, v)\to \wh \rho(v),$ we have that  $(0\times V\times \wh N)\cup ([0,1]\times \wh P)$ is a retract of $[0,1]\times V\times \wh N$ by Proposition~\ref{propretract}.
Thus we can extend $\wh F'$ to a continuous $\wh F:[0,1]\times V\times \wh N\to M.$ Finally for all $t\in [0,1], v\in V,$ we put the $\wh N\to M$ map in the diagram $H'_{l,a}(t,v)$ to be $\wh F|_{t\times v\times \wh N}:\wh N\to M.$ 

Perform this procedure for all vertices $\wh N$ in $\Gamma_a\setminus \Gamma_{a-1}$ that are connected to $\check \pt$ by an edge. Repeat the whole process for all the $\check \pt$ vertices in $\Gamma_a\setminus \Gamma_{a-1}$ to get the homotopy $H'_{l,a}.$

This finishes the proof of statement $1$ of Theorem~\ref{transversalitytheorem}. 

{\bf To prove statement $2$\/} one has to require that the homotopies $\mathcal H_l, \mathcal H'_l, H_{l,a}, H'_{l,a}$ are homotopies relative $\partial V.$
The construction of the homotopies is almost the same as the one we used in the proof of statement $1.$ The difference is that one uses different  retracts of $[0,1]\times V\times N$ that are obtained by adding $[0,1]\times \partial V\times N$ to the retracts used in the proof of statement $1.$

{\bf To prove statement $3$\/} one does the following. Choose the vertex $N_{\kappa}$ from which we started the construction in the proof of statement $1$ to be $N_i$ from statement $3.$ One has to modify the homotopy $\mathcal H_0$ so that the adjoint map $V\times N_i\to M$ arising from $\mathcal H_0|_{1\times V}$ that corresponds to the edge $N_{\kappa}\to M$ is transverse to $\ov g_j:V'\times N_j\to M.$

This is achieved as follows. By Whitney Approximation on Manifolds~\cite[Theorem 10.21]{Lee} there is a homotopy $h_1:[0,1]\to C(V\times N_i, M)$ of $\ov f_i$ to a smooth map $\ov f''_i.$
The result of Biasi and Saeki~\cite[Theorem 2.4, Remark 2.7]{BiasiSaeki}
implies that there is a homotopy $h_2:[0,1]\to C^{\infty}(V\times N_i, M)$ of $\ov f''_i$ to a $C^{\infty}$-close map transverse to $\ov g_{j}.$ Moreover their result says that the space of smooth maps $V\times N_{i}\to M$ that are transversal to $\ov g_{j}$ is residual in $C^{\infty}(V\times N_{i}, M).$ 
Define a homotopy $h:[0,1]\to C(V\times N_i, M)$ via $h(t)=h_1(2t)$ for $t\in [0,\frac{1}{2}];$ $h(t)=h_2(1-2t)$ for $t\in [\frac{1}{2}, 1].$ To get the proof one uses this homotopy $h$ as the homotopy $h$ in the construction of $H_{0,0}.$

{\bf To get statement $4$\/} One combines the ideas of the proofs of statements $2$ and $3.$ Namely one requires that all the homotopies $\mathcal H_l, \mathcal H'_l, H_{l,a}, H'_{l,a}$ are homotopies relative boundary. To construct them one uses different  retracts of $[0,1]\times V\times N$ that are obtained by adding $[0,1]\times \partial V\times N$ to the retracts used in the proof of statement $1.$

To get the desired transversality condition one starts the construction of the homotopies from $N_{\kappa}=N_i.$ The homotopy $\mathcal H_0$ should be such  that the adjoint map $V\times N_i\to M$ arising from $\mathcal H_0|_{1\times V}$ that corresponds to the edge $N_{\kappa}\to M$ is transverse to $\ov g_j:V'\times N_j\to M.$

This is achieved as follows. By Whitney Approximation on Manifolds~\cite[Theorem 10.21]{Lee} there is a homotopy $h_1:[0,1]\to C^0(V\times N_i, M)$ relative $\partial V \times N_i$ of $\ov f_i$ to a smooth map $\ov f''_i.$ Since $\ov f_i:\partial V\times N_i\to M$ is transverse to $\ov g_j:V'\times N_j\to M$ and $h_1$ is a homotopy relative boundary, we get that 
$\ov f''_i:\partial V\times N_i\to M$ is transverse to $\ov g_j:V'\times N_j\to M.$ 
A result similar to the one of Biasi and Saeki~\cite[Theorem 2.4, Remark 2.7]{BiasiSaeki}
implies that there is a homotopy $h_2:[0,1]\to C^{\infty}(V\times N_i, M)$ of 
$f''_i$ to a $C^{\infty}$-close map transverse to $\ov g_{j}$ such that for all $t\in [0,1]$ we have $h_2(t)|_{\partial V\times N_i}=h_2(0)|_{\partial V\times N_i}=\ov f''_i:\partial V\times N_i\to M.$
Now, similarly to the case of statement $3,$ one uses $h_1$ and $h_2$ to construct the homotopy $h.$ Then one uses this $h$ as the homotopy $h$ in the construction of $H_{0,0}.$ 

\qed

\begin{prop}\label{propretract}
Let $V, N$ be compact smooth connected manifolds, and let $\phi:[0,1]\times V\to N$ be a continuous map. 
Let $\Gamma(\phi)=\{(t, v, \phi(t, v))|t\in [0,1], v\in V\}\subset [0,1]\times V\times N$ be the graph of $\phi.$ Then
\begin{enumerate}
\item The subspace $\bigl((0\times V\times N)\cup \Gamma(\phi)\bigr)\subset [0,1]\times V\times N$ is a retract of $[0,1]\times V\times N$.
\item The subspace $\bigl((0\times V\times N)\cup ([0,1]\times \partial V\times N)\cup \Gamma(\phi)\bigr)\subset [0,1]\times V\times N$ is a retract of $[0,1]\times V\times N$.
\end{enumerate}
\end{prop}

\pp We prove statement $1$. Statement $2$ is proved similarly.

Let $\mathcal M$ be the class of metrizable spaces. By~\cite[page 98]{Borsuk} every manifold is an $\ANR$ for $\mathcal M.$ Hence $0\times V\times N,$ $[0,1]\times V\times N$ and $\Gamma(\phi),$ that is homeomorphic to $[0,1]\times V,$ are $\ANR$. Since $(0\times V\times N)\cap \Gamma(\phi)$ is  homeomorphic to the graph of $\phi|_{0\times V}:V\to N,$ it is a manifold and hence an $\ANR.$

By~\cite[page 33]{Borsuk} the class $\mathcal M$ is weakly hereditary. Thus by~\cite[Theorem 3.1 and Theorem 3.2, pages 83-84]{Borsuk} a space $Y$ is an $\ANE$ for the class $\mathcal M$ if and only if it is an $\ANR$ for the class $\mathcal M.$ Now since $(0\times V\times N), \Gamma(\phi)$ are closed subspaces of $(0\times V\times N)\cup \Gamma(\phi),$ and the spaces $(0\times V\times N),$ $\Gamma(\phi),$ $(0\times V\times N)\cap \Gamma(\phi)$ are $\ANR,$ then  $(0\times V\times N)\cup \Gamma(\phi)$ is an $\ANR,$ by~\cite[Proposition 10.1]{Borsuk}.

Clearly $0\times V\times N$ is a strong deformation retract of both $[0,1]\times V\times N$ and $(0\times V\times N)\cup \Gamma(\phi).$ Thus the inclusion $i:(0\times V\times N)\cup \Gamma(\phi)\to [0,1]\times V\times N$ induces an isomorphism $i_*:\pi_n\bigl((0\times V\times N)\cup \Gamma(\phi)\bigr)\to \pi_n \bigl([0,1]\times V\times N\bigr)$ for all $n.$ Since 
$(0\times V\times N)\cup \Gamma(\phi)$ is a closed connected $\ANR$ subspace of a connected $\ANR$ space $[0,1]\times V\times N,$ we get that $(0\times V\times N)\cup \Gamma(\phi)$ is a deformation retract of 
$[0,1]\times V\times N,$ by~\cite[Theorem 8.2, page 218]{Borsuk}.
\qed.

\begin{remark}[Bordisms and nice maps]\label{nicebordismmaps}
Statements $1$ and $2$ of Theorem~\ref{transversalitytheorem} imply that every bordism class in $\bor_*(G_{\frak N, M}^{\Gamma})$ is realizable by a nice map and that moreover if two nice maps are bordant then the bordism between them can be chosen to be a nice map. 
%In particular this shows that %the special bordism group of $G_{\frak N, M}$ that we %used in the first version of this text is the same as %the standard bordism group of $G_{\frak N, M}$ 
\end{remark}

\begin{remark}[satisfying multiple transversality conditions]\label{pletransversalityconditions}
If we are given $n\geq 2$ nice maps $g_r:V_r^{j_r}\to G_{\frak N, M}^{\Gamma_r}, r=1,\cdots, n$ and $n$ positive integers $j_r\leq \nu(\Gamma_r), r=1, \cdots n,$ 
then in the statements $3$ and $4$ of Theorem~\ref{transversalitytheorem} $f$ is homotopic to a nice map $f'$ such that $\ov f'_i:V\times N_i\to M$ is transverse to all the adjoint maps $\ov g_{r, j_r}:V_r\times N_{j_r}\to M.$ This is since, as it was shown by Biasi and Saeki~\cite{BiasiSaeki}, the space of smooth maps $V\times N_{i}\to M$ that are transversal to each $\ov g_{r, j_r}$ is residual in $C^{\infty}(V\times N_{i}, M).$ 
\end{remark}

\section{Algebras in the case where $\frak N$ consists of manifolds of the same even dimension $\frak n.$}\label{EvenDimensional}.

\begin{thm}\label{main2}
For $\frak N$ consisting of even-dimensional manifolds of the same dimension $\frak n$ and $\eta_i\in \bor_{j_i}(G_{\frak N, M})_ \Q, i=1,2,3,$ the operations $\{\cdot, \cdot\}$ and $\star$ satisfy the
following identities and have the following properties 
\begin{description}
\item[1]
graded commutativity 
$$\{\eta_1, \eta_2\}=(-1)^{(|\eta_1|-\frak m)(|\eta_2|-\frak m)}\{\eta_2, \eta_1\};$$
\item[2] graded Leibnitz identity 
$$
\{\eta_1, \eta_2\star\eta_3\}=\{\eta_1, \eta_2\}\star \eta_3+(-1)^{|\eta_2|(|\eta_1|-\frak m)}\eta_2\star\{\eta_1, \eta_3\}.
$$ 
\end{description}
\end{thm}

\begin{remark}
It is not hard to see that when $\frak n$ is even, the operation $\{\cdot, \cdot\}$ does not satisfy any graded Jacobi identity. For this, consider the elements of the form $1\omega_i$, where $\omega_i$ $i=1,2,3$ are $0$-dimensional bordism classes corresponding to the trivial (with no $\pt$-vertices) commutative diagrams of the form $\phi_i:N_i^{\frak n}\to M^{2\frak n}$. Assume moreover that $\phi_i$ and $\phi_j$, $ i\neq j$ have exactly one (transverse) intersection point.
\end{remark}

\section{Proof of Theorem~\ref{skewsymmetryingredient}}\label{proofskewsymmetryingredient}

\pp 
For $i=1,2$ put $N_{k_i}\in \frak N$ to be the manifold indexed by $k_i$
in the commutative diagrams from $G_{\frak N, M}^{\Gamma_i}.$ Using Theorem~\ref{transversalitytheorem} we choose nice maps $f_i: V_i^{j_i}\to G_{\frak N,M}^{\Gamma_i}$ realizing $\omega_i, i=1, 2$ so that the adjoint maps $\ov f_{1, k_1}:V_1\times N_{k_1}\to M$ and $\ov f_{2,k_2}:V_2\times N_{k_2}\to M$ (introduced in Section~\ref{Aoperation}) are transverse.

Put $\frak L=V_1\times N_{k_1}\times V_2\times N_{k_2}$ and $\frak R=V_2\times N_{k_2}\times V_1\times N_{k_1}.$
Put $L\subset \frak L$ to be the oriented $(j_1+j_2+2\frak n-\frak m)$-dimensional manifold $W$ in the definition of $A_{k_1, k_2}(\omega_1, \omega_2).$ Clearly 
\begin{equation}\label{Lskeweq}
L=\{(v_1, n_{k_1}, v_2, n_{k_2})\in \frak L| \ov f_{1, k_1}(v_1, n_{k_1})=\ov f_{2, k_2}(v_2, n_{k_2})\}.
\end{equation}
Put $g_L:L\to G_{\frak N, M}^{B_{k_1, k_2}(\Gamma_1, \Gamma_2)}$ to be the nice map constructed according to the
definition of the $A_{k_1, k_2}$-operation, so that 
\begin{equation}[L, g_L]=A_{k_1, k_2}(\omega_1, \omega_2)\in \bor_{(j_1+j_2+2\frak n-\frak m)}\bigl(G_{\frak N, M}^{B_{k_1, k_2}(\Gamma_1, \Gamma_2)}\bigr).
\end{equation}

Put $R\subset \frak R$ to be the oriented $(j_1+j_2+2\frak n-\frak m)$-dimensional manifold $W$ in the definition of $A_{k_2, k_1}(\omega_2, \omega_1).$ Clearly 
\begin{equation}\label{Rskeweq}
R=\{(v_2, n_{k_2}, v_1, n_{k_1})\in \frak R| \ov f_{2, k_2}(v_2, n_{k_2})=\ov f_{1, k_1}(v_1, n_{k_1})\}.
\end{equation}
Let $g_R:L\to R_{\frak N, M}^{B_{k_2, k_1}(\Gamma_2, \Gamma_1)}$ be the nice  map constructed according to the 
definition of $A_{k_2, k_1},$ so that $[R, g_R]=A_{k_2, k_1}(\omega_2, \omega_1)\in \bor_{(j_1+j_2+2\frak n-\frak m)}\bigl(G_{\frak N, M}^{B_{k_2, k_1}(\Gamma_2, \Gamma_1)}\bigr).$

Let $\Theta:\frak L=V_1\times N_{k_1}\times V_2\times N_{k_2}\to \frak R=V_2\times N_{k_2}\times V_1\times N_{k_1}$ be the factor permuting diffeomorphism defined by $\Theta(v_1, n_{k_1}, v_2, n_{k_2})=(v_2, n_{k_2}, v_1, n_{k_1}),$ for all $v_1\in V_1, n_{k_1}\in N_{k_1}, v_2\in V_2, n_{k_2}\in N_{k_2}.$
Comparing equations~\eqref{Lskeweq} and~\eqref{Rskeweq} that define $L$ and $R,$ we see that $\Theta|_{L}:L\to R$ is a diffeomorphism.

Take $l=(v_1, n_{k_1}, v_2, n_{k_2})\in L$ and $r=\Theta(l)=(v_2, n_{k_2}, v_1, n_{k_1})\in R.$
From the definition of $A_{k_1, k_2}$ it is clear that the commutative diagram $g_L(l)$
is obtained by the following sequence of operations:
\begin{description}
\item[1]  take the disjoint union $f_1(v_1)\sqcup f_2(v_2)$ of commutative diagrams;
\item[2] add the extra $\pt$-space mapped to $n_{k_1}\in N_{k_1}$ and to $n_{k_2}\in N_{k_2}$; 
\item[3] identify the copies of $M$ in the commutative diagrams $f_1(v_2), f_2(v_2);$ and
\item[4] change the indices of the $N$-manifolds and the labels of the $\pt$-spaces.
\end{description}
Similarly, the commutative diagram $g_R(r)=g_R(\Theta(l))$
is obtained by the following sequence of operations
\begin{description}
\item[1] take the  disjoint union $f_2(v_2)\sqcup f_1(v_1)$ of commutative diagrams; 
\item[2] add the extra $\pt$-space mapped to $n_{k_2}\in N_{k_2}$ and to $n_{k_1}\in N_{k_1}$; 
\item[3] identify the copies of $M$ in the commutative diagrams $f_1(v_2), f_2(v_2);$ and
\item[4] change the indices of the $N$-manifolds and the labels of the $\pt$-spaces.  
\end{description}
So, up the change of the indices of $N$-manifolds and of the labels of the $\pt$-spaces, the commutative diagrams $g_L(l)$ and $g_R(\Theta(l))$ are the same.
Using Statement $2$ of Proposition~\ref{propgraph} and Definition~\ref{groupaction}, we see that $\wt \alpha(g_R(\Theta(l)))=g_L(l).$
Thus we have
$[L, g_L]=[L, \wt \alpha\circ g_R\circ \Theta]= \alpha_*([L, g_R\circ \Theta])=\pm 1\alpha_*([R, g_R])\in \bor_{j_1+j_2+2\frak n-\frak m}\bigl(G_{\frak N, M}^{A_{k_1, k_2}(\Gamma_1, \Gamma_2}\bigr).
$

The manifolds $L$ and $R$ are oriented as it is described in the definition of $A_{j_1, j_2}^{k_1, k_2},$ and the statement of the Theorem would follow immediately if we show that 
$(\Theta|_{L})_*:TL\to (-1)^{(j_1+\frak n_{k_1})(j_2+\frak n_{k_2})+\frak m} TR$ is orientation preserving. 
{\it (Here and below a product of a sign and an oriented bundle denotes the initial oriented bundle if the sign is $+1,$ and it denotes the initial bundle with the changed orientation if the sign is $-1.$)\/}
Let us show that $(\Theta|_{L})_*$ is indeed orientation preserving. 

{\it In the proof we will use the following notations.\/} For an oriented space $T_xX,$ {\it we denote by $\Or(T_xX)$\/} its orientation. 
For $k>0,$ a {\it $k$-frame in $T_xX$\/} is an ordered sequence of $k$-vectors in $T_xX,$ and we call $k$ {\it the size of the frame.\/} {\it A $0$-frame in $T_xX$\/} is a number $\pm 1.$ {\it An orientation frame in $T_xX$\/} is a frame that is a basis of $T_xX.$
For an orientation frame $\Fr\subset T_xX,$ {\it we denote by $\Or(\Fr)$\/} the orientation of $T_xX$ given by the frame. 
(A $0$-frame is a number $\pm 1.$ So if $X$ is $0$-dimensional and $\Fr\subset T_xX$ is an orientation $0$-frame, we can still write $\Or(\Fr)$ for the orientation of $T_xX$ given by the $0$-frame $\Fr.$)  For an oriented space $T_xX$ we write $\Fr^+(T_xX)$ for a positive orientation frame of $T_xX.$

For a sequence $\Fr_1, \cdots, \Fr_j\subset T_xX$ of frames possibly of different sizes, {\it we denote by $\{\Fr_1, \cdots, \Fr_j\}$\/} the ordered sequence of vectors in $T_xX$ and signs that is obtained by concatenating the frames. That is to get $\{\Fr_1, \cdots, \Fr_j\}$ take the ordered sequence of vectors in $\Fr_1$ (or the corresponding sign $\pm 1$ if $\Fr_1$ is a $0$-frame); followed by the ordered sequence of vectors in $\Fr_2$ (or the corresponding sign $\pm 1$ if $\Fr_2$ is a $0$-frame), etc; followed by the ordered sequence of vectors in $\Fr_j$
(or by the corresponding sign $\pm 1$ if $\Fr_j$ is a $0$-frame). By a slight abuse of notation we 
would sometimes call $\{\Fr_1, \cdots, \Fr_j\}$ {\it a frame in $T_xX.$\/} We call such a frame 
a {\it $k$-frame (or a frame of size $k$)\/} if the number of vectors in it is $k.$

If the union of all the vectors in $\{\Fr_1, \cdots, \Fr_j\}$ is a basis of $T_xX,$ then {\it we write 
$\Or(\{\Fr_1, \cdots, \Fr_j\})$\/} for the orientation of $T_xX$ that is the product of all the signs $\pm 1$ in $\{\Fr_1, \cdots, \Fr_j\}$ and of the orientation of $T_xX$ given by the frame obtained by deleting all the signs from $\{\Fr_1, \cdots, \Fr_j\}$ and keeping the order of the vectors intact.

Take $l=(v_1, n_{k_1}, v_2, n_{k_2})\in L.$ 
Put $m=\ov f_{1, k_1}(v_1, n_{k_1})=\ov f_{2, k_2}(v_2, n_{k_2})\in M.$ 
Let $\{ \vec {m_1}, \cdots, \vec {m_{\frak m}}\}\subset T_mM$ be an orthogonal positive frame, then $\{ \vec {m_1}\oplus \vec {m_1}, \cdots, \vec {m_{\frak m}}\oplus \vec {m_{\frak m}}\}$ is an orthogonal positive frame of $T_{m\times m}\Delta$  and $\Fr_{m\times m}^{\perp}=\{ -\vec {m_1}\oplus \vec {m_1}, \cdots, -\vec {m_{\frak m}}\oplus \vec {m_{\frak m}}\}$ is an orthogonal  positive $\frak m$-frame in the fiber of the normal bundle $(T\Delta)^{\perp}$ over the point $m\times m.$

Recall that  $\pr:T(M\times M)|_{\Delta}=T\Delta\oplus (T\Delta)^{\perp}\to (T\Delta)^{\perp}$ is the orthogonal projection. Put $\Fr_L^{\perp}\subset T_l\frak L$ to be the unique $\frak m$-frame orthogonal to $T_lL$ such that 
$\pr \bigl((\ov f_{1, k_1}\times \ov f_{2, k_2})_*(\Fr_L^{\perp})\bigr)=\Fr_{m\times m}^{\perp}.$
Choose a positive $(j_1+j_2+2\frak n-\frak m)$-frame $\Fr_L$ in $T_lL.$ By our orientation convention 
\begin{equation}\label{skeweq1}
\Or(\{\Fr_L, \Fr_L^{\perp}\})=\Or(T_l\frak L).
\end{equation}

Put $r=\Theta(l)=(v_2, n_{k_2}, v_1, n_{k_1})$ and put $\Fr R^{\perp}\subset T_r\frak R$ to be the unique $\frak m$-frame orthogonal to $T_rR$ such that 
$\pr \bigl((\ov f_{2, k_2}\times \ov f_{1, k_1})_*(\Fr_R^{\perp})\bigr)=\Fr_{m\times m}^{\perp}.$
Choose a positive $(j_1+j_2+2\frak n-\frak m)$-frame $\Fr_R$ in $T_rR.$ We have \begin{equation}\label{skeweq2}
\Or(\{\Fr_R, \Fr_R^{\perp}\})=\Or(T_r\frak R).
\end{equation}

The vector spaces $T_l\frak L$ and $T_r\frak R$ are naturally decomposed as the direct sums of the linear subspaces $T_{v_1}V_1, T_{v_2}V_2, T_{n_{k_1}}N_{k_1}, T_{n_{k_2}}N_{k_2}.$ 
For a vector $\vec w$ in $T_l\frak L$ or in $T_r\frak R$ we will denote by $\vec w_{V_1}, \vec w_{V_2}, \vec w_{N_{k_1}}, \vec w_{N_{k_2}}$
the components of $\vec w$ in the corresponding linear subspaces.

Let $\vec w^i$  be the $i$-th vector in $\Fr^{\perp}_{L}.$ By definition of 
$\Fr^{\perp}_{L}$ we have $\pr\bigl((\ov f_{1, k_1}\times \ov f_{2, k_2})_*(\vec w^i_{V_1}, \vec w^i_{N_{k_1}}, \vec w^i_{V_2}, \vec w^i_{N_{k_2}})\bigr)= -\vec {m_i}\oplus \vec {m_i}.$ Since $\Theta$ exchanges the factors we have $\Theta_*(\vec w^i)_{V_1}=\vec w^i_{V_1}, \Theta_*(\vec w^i)_{N_{k_1}}=\vec w^i_{N_{k_1}},$ $\Theta_*(\vec w^i)_{V_2}=\vec w^i_{V_2}, \Theta_*(\vec w^i)_{N_{k_2}}=\vec w^i_{N_{k_2}}.$
Thus $\pr\bigl((\ov f_{2, k_2}\times \ov f_{1, k_1})_*(\Theta_*(\vec w^i)_{V_2}, \Theta_*(\vec w^i)_{N_{k_2}}, \Theta_*(\vec w^i)_{V_1}, \Theta_*(\vec w^i)_{N_{k_1}})\bigr)= \vec {m_i}\oplus -\vec {m_i}.$

Let $\vec u^i$ be the $i$-th vector in $\Fr_{R}^{\perp}.$
By definition of $\Fr_{R}^{\perp}$ we get $\pr\bigl((\ov f_{2, k_2}\times \ov f_{1, k_1})_*(\vec u^i_{V_2}, \vec u^i_{N_{k_2}}, \vec u^i_{V_1}, \vec u^i_{N_{k_1}})\bigr)=-\vec {m_i}\oplus \vec {m_i}.$ 

Thus $\pr\bigl((\ov f_{1, k_1}\times \ov f_{2, k_2})_*(-\vec u^i_{V_1}, -\vec u^i_{N_{k_1}}, -\vec u^i_{V_2}, -\vec u^i_{N_{k_2}})\bigr)=\vec {m_i}\oplus -\vec {m_i}.$ Since $\Theta$ preserves the Riemannian metric, $\Theta_*(\vec w^i)$ is orthogonal to $T_rR$ and we have $\Theta_*(\vec w^i)=-\vec u^i$ for all $i=1, \cdots, \frak m.$ Hence 
\begin{equation}\label{skeweq3}
\Theta_*( \Fr^{\perp}_{L})=-\Fr^{\perp}_R, 
\end{equation}
where $-\Fr^{\perp}_R$ denotes the frame obtained by multiplying each vector in $\Fr^{\perp}_R$ by $(-1)$ and keeping the order of the vectors in the frame intact.

Since $\Theta$ is defined by permuting the factors we get
\begin{equation}\label{skeweq4}
\Or\bigl(\Theta_*(\Fr^+(T_l\frak L))\bigr)=(-1)^{(j_1+\frak n)(j_2+\frak n)}\Or(T_r\frak R).
\end{equation}
Equations~\eqref{skeweq1},~\eqref{skeweq3},~\eqref{skeweq4} imply
\begin{equation}\label{skeweq5}
\Or(\{\Theta_*(\Fr_L), \Theta_*(\Fr_L^{\perp})\})=\Or(\{\Theta_*(\Fr_L), -\Fr_R^{\perp}\}) =(-1)^{(j_1+\frak n)(j_2+\frak n)}\Or(T_r\frak R).
\end{equation}
 Since  $\Fr_R^{\perp}$ is an $\frak m$-frame, equation~\eqref{skeweq5} implies
\begin{equation}\label{skeweq6}
\Or(\{\Theta_*(\Fr_L), \Fr_R^{\perp}\})=(-1)^{(j_1+\frak n)(j_2+\frak n)+\frak m}\Or(T_r\frak R).
\end{equation}

Equations~\eqref{skeweq2},~\eqref{skeweq6} identities $\Or(\Fr_R)=\Or(T_rR), \Or(\Fr_L)=\Or(T_lL),$
and the fact that $\Theta_*(\Fr_L)$ is an orientation frame in $T_rR$ imply that
$\Or(\Theta_*(\Fr_L))=(-1)^{(j_1+\frak n)(j_2+\frak n)+\frak m}\Or(T_rR).$ \qed

\section{Proof of the Theorem~\ref{Jacobiingredient}}\label{proofJacobiingredient}

\pp Throughout this proof $N_{k_i}\in \frak N$ denotes the manifold indexed by $k_i$ in the commutative diagrams from $G_{\frak N, M}^{\Gamma_i}, i=1,2,3;$ and $N_{\ov k_2}\in \frak N$ denotes the manifold indexed by $\ov k_2$ 
in the commutative diagrams from $G_{\frak N, M}^{\Gamma_2}.$

Choose nice maps  $f_i: V_i^{j_i}\to G_{\frak N,M}^{\Gamma_i}$ realizing $\omega_i, i=1, 2,3.$ Using Theorem~\ref{transversalitytheorem} deform 
$f_2$ so that the adjoint maps $\ov f_{1, k_1}:V_1\times N_{k_1}\to M$ and $\ov f_{2, k_2}:V_2\times N_{k_2}\to M$ are transverse. (These maps are defined in~\ref{Aoperation}.)

Put $\frak L=V_1\times N_{k_1}\times V_2\times N_{k_2}\times N_{\ov k_2}\times V_3\times N_{k_3}$ and put $\frak R=V_2\times N_{\ov k_2}\times V_3\times N_{k_3} \times N_{k_2}\times V_1\times N_{k_1}.$

Put $P_{k_1, k_2}\subset V_1\times N_{k_1}\times V_2\times N_{k_2}$ to be the manifold $W$ in the definition of $A_{k_1, k_2}(\omega_1, \omega_2).$ Thus 
\begin{equation}
P_{k_1, k_2}=\{(v_1, n_{k_1}, v_2, n_{k_2})| \ov f_{1, k_1}(v_1, n_{k_1})=\ov f_{2, k_2}(v_2, n_{k_2})\}.
\end{equation}

Using Theorem~\ref{transversalitytheorem} and Remark~\ref{pletransversalityconditions} deform $f_3$ so that the adjoint map 
$\ov f_{3, k_3}:V_3\times N_{k_3}\to M$ is transverse to the adjoint maps $\ov f_{2, \ov k_2}:V_2\times N_{\ov k_2}\to M$ and $P_{k_1, k_2}\times N_{\ov k_2}\to M, ((v_1, n_{k_1}, v_2, n_{k_2}), n_{\ov k_2})\to \ov f_{2, k_2}(v_2, n_{\ov k_2}).$

Put $S_{k_1, k_2}=P_{k_1, k_2}\times N_{\ov k_2}\times V_3\times N_{k_3}$ to be the submanifold of $\frak L.$  
Put $L$ to be the oriented $(j_1+j_2+j_3+4\frak n-2\frak m)$-dimensional manifold $W$ constructed according to the definition of $A_{\nu (\Gamma_1)+\ov k_2, k_3}\bigl (A_{k_1, k_2}(\omega_1, \omega_2), \omega_3\bigr).$ Clearly $L\subset S_{k_1, k_2}\subset \frak L$ and 
\begin{equation}\label{LJacobieq}
\begin{split}
L=\bigl \{(v_1, n_{k_1}, v_2, n_{k_2}, n_{\ov k_2}, v_3, n_{k_3})\in \frak L \big | \ov f_{1, k_1}(v_1, n_{k_1})=\ov f_{2, k_2}(v_2, n_{k_2}) \\ \text{ and } \ov f_{2, \ov k_2}(v_2, n_{\ov k_2})=\ov f_{3, k_3}(v_3, n_{k_3})\bigr\}.
\end{split}
\end{equation}
Let $g_L:L\to G_{\frak N, M}^{B_{\nu(\Gamma_1)+\ov k_2, k_3}(B_{k_1, k_2}(\Gamma_1, \Gamma_2), \Gamma_3)} $ be the map constructed according to the definitions of $A_{k_1, k_2}$ and $A_{\nu (\Gamma_1)+\ov k_2, k_3},$
so that 
\begin{equation}[L, g_L]=A_{\nu (\Gamma_1)+\ov k_2, k_3}\bigl (A_{k_1, k_2}(\omega_1, \omega_2), \omega_3\bigr)\in \bor_*\Bigl(G_{\frak N, M}^{B_{\nu(\Gamma_1)+\ov k_2, k_3}(B_{k_1, k_2}(\Gamma_1, \Gamma_2), \Gamma_3)}\Bigr).
\end{equation}

Recall that $\ov f_{3, k_3}$ and $\ov f_{2, \ov k_2}$ are transverse.
Put $P'_{\ov k_2, k_3}\subset V_2\times N_{\ov k_2}\times V_3\times N_{k_3}$ to be the manifold $W$ in the definition of $A_{\ov k_2, k_3} (\omega_2, \omega_3).$ Thus \begin{equation}
P'_{\ov k_2, k_3}=\{(v_2, n_{\ov k_2}, v_3, n_{k_3})\in V_2\times N_{\ov k_2}\times V_3\times N_{k_3}| \ov f_{2, \ov k_2}(v_2, n_{\ov k_2})=\ov f_{3, k_3}(v_3, n_{k_3})\}.
\end{equation}
Put $S'_{\ov k_2, k_3}=P'_{\ov k_2, k_3}\times N_{k_2}\times V_1\times N_{k_1}\subset \mathfrak R.$

Let $P'_{\ov k_2, k_3}\times N_{k_2}\to M, ((v_2, n_{\ov k_2}, v_3, n_{k_3}), n_{k_2})\to \ov f_{2, k_2}(v_2, n_{k_2})$ and $\ov f_{1, k_1}:V_1\times N_{k_1}\to M$ be the natural maps we use to define $A_{k_2, k_1}(A_{\ov k_2, k_3}(\omega_2, \omega_3), \omega_1).$ The product of these maps gives the map $S'_{\ov k_2, k_3}\to M\times M.$

Put $R\subset S'_{\ov k_2, k_3}\subset \frak R$ to be the preimage of $\Delta$ under this map. 
Clearly
\begin{equation}\label{RJacobieq}
\begin{split}
R=\bigl \{  (v_2, n_{\ov k_2}, v_3, n_{k_3}, n_{k_2}, v_1, n_{k_1})\in \frak R
 \big| \ov f_{2, \ov k_2}(v_2, n_{\ov k_2})=\ov f_{3, k_3}(v_3, n_{k_3})\\ \text{ and } \ov f_{2, k_2}(v_2, n_{k_2})=\ov f_{1, k_1}(v_1, n_{k_1})\bigr \}.
\end{split}
\end{equation}

Put $\Theta: \frak L\to \frak R$ to be the diffeomorphism permuting the factors defined by 
$\Theta(v_1, n_{k_1}, n_2, n_{k_2},n_{\ov k_2}, v_3, n_{k_3})=(v_2, n_{\ov k_2}, v_3, n_{k_3}, n_{k_2}, v_1,  n_{k_1}),$ for all $v_1\in V_1, v_2\in V_2, v_3\in V_3, n_{k_1}\in N_{k_1}, n_{k_2}\in N_{k_2}, n_{\ov k_2}\in N_{\ov k_2}, n_{k_3}\in N_{k_3}.$ 

Comparing equations~\eqref{LJacobieq} and~\eqref{RJacobieq} that define $L$ and $R,$ we see that $R$ is a submanifold and $\Theta|_L:L\to R$ is a diffeomorphism. Proposition~\ref{automatictransversality} says that the maps $P'_{\ov k_2, k_3}\times N_{k_2}\to M, ((v_2, n_{\ov k_2}, v_3, n_{k_3}), n_{k_2})\to \ov f_{2, k_2}(v_2, n_{k_2})$ and $\ov f_{1, k_1}:V_1\times N_{k_1}\to M$ are transverse. Hence $R$ is in fact the manifold $W$ constructed according to the definition of  $A_{k_2, k_1}(A_{\ov k_2, k_3} (\omega_2, \omega_3), \omega_1\bigr).$ 

We orient it according to the definition of $A_{k_2, k_1}$ and $A_{\ov k_2, k_3}$ and we put 
$g_R:R\to G_{\frak N, M}^{B_{k_2, k_1}(B_{\ov k_2, k_3}(\Gamma_2, \Gamma_3), \Gamma_1)} $ to be the corresponding nice map so that 
\begin{equation}[R, g_R]=A_{k_2, k_1}\bigl (A_{\ov k_2, k_3}(\omega_2, \omega_3), \omega_1\bigr) \in \bor_*\Bigl(G_{\frak N, M}^{B_{k_2, k_1}(B_{\ov k_2, k_3}(\Gamma_2, \Gamma_3), \Gamma_1)} \Bigr).
\end{equation}

Take 
\begin{equation} 
\begin{split}
l=(v_1, n_{k_1}, n_2, n_{k_2},n_{\ov k_2}, v_3, n_{k_3})\in L  \text{ and }  \\ r=\Theta(l)=(v_2, n_{\ov k_2}, v_3, n_{k_3}, n_{k_2}, v_1,  n_{k_1})\in R.
\end{split}
\end{equation}
Applying twice the definition of the $A$ operation, we see that 
the commutative diagram $g_L(l)$
is obtained by the following sequence of steps:
\begin{description}
\item[1]  take the disjoint union $\bigl(f_1(v_1)\sqcup f_2(v_2)\bigr)\sqcup f_3(v_3)$ of commutative diagrams;
\item[2] add the extra $\pt$-space that is mapped to $n_{k_1}\in N_{k_1}$ and to $n_{k_2}\in N_{k_2}$; 
\item[3] add the extra $\pt$-space that is mapped to $n_{\ov k_2}\in N_{\ov k_2}$ and to $n_{k_3}\in N_{k_3};$ 
\item[4] identify the copies of $M$ in the commutative diagrams $f_1(v_1), f_2(v_2), f_3(v_3);$ and
\item[5] change the indices of the $N$-manifolds and the labels of the $\pt$-spaces. 
\end{description}
Similarly, the commutative diagram $g_R(r)=g_R(\Theta(l))$
is obtained by the following sequence of steps
\begin{description}
\item[1]  take the disjoint union $\bigl(f_2(v_2)\sqcup f_3(v_3)\bigr)\sqcup f_1(v_1)$ of commutative diagrams;
\item[2] add the extra $\pt$-space that is mapped to $n_{\ov k_2}\in N_{\ov k_2}$ and to $n_{k_3}\in N_{k_3};$ 
\item[3] add the extra $\pt$-space that is mapped to $n_{k_1}\in N_{k_1}$ and to $n_{k_2}\in N_{k_2}$; 
\item[4] identify the copies of $M$ in the commutative diagrams $f_1(v_1), f_2(v_2), f_3(v_3);$ and
\item[5] change the indices of the $N$-manifolds and the labels of the $\pt$-spaces.
\end{description}

So, up the change of the indices of $N$-manifolds and of the labels of the $\pt$-spaces, the commutative diagrams $g_L(l)$ and $g_R(\Theta(l))$ are the same.
Using Statement $1$ of Proposition~\ref{propgraph} and Definition~\ref{groupaction}, we see that $\wt \alpha(g_R(\Theta(l)))=g_L(l).$
Thus 
%\begin{equation}
$[L, g_L]=[L, \wt \alpha\circ g_R\circ \Theta]=\alpha_*([L, g_R\circ \Theta])=\pm 1\alpha_*([R, g_R])\in \bor_*\bigl(G_{\frak N, M}^{B_{\nu(\Gamma_1)+\ov k_2, k_3}(B_{k_1, k_2}(\Gamma_1, \Gamma_2), \Gamma_3)}\bigr).$
%\end{equation}

The manifolds $L$ and $R$ are oriented as it is described in the definition of $A_{j_1, j_2}^{k_1, k_2}$ and the statement of the Theorem would follow immediately if we show that 
$$(\Theta|_{L})_*:TL\to (-1)^{\sigma} TR
$$ 
is orientation preserving. 
(Recall that a product of a sign and an oriented bundle denotes the initial oriented bundle if the sign is $+1,$ and it denotes the initial bundle with the changed orientation if the sign is $-1.$)

{\it Below we provide a computation showing that this map of bundles is indeed orientation preserving. For orientations and frames we use the notation introduced in the proof of Theorem~\ref{skewsymmetryingredient}.\/}
Take $l=(v_1, n_{k_1}, v_2, n_{k_2}, n_{\ov k_2}, v_3, n_{k_3})\in L.$ Put $m=\ov f_{1, k_1}(v_1, n_{k_1})=\ov f_{2, k_2}(v_2, n_{k_2})\in M$ and $m'=\ov f_{2, k_2}(v_2, n_{\ov k_2})=\ov f_{3, k_3}(v_3, n_{k_3})\in M.$ Put $\{ \vec {m_1}, \cdots, \vec {m_{\frak m}}\}$ to be an orthogonal positive $\frak m$-frame in $T_mM$ and $\{ \vec {m'_1}, \cdots, \vec {m'_{\frak m}}\}$ to be an orthogonal positive $\frak m$-frame in $T_{m'}M.$ Then $\{ \vec {m_1}\oplus \vec {m_1}, \cdots, \vec {m_{\frak m}}\oplus \vec {m_{\frak m}}\}$ and $\{ \vec {m'_1}\oplus \vec {m'_1}, \cdots, \vec {m'_{\frak m}}\oplus \vec {m'_{\frak m}}\}$ are positive $\frak m$-frames in $T_{m\times m}\Delta$ and in $T_{m'\times m'}\Delta,$ respectively. 
Hence $\Fr_{m\times m}^{\perp}=\{ - \vec {m_1}\oplus \vec {m_1}, \cdots, - \vec {m_{\frak m}}\oplus \vec {m_{\frak m}}\}$ and $\Fr_{m'\times m'}^{\perp}=\{ -\vec {m'_1}\oplus \vec {m'_1}, \cdots, - \vec {m'_{\frak m}}\oplus \vec {m'_{\frak m}}\}$ are  positive orientation frames in the fibers of $(T\Delta)^{\perp}$ over $m\times m$ and  over $m'\times m',$ respectively.

Put $\Fr_{k_1, k_2}$ to be the unique $\frak m$-frame in  $T_{(v_1, n_{k_1}, v_2, n_{k_2})} (V_1\times N_{k_1}\times V_2\times N_{k_2})$ orthogonal to $T_{(v_1, n_{k_1}, v_2, n_{k_2})}P_{k_1, k_2}$ and satisfying $\pr \bigl((\ov f_{1, k_1}\times \ov f_{2, k_2})_*(\Fr_{k_1, k_2})\bigr)=\Fr^{\perp}_{m\times m}.$ By our orientation convention we have
\begin{equation}\label{jaceq1}
\Or\bigl(\Fr^+\bigl(T_{(v_1, n_{k_1}, v_2, n_{k_2})}
P_{k_1, k_2}\bigr), \Fr_{k_1, k_2}\bigr)=\Or\bigl(T_{(v_1, n_{k_1}, v_2, n_{k_2})}(V_1\times N_{k_1}\times V_2\times N_{k_2})\bigr). 
\end{equation}
Put 
$\wt \Fr_{k_1, k_2}$ to be the $\frak m$-frame in $T_l\frak L$ that is the image of $\Fr_{k_1, k_2}$ under the differential of the inclusion 
$V_1\times N_{k_1}\times V_2\times N_{k_2}\to V_1\times N_{k_1}\times V_2\times N_{k_2}\times n_{\ov k_2}\times v_3\times n_{k_3}\subset \mathfrak L.$ From~\eqref{jaceq1} we get 
\begin{equation}\label{jaceq2}
\Or\bigl(\Fr^+\bigl(T_l(P_{k_1, k_2}\times N_{\ov k_2}\times V_3\times N_{k_3})\bigr), \wt \Fr_{k_1,k_2}\bigr)=(-1)^{\sigma_1}\Or(T_l\frak L),
\end{equation}
where $\sigma_1=\frak mj_3.$

Put $\Fr_{L}$ to be a positive orientation $(j_1+j_2+j_3+4\frak n-2\frak m)$-frame of $T_lL.$ We will often 
regard it as a frame in $T_l\frak L$ obtained by the inclusion $T_lL\to T_l\frak L.$
Put $\Fr_{L}^{\perp}\subset T_l(P_{k_1, k_2}\times N_{\ov k_2}\times V_3\times N_{k_3})$ to be the unique $\frak m$-frame orthogonal to $T_lL$ such that it is mapped to $\Fr_{m'\times m'}^{\perp}$ under the composition of $\pr$ and of the differential of the map $P_{k_1, k_2}\times N_{\ov k_2}\times V_3\times N_{k_3}\to M\times M$ used to define $L.$ Then 
\begin{equation}\label{jaceq3}
\Or(\{\Fr_{L}, \Fr_{L}^{\perp}\})=\Or\bigl(T_l(P_{k_1, k_2}\times N_{\ov k_2}\times V_2\times N_{k_3})\bigr). 
\end{equation}
From~\eqref{jaceq2} and~\eqref{jaceq3} we get
\begin{equation}\label{jaceq4}
\Or\bigl(\{\Fr_{L}, \Fr_{L}^{\perp}, \wt \Fr_{k_1, k_2}\}\bigr)=(-1)^{\sigma_1}\Or(T_l \frak L).
\end{equation}

Put $r=\Theta(l)=(v_2, n_{\ov k_2}, v_3, n_{k_3}, n_{k_2}, v_1, n_{k_1})\in R.$
We define $\Fr'_{\ov k_2, k_3}$ to be the unique $\frak m$-frame in 
$T_{(v_2, n_{\ov k_2}, v_3, n_{k_3})} (V_2\times N_{\ov k_2}\times V_3\times N_{k_3})$ that is orthogonal to $T_{(v_2, n_{\ov k_2}, v_3, n_{k_3})}P'_{\ov k_2, k_3}$ and
that satisfies $\pr \bigl((\ov f_{2, \ov k_2}\times \ov f_{3, k_3})_*(\Fr^{\perp}_{\ov k_2, k_3})\bigr)=\Fr^{\perp}_{m'\times m'}.$ 
Put $\wt \Fr'_{\ov k_2,k_3}$ to be the $\frak m$-frame in $T_r\mathfrak R$ that is the image of $\Fr'_{\ov k_2, k_3}$ under the differential of the inclusion $V_2\times N_{\ov k_2}\times V_3\times N_{k_3}\to V_2\times N_{\ov k_2}\times V_3\times N_{k_3}\times n_{k_2}\times v_1\times n_{k_1}\subset \frak R.$

Put $\Fr_{R}$ to be a positive orientation $(j_1+j_2+j_3+4\frak n-2\frak m)$-frame of $T_rR$ that we will often identify with the frame in $T_r\frak R$ obtained by the inclusion $T_rR\to T_r\frak R.$
Put $\Fr_{R}^{\perp}\subset T_r(P'_{\ov k_2, k_3}\times N_{k_2}\times V_1\times N_{k_1})$ to be the unique $\frak m$-frame orthogonal to $T_rR$ such that it is mapped to $\Fr_{m\times m}^{\perp}$ under the composition of $\pr$ and of the differential of the map $P'_{\ov k_2, k_3}\times N_{k_2}\times V_1\times N_{k_1}\to M\times M$ that was used to define $R.$ Similarly to the computations above we get that 
\begin{equation}\label{jaceq5}
\Or\bigl(\{\Fr_{R}, \Fr_{R}^{\perp}, \wt \Fr'_{\ov k_2, k_3}\}\bigr)=(-1)^{\sigma_2}\Or(T_r\frak R),
\end{equation}
where $\sigma_2=\frak m j_1.$

Since $\Theta$ is defined as the permutation of the factors, it is easy to see that 
\begin{equation}\label{jaceq6}
\Or\bigl(\Theta_*(\Fr^+(T_l\frak L))\bigr)=(-1)^{\sigma_3}\Or(T_r\frak R),
\end{equation}
where $\sigma_3=(j_1+\frak n)(j_2+j_3+\frak n)+\frak n j_3.$

Using~\eqref{jaceq4} and~\eqref{jaceq6} we get 
\begin{equation}\label{jaceq7}
\Or \{\Theta_*(\Fr_l), \Theta_*(\Fr_{L}^{\perp}), \Theta_*(\wt \Fr_{k_1, k_2})\}=(-1)^{\sigma_1+\sigma_3}\Or(T_r\frak R).
\end{equation}
Since the frames $\wt \Fr_{k_1, k_2}$ and $\Fr_{L}^{\perp}$ consist of $\frak m$ vectors each, equation~\eqref{jaceq7} implies
\begin{equation}\label{jaceq8} \Or\bigl(\{\Theta_*(\Fr_L),  \Theta_*(\wt \Fr_{k_1, k_2}), \Theta_*(\Fr_{L}^{\perp}),\}\bigr)=(-1)^{\sigma_1+\sigma_3+\frak m^2}\Or(T_r\frak R).
\end{equation}

Put $-Fr_{R}^{\perp}$ to be the frame obtained by plying each vector in $Fr_{R}^{\perp}$ by $-1$ and keeping the order of the vectors intact. Propositions~\ref{framehomotopyone} and~\ref{framehomotopytwo} imply that the 
frames $\{\Theta_*(\Fr_L),  \Theta_*(\wt \Fr_{k_1, k_2}), \Theta_*(\Fr_{L}^{\perp}),\}$
and $\{\Theta_*(\Fr_L), -\Fr_{R}^{\perp}, \wt \Fr'_{\ov k_2, k_3}\}$ 
give equal orientations of $T_r\frak R.$ 
Thus 
\begin{equation}\label{jaceq9}
\Or\bigl(\{\Theta_*(\Fr_L), -\Fr_{R}^{\perp}, \wt \Fr'_{\ov k_2, k_3}\}\bigr)=
(-1)^{\sigma_1+\sigma_3+\frak m^2}
\Or(T_r\frak R).
\end{equation}
Hence 
\begin{equation}\label{jaceq10}
\Or\bigl(\{\Theta_*(\Fr_L), \Fr_{R}^{\perp}, \wt \Fr'_{\ov k_2, k_3}\}\bigr)=
(-1)^{\sigma_1+\sigma_3}\Or(T_r\frak R).
\end{equation}

Use equations~\eqref{jaceq5},~\eqref{jaceq10}, identities $\Or(\Fr_R)=\Or(T_rR), \Or(\Fr_L)=\Or(T_lL)$ and the fact that $\Theta_*(\Fr_L)$ gives an orientation of $T_rR,$ to get
\begin{equation}
\Or(\Theta_*(\Fr_L))=(-1)^{\sigma_1+\sigma_2+\sigma_3}\Or(T_rR)=
(-1)^{\sigma}\Or(T_rR).
\end{equation}
This  finishes the proof of Theorem~\ref{Jacobiingredient} modulo the proofs of Propositions~\ref{automatictransversality},~\ref{framehomotopyone} and~\ref{framehomotopytwo}.\qed

\begin{prop}\label{automatictransversality} The maps $\ov f_{1, k_1}:V_1\times N_{k_1}\to M$ and 
$P'_{\ov k_2, k_3}\times N_{k_2}\to M,$ $((v_2, n_{\ov k_2}, v_3, n_{k_3}), n_{k_2})\to \ov f_{2, k_2}(v_2, n_{k_2})$ 
from the proof of Theorem~\ref{Jacobiingredient} are transverse.
\end{prop}

\pp Put $S_{\ov k_2, k_3}=\Theta^{-1}(S'_{\ov k_2, k_3})\subset \mathfrak L$ and $S'_{k_1, k_2}=\Theta(S_{k_1, k_2})\subset \mathfrak R$ to be the submanifolds.
Clearly $L=S_{k_1, k_2}\cap S_{\ov k_2, k_3}\subset \mathfrak L.$ Let us show that $S_{k_1, k_2}$ and $S_{\ov k_2, k_3}$ are transverse.

Choose $l=(v_1, n_{k_1}, v_2, n_{k_2}, n_{\ov k_2}, v_3, n_{k_3})\in L=S_{k_1, k_2}\cap S_{\ov k_2, k_3}.$
The vector spaces $T_l\frak L$ and $T_{\Theta(l)}\frak R$ are naturally decomposed as the direct sums of the linear subspaces $T_{v_1}V_1,$ $T_{v_2}V_2,$ $T_{v_2}V_3,$ $T_{n_{k_1}}N_{k_1},$ $T_{n_{k_2}}N_{k_2},$ $T_{n_{\ov k_2}}N_{\ov k_2},$ $T_{n_{k_3}}N_{k_3}.$ 
For a vector $\vec w$ in $T_l\frak L$ or in $T_{\Theta(r)}\frak R$ we will denote by $\vec w_{V_1},$ $\vec w_{V_2},$ $\vec w_{V_3},$ $\vec w_{N_{k_1}},$ $\vec w_{N_{k_2}},$ $\vec w_{N_{\ov k_2}},$ $\vec w_{N_{k_3}}$
the components of $\vec w$ in the corresponding linear subspaces. {\it Note that since $\Theta$ is defined by permuting the factors, the components of $\vec w\in T_l\frak L$ and of $\Theta_*(\vec w)\in T_{\Theta(l)}\frak R$ in the corresponding subspaces are equal. \/}

By the definition of $S_{\ov k_2, k_3},$  $\pr \bigl((\ov f_{2, \ov k_2}\times \ov f_{3, k_3})_*(\vec w_{V_2}, \vec w_{N_{\ov k_2}}, \vec w_{V_3}, \vec w_{N_{k_3}})\bigr)=\vec 0\in (T\Delta)^{\perp}\subset TM\times TM$ for every $\vec w\in T_lS_{\ov k_2, k_3}.$ On the other hand by definition of $L,$ there are $\frak m$ linearly independent vectors $\vec u^i, i=1, \cdots, \frak m$ in 
$T_{l}S_{k_1, k_2}$ such that $\pr \bigl((\ov f_{2, \ov k_2}\times \ov f_{3, k_3})_*(\vec u^i_{V_2}, \vec u^i_{N_{\ov k_2}}, \vec u^i_{V_3}, \vec u^i_{N_{k_3}})\bigr)\neq \vec 0, i=1, \cdots, \frak m.$ Since $S_{\ov k_2, k_3}$ is a codimension $\frak m$ submanifold of $\frak L,$ the linear span of $T_lS_{k_1, k_2}$ and $T_lS_{\ov k_2, k_3}$ is $T_l\frak L.$  Hence the intersection $S_{k_1, k_2}\cap S_{\ov k_2, k_3}$ is indeed transverse.

Since $\Theta$ is a diffeomorphism, we get that 
$S'_{k_1, k_2}=\Theta(S_{k_1, k_2})$ and $S'_{\ov k_2, k_3}=\Theta(S_{\ov k_2, k_3})$ intersect transversally in $\mathfrak R.$
Moreover $R$ is the transverse intersection $S'_{\ov k_2, k_3}\cap S'_{k_1, k_2}.$

To prove the Proposition it suffices to show that the natural product map $S'_{\ov k_2, k_3}=(P'_{\ov k_2, k_3}\times N_{k_2})\times (V_1\times N_{k_1})\to M\times M$ is transverse to the diagonal $\Delta.$
$R$ is a codimension $\frak m$ submanifold of 
$S'_{\ov k_2, k_3}.$ Thus to show that $S'_{\ov k_2, k_3}\to M\times M$ is transverse to $\Delta,$ it suffices to prove that if $r\in R$ and $\vec w'\in T_rS'_{\ov k_2, k_3}\setminus T_rR,$ then $\pr (\ov f_{2, k_2}\times \ov f_{1, k_1})_*(\vec w'_{V_2}, \vec w'_{N_{k_2}}, \vec w'_{V_1}, \vec w'_{N_{k_1}})\neq \vec 0\in (T\Delta)^{\perp}.$ Since $R$ is the transverse intersection $S'_{\ov k_2, k_3}\cap S'_{k_1, k_2},$ we have $\vec w'\not \in T_rS'_{k_1, k_2}.$
Put $\vec w=(\Theta^{-1})_*(\vec w')\in T_{\Theta^{-1}(r)}S_{\ov k_2, k_3}.$ Clearly $\vec w\not \in T_{\Theta^{-1}(r)}S_{k_1, k_2}.$

Since $\Theta_*$ preserves the vector components, it suffices to show that 
$\pr \bigl((\ov f_{2, k_2}\times \ov f_{1, k_1})_*(\vec w_{V_2}, \vec w_{N_{k_2}}, \vec w_{V_1}, \vec w_{N_{k_1}})\bigr)\neq \vec 0\in (T\Delta)^{\perp},$
or which is the same that $\pr \bigl((\ov f_{1, k_1}\times \ov f_{2, k_2})_*(\vec w_{V_1}, \vec w_{N_{k_1}}, \vec w_{V_2}, \vec w_{N_{k_2}})\bigr)\neq \vec 0\in (T\Delta)^{\perp}.$
Recall that $S_{k_1, k_2}=P_{k_1, k_2}\times N_{\ov k_2}\times V_3\times N_{k_3},$ where $P_{k_1, k_2}$ is the transverse preimage of $\Delta$ under $\ov f_{1, k_1}\times \ov f_{2, k_2}:V_1\times N_{k_1}\times V_2\times N_{k_2}\to M\times M.$ So for every $\vec u\in T_{\Theta^{-1}(r)}\frak L\setminus T_{\Theta^{-1}(r)}S_{k_1, k_2}$ we have 
$\pr \bigl((\ov f_{1, k_1}\times \ov f_{2, k_2})_*(\vec u_{V_1}, \vec u_{N_{k_1}}, \vec u_{V_2}, \vec u_{N_{k_2}})\bigr)\neq \vec 0\in (T\Delta)^{\perp}.$ Put $\vec u=\vec w$ to get the desired statement. So $S'_{\ov k_2, k_3}\to M\times M$ is indeed transverse to the diagonal $\Delta.$ 

\qed

\begin{prop}\label{framehomotopyone} We have the equality of orientations of $T_r\frak R$ 
$$\Or(\{\Theta_*(\Fr_L),  \Theta_*(\wt \Fr_{k_1, k_2}),\Theta_*(\Fr_{L}^{\perp})\})=
\Or(\{\Theta_*(\Fr_{L}),-\Fr_{R}^{\perp},
\Theta_*(\Fr_{L}^{\perp})\}).$$
\end{prop}

\pp

Let $\vec u^i$ be the $i$-th vector in $\Fr_{R}^{\perp}.$ By definition of $\Fr_{R}^{\perp}$ we have $\pr\bigl((\ov f_{2, k_2}\times \ov f_{1, k_1})_*(\vec u^i_{V_2}, \vec u^i_{N_{k_2}}, \vec u^i_{V_1}, \vec u^i_{N_{k_1}})\bigr)=-\vec {m_i}\oplus \vec {m_i}.$ Thus 
\begin{equation}\label{lineq-1lem1}
\pr\bigl((\ov f_{1, k_1}\times \ov f_{2, k_2})_*(-\vec u^i_{V_1}, -\vec u^i_{N_{k_1}}, -\vec u^i_{V_2}, -\vec u^i_{N_{k_2}})\bigr)=-\vec {m_i}\oplus \vec {m_i}.
\end{equation}
Let $\vec w^i$  be the $i$-th vector in $\wt \Fr_{k_1, k_2}.$ By definition of $\wt \Fr_{k_1, k_2}$ we get $\pr\bigl((\ov f_{1, k_1}\times \ov f_{2, k_2})_*(\vec w^i_{V_1}, \vec w^i_{N_{k_1}}, \vec w^i_{V_2}, \vec w^i_{N_{k_2}})\bigr)= -\vec {m_i}\oplus \vec {m_i}.$ Since $\Theta$ is defined by permuting the factors, $\Theta_*$ preserves the components of the vectors in the subspaces $T_{v_1}V_1,$ $T_{v_2}V_2,$ $T_{v_2}V_3,$ $T_{n_{k_1}}N_{k_1},$ $T_{n_{k_2}}N_{k_2},$ $T_{n_{\ov k_2}}N_{\ov k_2},$ $T_{n_{k_3}}N_{k_3}.$  Hence 
\begin{equation}\label{lineq0lem1}
\pr\Bigl((\ov f_{1, k_1}\times \ov f_{2, k_2})_*(\Theta_*(\vec w^i)_{V_1}, \Theta_*(\vec w^i)_{N_{k_1}}, \Theta_*(\vec w^i)_{V_2}, \Theta_*(\vec w^i)_{N_{k_2}})\Bigr)= -\vec {m_i}\oplus \vec {m_i}.
\end{equation}

Let $\Fr'_t, t\in [0,1],$ be the time-dependent $\frak m$-frame whose $i$-th vector at time $t$ is $\vec z^{i,t}=t\Theta_*(\vec w_i)+(1-t)(-\vec u^i).$ Since $\pr$ and  $(\ov f_{1, k_1}\times \ov f_{2, k_2})_*$ are linear, equations~\eqref{lineq-1lem1},~\eqref{lineq0lem1} imply that 
\begin{equation}\label{lineq1lem1}
\pr\bigl((\ov f_{1, k_1}\times \ov f_{2, k_2})_*(\vec z^{i,t}_{V_1}, \vec z^{i,t}_{N_{k_1}}, \vec z^{i,t}_{V_2}, \vec z^{i,t}_{N_{k_2}})\bigr)= -\vec {m_i}\oplus \vec {m_i},
\end{equation}
for all $i\in \{1, \cdots, \frak m\}$ and $t\in [0,1].$

Since $L\subset P_{k_1, k_2}\times N_{\ov k_2}\times V_3\times N_{k_3}$ and $(\ov f_{1, k_1}\times \ov f_{2, k_2})(v_1, n_{k_1}, v_2, n_{k_2})\in \Delta,$ for every  $(v_1, n_{k_1}, v_2, n_{k_2})\in P_{k_1, k_2}\subset (V_1\times N_{k_1}\times V_2\times N_{k_2}),$ we have that
\begin{equation}
\begin{split}
(\ov f_{1, k_1}\times \ov f_{2, k_2})_*(\vec x_{V_1}, \vec x_{N_{k_1}}, \vec x_{V_2}, \vec x_{N_{k_2}})\in T_{m\times m}\Delta \text{ and } \\
\pr\bigl((\ov f_{1, k_1}\times \ov f_{2, k_2})_*(\vec x_{V_1}, \vec x_{N_{k_1}}, \vec x_{V_2}, \vec x_{N_{k_2}})\bigr)=\vec 0,
\end{split}
\end{equation}
for any vector $\vec x$ in the frames $\Fr_L$ and $\Fr_L^{\perp}.$ Since $\Theta_*$ preserves the components of the vectors, 
we have 
\begin{equation}\label{lineq2lem1}
\pr \bigl( (\ov f_{1, k_1}\times \ov f_{2, k_2})_*(\vec x_{V_1}, \vec x_{N_{k_1}}, \vec x_{V_2}, \vec x_{N_{k_2}})\bigr)=\vec 0,
\end{equation} for any vector $\vec x$ in the frames $\Theta_*(\Fr_L)$ and $\Theta_*(\Fr_L^{\perp}).$

Equations~\eqref{lineq1lem1},~\eqref{lineq2lem1} and the fact that the vectors in $\{\Theta_*(\Fr_L), \Theta_*(\Fr_{L}^{\perp})\}$ are linearly independent imply that
the frame  $\{\Theta_*(\Fr_L),  \Fr'_t, \Theta_*(\Fr_{L}^{\perp})\}$ is non-degenerate, for all $t\in [0,1].$ Thus 
$$\{\Theta_*(\Fr_L),  \Fr'_0, \Theta_*(\Fr_{L}^{\perp})\}=\{\Theta_*(\Fr_L),  -\Fr_R^{\perp}, \Theta_*(\Fr_{L}^{\perp})\} \text{ and }$$
$$\{\Theta_*(\Fr_L),  \Fr'_1, \Theta_*(\Fr_{L}^{\perp})\}=\{\Theta_*(\Fr_L),  \Theta_*(\wt \Fr_{k_1, k_2}), \Theta_*(\Fr_{L}^{\perp})\}$$
give equal orientations of $T_r\frak R.$ \qed

\begin{prop}\label{framehomotopytwo}
We have the equality of orientations of $T_r\frak R$
$$\Or(\{\Theta_*(\Fr_{L}), -\Fr_{R}^{\perp},\Theta_*(\Fr_{L}^{\perp})\})=\Or(\{\Theta_*(\Fr_{L}), -\Fr_{R}^{\perp}, \wt \Fr'_{\ov k_2, k_3}\}).$$
\end{prop}

\pp We use notation conventions of the proof of Proposition~\ref{framehomotopyone}.
Let $\vec w^i$  be the $i$-th vector in the frame $\wt \Fr'_{\ov k_2, k_3}.$ By definition of $\wt \Fr'_{\ov k_2, k_3}$ we have that 
\begin{equation}\label{lineq-1lem2}
\pr\bigl((\ov f_{2, \ov k_2}\times \ov f_{3, k_3})_*(\vec w^i_{V_2}, \vec w^i_{N_{\ov k_2}}, \vec w^i_{V_3}, \vec w^i_{N_{k_3}})\bigr)= -\vec {m'_i}\oplus \vec {m'_i}.
\end{equation}
Let $\vec u^i$ be the $i$-th vector in $\Fr_{L}^{\perp}.$ By definition of $\Fr_{L}^{\perp}$ we have that $\pr\Bigl((\ov f_{2, \ov k_2}\times \ov f_{3, k_3})_*(\vec u^i_{V_2}, \vec u^i_{N_{\ov k_2}}, \vec u^i_{V_3}, \vec u^i_{N_{k_3}})\Bigr)=-\vec {m'_i}\oplus \vec {m'_i}.$ Since $\Theta_*$ preserves the components of the vectors in the subspaces, we have
\begin{equation}\label{lineq0lem2}
\pr\bigl((\ov f_{2, \ov k_2}\times \ov f_{3, k_3})_*(\Theta_*(\vec u^i)_{V_2}, \Theta_*(\vec u^i)_{N_{\ov k_2}}, \Theta_*(\vec u^i)_{V_3}, \Theta_*(\vec u^i)_{N_{k_3}})\bigr)=-\vec {m'_i}\oplus \vec {m'_i}.
\end{equation}

Let $\Fr''_t, t\in [0,1],$ be the 
time-dependent $\frak m$-frame whose $i$-th vector at time $t$ is $\vec z^{i,t}=t\vec w_i+(1-t)\Theta_*(\vec u^i).$ Since $\pr$ and  $(\ov f_{2, \ov k_2}\times \ov f_{3, k_3})_*$ are linear, equations~\eqref{lineq-1lem2},~\eqref{lineq0lem2} imply that 
\begin{equation}\label{lineq1lem2}
\pr\bigl((\ov f_{2, \ov k_2}\times \ov f_{3, k_3})_*(\vec z^{i,t}_{V_2}, \vec z^{i,t}_{N_{\ov k_2}}, \vec z^{i,t}_{V_3}, \vec z^{i,t}_{N_{k_3}})\bigr)= -\vec {m'_i}\oplus \vec {m'_i},
\end{equation}
for all $i\in \{1, \cdots, \frak m\}$ and $t\in [0,1].$

By definition of $L$ we have that $(\ov f_{2, \ov k_2}\times \ov f_{3, k_3})(v_2, n_{\ov k_2}, v_3, n_{k_3})\in \Delta,$ for every $(v_1, n_{k_1}, v_2, n_{k_2}, n_{\ov k_2}, v_3, n_{k_3})\in L\subset \frak L.$ Thus $(\ov f_{2, \ov k_2}\times \ov f_{3, k_3})_*(\vec x_{V_2}, \vec x_{N_{\ov k_2}}, \vec x_{V_3}, \vec x_{N_{k_3}})\in T_{m\times m}\Delta$ 
and $\pr\bigl((\ov f_{2, \ov k_2}\times \ov f_{3, k_3})_*(\vec x_{V_2}, \vec x_{N_{\ov k_2}}, \vec x_{V_3}, \vec x_{N_{k_3}})\bigr)=\vec 0,$
for every $\vec x\in \Fr_L.$ Since $\Theta_*$ preserves the components of the vectors, we have that
\begin{equation}\label{lineq2lem2}
\pr\bigl((\ov f_{2, \ov k_2}\times \ov f_{3, k_3})_*(\vec x_{V_2}, \vec x_{N_{\ov k_2}}, \vec x_{V_3}, \vec x_{N_{k_3}})\bigr)=\vec 0,
\end{equation} 
for every $\vec x\in \Theta_*(\Fr_L).$

Since $(\ov f_{2, \ov k_2}\times \ov f_{3, k_3})(v_2, n_{\ov k_2}, v_3, n_{k_3})\in \Delta,$ for every $(v_2, n_{\ov k_2}, v_3, n_{k_3})\in P'_{\ov k_2, k_3}$ and since $R\subset P'_{\ov k_2, k_3}\times N_{k_2}\times V_1\times N_{k_1},$ we get that $(\ov f_{2, \ov k_2}\times \ov f_{3, k_3})_*(\vec x_{V_2}, \vec x_{N_{\ov k_2}}, \vec x_{V_3}, \vec x_{N_{k_3}})\in T_{m\times m}\Delta,$ for every $\vec x\in \Fr_R.$
Thus
\begin{equation}\label{lineq3lem2}
\pr\bigl((\ov f_{2, \ov k_2}\times \ov f_{3, k_3})_*(\vec x_{V_2}, \vec x_{N_{\ov k_2}}, \vec x_{V_3}, \vec x_{N_{k_3}})\bigr)=\vec 0,
\end{equation} 
for every $\vec x\in -\Fr_R.$

Since the vectors in $\{\Theta_*(\Fr_L),  -\Fr_{R}^{\perp}\}$ are linearly independent and 
equations \eqref{lineq1lem2}, \eqref{lineq2lem2}, \eqref{lineq3lem2} hold, we get that
the frame  $\{\Theta_*(\Fr_L),  -\Fr_{R}^{\perp}, \Fr''_t\}$ is non-degenerate, for all $t\in [0,1].$ Hence 
$$\{\Theta_*(\Fr_L),  -\Fr_{R}^{\perp}, \Fr''_0\}=\{\Theta_*(\Fr_{L}), -\Fr_{R}^{\perp},\Theta_*(\Fr_{L}^{\perp})\} \text{ and }$$ 
$$\{\Theta_*(\Fr_L),  -\Fr_{R}^{\perp}, \Fr''_1\}
=\{\Theta_*(\Fr_{L}), -\Fr_{R}^{\perp}, \wt \Fr'_{\ov k_2, k_3}\}$$
give equal orientations of $T_r\frak R.$ \qed

This finishes the proof of Theorem~\ref{Jacobiingredient}.

\section{One more bracket that can be constructed when $\frak N$ is a one-element set}\label{coincedenceLiebracket}
If $\frak N=\{N\}$ is a one-element set, then one can construct one more important operation $\lceil \cdot, \cdot \rceil$ on $\bor_*(G_{\frak N, M})\otimes \Q$. The operation is of degree $\frak n-\frak m$. It is constructed similar to $[\cdot, \cdot]$ but the operation $A$ is modified so that instead of the pullback submanifold $W\subset V_1\times N\times V_2\times N$ we take the intersection of $W$ with the diagonal 
subset $\{(v_1, n, v_2, n), n\in N, v_1\in V_1, v_2\in V_2\}$. To make the identities look better we make the following substitution: for odd $\frak n$ we put $\lfloor \omega_1, \omega_2 \rfloor=(-1)^{\frak m(|\omega_2|+1)}\lceil \omega_1, \omega_2 \rceil;$ for even $\frak n$ we put $\lfloor \omega_1, \omega_2 \rfloor=(-1)^{\frak m |\omega_2|}\lceil \omega_1, \omega_2 \rceil.$

At the moment we do not know how to generalize the transversality result of Appendix~\ref{transversality} so that it takes care of the construction of the operation $\lceil \cdot, \cdot \rceil$. Modulo this one can get the following result.

\begin{thm}\label{main3}
For odd $\frak n$ we have the following identities 
\begin{enumerate}
\item graded skew-symmetry $\lfloor \omega_1, \omega_2 \rfloor=-(-1)^{(|\omega_1|+\frak m+1)(|\omega_2|+\frak m+1)}\lfloor \omega_2, \omega_1\rfloor$
\item graded Jacobi identity $$\lfloor \omega_1, \lfloor \omega_2, \omega_3\rfloor\rfloor=\lfloor \lfloor \omega_1, \omega_2\rfloor, \omega_3\rfloor+(-1)^{(|\omega_1|+\frak m+1)(|\omega_2|+\frak m+1)}\lfloor \omega_2, \lfloor \omega_1, \omega_3\rfloor \rfloor$$
\item graded Leibnitz identity $$\lfloor \omega_1, \omega_2\star\omega_3\rfloor =\lfloor \omega_1, \omega_2\rfloor \star \omega_3+(-1)^{(|\omega_1|+\frak m+1)|\omega_2|}\omega_2
\star\lfloor \omega_1, \omega_3\rfloor.$$
\end{enumerate}
\end{thm}

\begin{thm}\label{main4}
For even $\frak n$ we have the following identities 
\begin{enumerate}
\item graded commutativity $$\lfloor \omega_1, \omega_2\rfloor=(-1)^{(|\omega_1|+\frak m)(|\omega_2|+\frak m)}\lfloor \omega_2, \omega_1\rfloor$$
\item graded Leibnitz identity $$ \lfloor \omega_1, \omega_2 \star \omega_3\rfloor=\lfloor \omega_1, \omega_2\rfloor\star \omega_3+(-1)^{(\frak m+|\omega_1|)|\omega_2|}\omega_@\star\lfloor \omega_1, \omega_3\rfloor$$
\item an identity relating the $\lfloor \cdot, \cdot \rfloor$ and $\{\cdot, \cdot \}$ operations

\begin{equation}
\begin{split}
(-1)^{|\omega_1||\omega_3|+\frak m|\omega_2|}\Bigl (  \{\lfloor\omega_1, \omega_2\rfloor, \omega_3\}+    \lfloor\{\omega_1, \omega_2\}, \omega_3\rfloor    \Bigr)+\\
(-1)^{|\omega_1||\omega_2|+\frak m|\omega_3|}\Bigl ( \{\lfloor \omega_2, \omega_3\rfloor, \omega_1\}+\lfloor \{\omega_2, \omega_3\},\omega_1\rfloor\Bigr)+\\(-1)^{|\omega_2||\omega_3|+\frak m |\omega_1|}\Bigl( \{\lfloor \omega_3, \omega_1\rfloor, \omega_2\}+\lfloor\{ \omega_3, \omega_1\}, \omega_2\rfloor\Bigr)=0
\end{split}
\end{equation}

\end{enumerate}
\end{thm}

The proofs of Theorems~\ref{main3} and ~\ref{main4} are similar to the proofs of Theorem~\ref{main1}.

\m

{\bf Acknowledgments.}
I am grateful to Moira Chas, Ralph Cohen, Alexander Dranishnikov, Tobias Ekholm, Michael Penakava, Dennis Sullivan and David Webb for the stimulating
discussions. 
I am extremely thankful to Yuli Rudyak for all his help and the very inspiring discussions.

I am thankful to Michael Penkava who suggested that identities in the algebra would become much nicer if we formulate them for another operation $\{\cdot, \cdot\}$ that is closely related to $[\cdot, \cdot].$

I am thankful to the anonymous referees for the many valuable suggestions on how to improve the exposition in the previous versions of this paper.

This work was partially supported by a grant from the Simons Foundation (\# 513272  to Vladimir 
Chernov).

\end{document}